%% file: normalreconstruction.tex
\documentclass[a4paper,10pt,fleqn]{article}
\usepackage[english]{babel}%

\usepackage[numbers]{natbib}%
\usepackage[bottom=3cm,top=2cm,left=1.5cm,right=1.5cm]{geometry}%
\usepackage{framed,multirow,url,hyperref,setspace,graphicx,bm,xcolor,siunitx,nicefrac,multicol}
\usepackage{amsmath,amssymb,amsthm}
\usepackage{stmaryrd,colortbl,multido}

\usepackage[hang,centerlast]{subfigure}

\usepackage[draft,author={Johannes},color=red,hoffset=-1cm]{pdfcomment}%

\usepackage{listings}%

\usepackage{notation}%
\usepackage{mathoperators}%
\usepackage{reconstructnotation}%

\input{normalreconstruction_definitions.tex}%

\definecolor{newcolor}{rgb}{.8,.349,.1}

\def\worktitle{Second-order accurate normal reconstruction from volume fractions on unstructured meshes with arbitrary polyhedral cells}%

\providecommand{\keywords}[1]{\textbf{\textit{Keywords---}} #1}

\allowdisplaybreaks

\begin{document}

\newcommand{\refapp}[1]{appendix~\ref{app:#1}}%
\let\paragraphold\paragraph%
\renewcommand{\paragraph}[1]{\paragraphold{#1.}}%

\title{\worktitle}%
\author{Johannes Kromer\href{https://orcid.org/0000-0002-6147-0159}{\includegraphics[height=10pt]{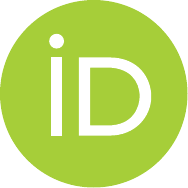}}, Fabio Leotta and Dieter Bothe\textsuperscript{$\dagger$}\href{https://orcid.org/0000-0003-1691-8257}{\includegraphics[height=10pt]{orcid_logo.pdf}}}%
\date{}%
\maketitle
\begin{center}
Mathematical Modeling and Analysis, Technische Universit\"at Darmstadt\\ Alarich-Weiss-Strasse 10, 64287 Darmstadt, Germany\\%
\textsuperscript{$\dagger$}Email for correspondence: \href{mailto:bothe@mma.tu-darmstadt.de?subject=Second-order\%20accurate\%20normal\%20reconstruction\%20from\%20volume\%20fractions\%20on\%20unstructured\%20meshes\%20with\%20arbitrary\%20polyhedral\%20cells}{bothe@mma.tu-darmstadt.de}
\end{center}

\begin{abstract}%
This paper introduces a novel method for the efficient second-order accurate computation of normal fields from volume fractions on unstructured polyhedral meshes. %
Locally, i.e.~in each mesh cell, an averaged normal is reconstructed by fitting a plane in a least squares sense to the volume fraction data of neighboring cells while implicitly accounting for volume conservation in the cell at hand. The resulting minimization problem is solved approximately by employing a \textsc{Newton}-type method. Moreover, applying the \textsc{Reynolds} transport theorem allows to assess the regularity of the derivatives. %
Since the divergence theorem implies that the volume fraction can be cast as a sum of face-based quantities, our method considerably simplifies the numerical procedure for applications in three spatial dimensions while demonstrating an inherent ability to robustly deal with unstructured meshes. %
We discuss the theoretical foundations, regularity and appropriate error measures, along with the details of the numerical algorithm. Finally, numerical results for convex and non-convex hypersurfaces embedded in cuboidal and tetrahedral meshes are presented, where we obtain second-order convergence for the normal alignment and symmetric volume difference. %
Moreover, the findings are substantiated by completely new deep insights into the minimization procedure. %
\end{abstract}

%
%
\input{01_introduction}%
%
%
\input{02_literature_review}%
%
%
\input{03_mathematical_details}%
%
%
\input{04_design_experiments}%
%
%
\input{05_numerical_results}%
%
%
\input{06_summary}%
%
%
\clearpage
\renewcommand{\bibname}{References}

\input{normalreconstruction.bbl}\begin{center}%
\textsc{Acknowledgment}\\[2ex]%
The authors gratefully acknowledge financial support provided by the German Research Foundation (DFG) within the scope of \href{www.sfbtrr75.de}{SFB-TRR 75 (project number 84292822)}. %

The figures in this manuscript were produced using the versatile and powerful library \href{https://ctan.org/topic/pstricks}{\texttt{pstricks}}. For further details and a collection of examples, the reader is referred to the book of \citet{pstricks_2008}.\\[12pt]%
\href{https://www.elsevier.com/authors/policies-and-guidelines/credit-author-statement}{\textsc{CRediT statement}}\\[2ex]%
\textbf{Johannes Kromer}: conceptualization, methodology, software, validation, investigation, data curation, visualisation, writing--original draft preparation, writing--reviewing and editing, supervision\\ %
\textbf{Fabio Leotta}: conceptualization, methodology, investigation, validation,  writing--reviewing and editing\\ %
\textbf{Dieter Bothe}: conceptualization, methodology, investigation, writing--reviewing and editing, supervision, funding acquisition, project administration%
\end{center}%
%
%
\input{99_appendix}
\end{document}

%% file: normalreconstruction_definitions.tex
\newcommand{\reffig}[1]{figure~\ref{fig:#1}}%
\newcommand{\reffigs}[1]{figures~\ref{fig:#1}}%
\newcommand{\reffigno}[1]{\ref{fig:#1}}%
\newcommand{\refFig}[1]{Figure~\ref{fig:#1}}%
\newcommand{\reftab}[1]{table~\ref{tab:#1}}%
\newcommand{\refEqn}[1]{Equation~\ref{eqn:#1}}%
\newcommand{\refeqs}[1]{eqs.~\ref{eqn:#1}}%
\newcommand{\refeqno}[1]{\ref{eqn:#1}}%
\makeatletter
\newcommand{\refeqn}{\@ifstar{\@ifstar\@@refeqn\@refeqn}{\@@@refeqn}}%
\def\@refeqn#1{\ref{eqn:#1}}%
\def\@@refeqn#1{eqs.~\ref{eqn:#1}}%
\def\@@@refeqn#1{eq.~\ref{eqn:#1}}%
\makeatother

\colorlet{shadecolor}{black!12}%
\newtheorem{notetheorem}{Note}[section]
\newenvironment{note}%
{\colorlet{shadecolor}{black!12}\begin{shaded}\begin{notetheorem}}%
{\end{notetheorem}\end{shaded}}%
\newcommand{\refnote}[1]{note~\ref{note:#1}}%

\newtheorem{remarktheorem}{Remark}[section]%
\newenvironment{remark}%
{\colorlet{shadecolor}{black!12}\begin{shaded}\begin{remarktheorem}}%
{\end{remarktheorem}\end{shaded}}%

\newtheorem{rationaletheorem}{Rationale}[section]%
{\colorlet{shadecolor}{black!12}\begin{shaded}\begin{rationaletheorem}}%
{\end{rationaletheorem}\end{shaded}}%

\newtheorem*{impltheo}{Implementation detail}
{\colorlet{shadecolor}{black!12}\begin{shaded}\begin{impltheo}}%
{\end{impltheo}\end{shaded}}%

%% file: 01_introduction.tex
\section{Introduction}\label{sec:introduction}%
The two-phase \textsc{Navier-Stokes} equations govern flows of immiscible fluid phases $\domain^-\foft$ and $\domain^+\foft$, separated by a sharp interface $\iface\foft$. For their numerical integration, especially in the presence of strong deformations and topological changes of the fluid phases, the Volume-of-Fluid (VOF) method of \citet{JCP_1981_vofm} has been successfully employed for decades. %
In a one-field formulation, the volume fraction $\polyvof\foftx$ encodes the assignment of a spatial point $\vx$ to the time-dependent fluid phases $\domain^-$ (\textit{interior}) and $\domain^+$ (\textit{exterior}) and, say, admits a value of one (zero) in the interior (exterior) phase. %
The evolution of the phases over time can be obtained from solving a passive advection equation of the volume fraction field, where the \textbf{P}iecewise \textbf{L}inear \textbf{I}nterface \textbf{C}alculation (PLIC) method of \citet{JCP_1998_rvt} can be considered one of the most common numerical schemes. %
For this \textit{geometric} transport, the approximative reconstruction of the interface from volume fractions constitutes a central element, in terms of both accuracy of the method and consumption of computational ressources. %
However, on a discrete computational mesh, the discontinuous character of the volume fraction makes it notoriously difficult to obtain accurate approximations of the normal field, since the majority of the available algorithms requires the computation of spatial derivatives in one form or another. %
While the disadvantageous influence of the discontinuity can be alleviated by smoothing the volume fractions, e.g., by applying a convolution with an appropriate kernel~\cite{FDI_1999_aaco}, %
doing so gives rise to numerical artefacts of its own \cite{IJMF_2002_dfdu,EJM_2002_anvo}. %
Furthermore, the accuracy of the geometric transport crucially depends on an accurate representation of the interface, i.e.~the normal field. In other words, the geometric transport is prone to amplify noise, whose magnitude and structure may be affected, among others, by the numerical transport scheme (\textit{split} \cite{JNA_1968_otca,JCP_2007_irwl} or \textit{un-split} \cite{JCP_2014_acff,JCP_2018_aeuf}; cf.~the review of \citet{JCP_2020_uusg}), the underlying mesh and the quality of the initial volume fractions. %
These mutual dependencies highlight the necessity for accurate interface reconstruction. %
The desired applicability to unstructured meshes with arbitrary polyhedral cells requires to restrain from computing finite differences of volume fractions and, hence, suggests to formulate an integration-based minimization approach. A prominent representative of this class, namely the \textbf{L}east Squares \textbf{V}OF \textbf{I}nterface \textbf{R}econstruction \textbf{A}lgorithm (LVIRA), was proposed by \citet{ISCFD_1991_avof} for structured square meshes in two spatial dimensions. In essence, the goal of the present paper is to extend their work to meshes composed of \textit{arbitrary polyhedra} in three spatial dimensions. %
In a qualitative sense, we may formulate the problem as follows: %
%
\begin{note}[Conceptual problem formulation]\label{note:conceptual_problem_formulation}%
For a given neighborhood of cells find a plane with normal $\plicnormal$ that minimizes the quadratic deviation of the induced and prescribed volume fractions, while matching the volume fraction in the center exactly. %
\end{note}%
Before subsection~\ref{subsec:strategy} provides a formal definition, the upcoming subsection~\ref{subsec:notation} introduces the relevant notation. %

%
%
\subsection{Notation}\label{subsec:notation}%
Throughout this paper, $\abs{\domain}$ denotes the \textsc{Lebesgue} measure (length, area and volume in $\setR$, $\setR^2$ and $\setR^3$) and number of elements of a set $\domain$, respectively. E.g., a cuboid with lower left corner $[-1,0,7]\transpose$ and upper right corner $[1,3,8]\transpose$ admits a volume of $\abs{[-1,1]\times[0,3]\times[7,8]}=6$. For finite sets, $\abs{M}$ denotes the number of elements of $M$, e.g.~$\abs{\set{1,3,5,6}}=4$. %
%
%
\paragraph{Meshes}
Let $\mathcal{M}_\domain=\set{\vfcell_k}_{k=1}^{N_\domain}$ be a mesh composed of $N_\domain$ polyhedral cells $\vfcell_k$, which form a disjoint decomposition of some finite domain $\domain\subset\setR^3$, i.e.\ $\domain=\bigcup_k\vfcell_k$ with $\vfcell_i\cap\vfcell_j=\emptyset$ for all $i\neq j$. %
%
%
In terms of data representation, a mesh $\mathcal{M}_\domain$ consists of a list of vertices\footnote{The meshes under consideration here are assumed to have no hanging nodes, i.e.~each vertex belongs to at least two faces.} $\set{\vx^\domain_i}$ as well as a list of planar polygonal faces $\set{\vfface_f}$. %
Each face $\vfface_f$ is in turn described by a list of $N^{\vfface}_f$ vertex labels $\set{i_{f,m}}_{m=1}^{N^{\vfface}_f}$, which are ordered counter-clockwise with respect to the unit normal $\vn^\vfface_f$. %
\begin{remark}[Labels]%
The identification of entities (e.g.~vertices or faces) associated to the mesh resorts to indices, which in practice are often referred to as \textit{labels}. Throughout this manuscript, these terms will be used synonymously. %
\end{remark}
Introducing $\vx^\vfface_{f,m}\defeq\vx^\domain_{i_{f,m}}$ (\textit{mesh to face}) for ease of notation, the $m$-th edge $\vfedge_{f,m}$ on face $\vfface_f$ is spanned by the line segment $\operatorname{conv}\fof{\vx^\vfface_{f,m},\vx^\vfface_{f,m+1}}$ with the co-normal %
$$\vN_{f,m}=\frac{\cross{\brackets{\vx^\vfface_{f,m+1}-\vx^\vfface_{f,m}}}{\vn^\vfface_f}}{\norm{\vx^\vfface_{f,m+1}-\vx^\vfface_{f,m}}_2},$$%
where the vertex labels are continued periodically, i.e.\ $\vx^\vfface_{f,N^\vfface_f+1}=\vx^\vfface_{f,1}$. 
\refFig{notation_illustration} illustrates the setup and the relevant quantities. %
\begin{figure}[htbp]
\null\hfill%
\includegraphics[page=1]{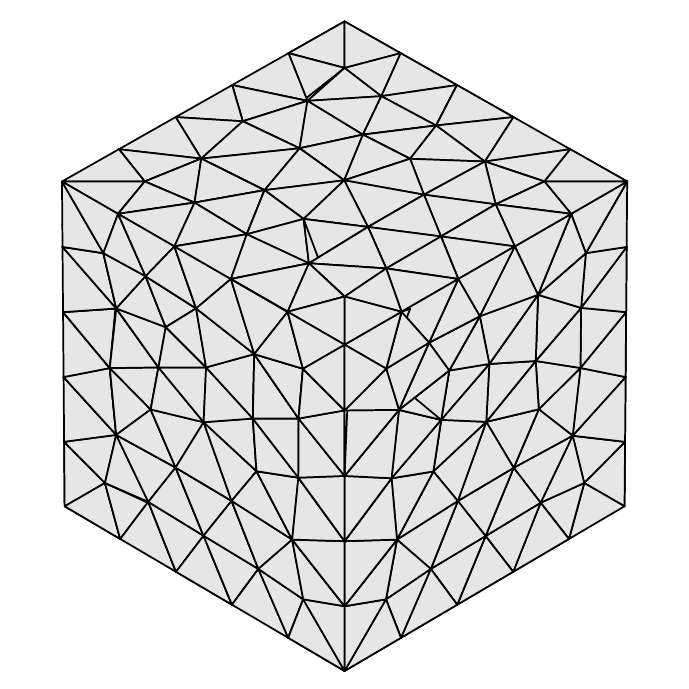}%
\hfill%
\includegraphics[page=2]{mesh_illustration}%
\hfill\null%
\caption{Illustration of a tetrahedron mesh (left) and single cell (right) with some of the relevant quantities.}%
\label{fig:notation_illustration}%
\end{figure}

Each cell $\vfcell_k$ is described by its boundary in terms of a list $\set{(f_{k,l},\omega_{k,l})}_{l=1}^{N^\vfface_k}$ of face labels $f_{k,l}$ and orientation flags\footnote{Since at most two cells share a face, the orientation flag sometimes is referred to as \textit{ownership}, because the orientation is chosen to coincide with the outer normal of the cell with, say, the smaller label, the so-called \textit{owner}.} %
$\omega_{k,l}\in\set{\pm1}$, where the latter ensure that the product $\omega_{k,l}\vn^\vfface_{f_{k,l}}$ is an \textit{outer} normal of $\vfcell_k$. %
To ensure applicability of the \textsc{Gaussian} divergence theorem, both the boundary of the faces $\partial\vfface_f$ and the cells $\partial\vfcell_k=\bigcup_l{\vfface_{f_{k,l}}}$ are assumed to admit no self-intersections. This is not a relevant restriction since objects of the latter class have no relevance for the desired application within finite-volume based methods. %
%
%
\paragraph{Neighborhood}\label{subsec:notation_neighborhod}%
We assign a neighborhood $\mathcal{N}_k$ (synonymously referred to as \textit{stencil}) to each cell $\vfcell_k$ (referred to as its \textit{center}), where the present work resorts to two different types: %
\begin{description}
\item[face-based:] The neighborhood contains all cells that share a face with the center cell, i.e. %
$$\neighborhood{k}^{\mathrm{face}}\defeq\set*{\vfcell_u\in\mathcal{M}_\domain:\exists\vfface_z\in\set{\vfface_f}\,s.t.\,\vfcell_u\cap\vfcell_k=\vfface_z}.$$%
The number of faces of the center cell resembles an upper bound for the neighborhood size, i.e.\ $\abs{\mathcal{N}_k^{\mathrm{face}}}\leq N^{\vfface}_k$. %
\item[edge-based:] The neighborhood contains all cells $\vfcell_u$ which contain at least one face $\vfface_m$ that shares at least one edge with one of the faces $\vfface_n$ of the stencil center, i.e.%
$$\neighborhood{k}^{\mathrm{edge}}\defeq%
\set*{\vfcell_u\in\mathcal{M}:\exists\brackets{\vfface_m\subset\vfcell_u,\vfface_n\subset\vfcell_k}\,s.t.\,\abs{\set{i_{m,l}}\cap\set{i_{n,l}}}\geq2}\supset\mathcal{N}_k^{\mathrm{face}}.$$%
Note that the above definition could be written less extensive resorting to edges instead of intersections of faces. However, we restrain from doing so since edges are typically not used in the discretization schemes and hence not explicitly defined as entities. %
\end{description}
%
As we shall see in section~\ref{sec:numerical_results} below, the size of the neighborhood significantly affects the reconstruction quality. Obviously, one could employ other types of neighborhood, e.g., vertex-based or iterated ones (neighbors of neighbors). However, there are two points to be kept in mind: %
\begin{enumerate}
\item Extending the neighborhood beyond the requirements of other components of the embedding numerical scheme (gradient/flux computation) implies additional costs, in terms of either consumption of memory for static or computational resources for dynamic meshes. Especially, parallel applications would experience a severe loss in performance. %
\item In a geometrical sense, the normal $\nS$ of the hypersurface $\iface$ is a local quantity associated to a point on $\iface$. The discrete representation of the bulk phases by volume fractions according to %
$$\polyvof\fof{\vy}=\frac{1}{\abs{\ball[R_\polyvof]{\vy}}}\int_{\ball[R_\polyvof]{\vy}}{\heaviside{-\lvlset_\iface\fofx}\dvol},$$%
corresponds to a spatial averaging, where $\ball[R]{\vy}$ and $\lvlset_\iface$, respectively, denote a domain with characteristic length $R$ centered around $\vy$ and the level-set such that $\iface=\set{\vx\in\setR^3:\lvlset_\iface\fofx=0}$. %
The associated discrete normal reads 
$$\hat{\vn}_\iface\fof{\vy}=\frac{1}{\abs{\ball[R_\vn]{\vy}\cap\iface}}\int_{\ball[R_\vn]{\vy}\cap\iface}{\nS\darea}.$$ %
While a larger neighborhood for the normal reconstruction ($R_\vn>R_\polyvof$) may be advantageous in terms of numerical robustness, beyond the inconsistency there is also a smoothing effect, which may become disadvantageous for capillary flows \cite{IJMF_2002_dfdu,EJM_2002_anvo}. %
\end{enumerate}
%
%
\paragraph{Summation}%
Throughout this work, the subscripts $k$, $f$ and $m$ always refer to cells, faces and edges/vertices, respectively. Wherever unambiguously possible, we formally abbreviate the labels by the summation indices, i.e.~we omit the sets over which the summation index runs.. E.g., $\sum_{f\in\set{f_{k,l}}}$ (all faces associated to cell $\vfcell_k$) becomes $\sum_{f}$ and $\sum_{m\in\set{i_{f,m}}}$ (all edges associated to face $\vfface_f$) becomes $\sum_m$. %
Since the proposed algorithm traverses the mesh $\mathcal{M}_\domain$ cell by cell, the remainder of this manuscript considers a single cell %
$\vfcell_k$ (subsection~\ref{subsec:computation_volume_fractions}) %
and neighborhood $\neighborhood{k}$ (subsection~\ref{subsec:minimization_strategy}), respectively, where we omit the label for ease of notation (i.e., $\vfcell_k\rightarrow\vfcell$ and $\neighborhood{k}\rightarrow\neighborhood{}$). %
With $\mathcal{I}=\set{u\in\setN:\vfcell_u\in\neighborhood{}}$ denoting the list of cell labels associated to the neighborhood $\neighborhood{}$, the sum $\sum_{k\in\mathcal{I}}$ (all cells in the neighborhood $\neighborhood{}$) becomes $\sum_k$, where $k=0$ corresponds to the center. %
\begin{remark}[Relevance for applications]%
In many common flowsolvers, e.g., OpenFOAM, each face $\vfface_f$ of the mesh $\mathcal{M}_\domain$ has assigned to it the labels of its owner and neighbor. For the divergence-based evaluation of averages, the algorithm traverses the list of faces, whose $\omega$-weighted contributions are added to the value associated to the respectively owning ($\omega=1$) and neighboring ($\omega=-1$) cell. Hence, schemes of this type get along without an explicit assignment of faces to cells. For other applications, such as, e.g., the initialization of volume fractions \cite{JCP_2021_toai}, the contributions of faces $\vfface_f$ to the owning and neighboring cell do not coincide in general, necessitating an assignment of faces to cells. %
\end{remark}%

%
%
\paragraph{Plane parametrization}%
For mathematical convenience, the sought plane $\plicplane$ will be parametrized by the signed distance $\signdist$ to some arbitrary but spatially fixed reference $\xbase$ and its normal $\plicnormal$ given in spherical angles, i.e. %
\begin{align}%
\plicplane\fof{s,\varphi,\theta}=\set{\vx\in\setR^3:\pliclvlset{}\fof{\vx}=0}%
\quad\text{with}\quad%
\pliclvlset{}\fof{\vx}=\iprod{\vx-\xbase}{\plicnormal}-\signdist%
\quad\text{and}\quad%
\plicnormal=\brackets*[s]{\begin{matrix}\cos\varphi\sin\theta\\\sin\varphi\sin\theta\\\cos\theta\end{matrix}},
\label{eqn:plane_parametrization}%
\end{align}
where the parameter domain $\unitsphere=[0,2\pi)\times[0,\pi]$ will henceforth be referred to as \textit{unit sphere}. %
Accordingly, the derivatives of the normal are to be interpreted with respect to $\varphi$ and $\theta$. %
The negative half-space associated to $\plicplane$ reads $\neghalfspace{\signdist,\plicnormal}=\set{\vx\in\setR^3:\iprod{\vx-\xbase}{\plicnormal}\leq\signdist}$ and induces the volume fraction %
\begin{align}
\polyvof\fof{\signdist,\plicnormal,\vfcell_k}=\polyvof_k\fof{\signdist,\plicnormal}=\frac{\abs{\vfcell_k\cap\neghalfspace{\signdist,\plicnormal}}}{\abs{\vfcell_k}}.%
\label{eqn:volume_fraction_halfspace}%
\end{align}
\refFig{plane_parametrization} provides an illustration in two spatial dimensions. %
\begin{figure}[htbp]%
\null\hfill%
\includegraphics[page=1]{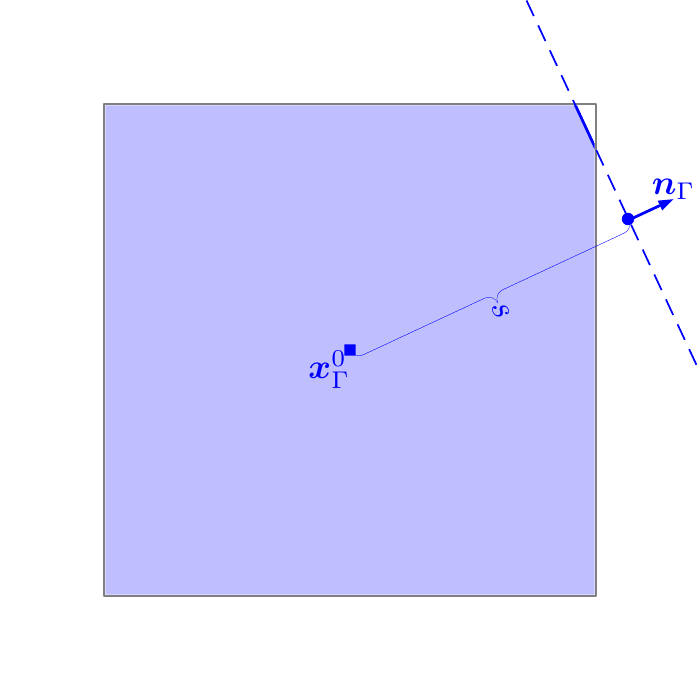}%
\hfill%
\includegraphics[page=2]{plane_parametrization}%
\hfill\null%
\caption{Cell $\vfcell_k$ intersected by plane $\plicplane$ with normal $\plicnormal$, reference point $\xbase$ and signed distance $\signdist$ for two different orientations. Note that base point of $\plicplane$ (\textcolor{blue}{$\bullet$}/\textcolor{cyan}{$\bullet$}) produced by the parametrization does not coincide with the PLIC centroid in general. In fact, the base point is not necessarily contained in $\plicplane\cap\vfcell_k$ (right panel).}%
\label{fig:plane_parametrization}
\end{figure}%
\paragraph{Averages}%
For an integrable function $\phi:\setR^3\mapsto\setR^n$ and a non-empty finite domain $\domain\subset\setR^3$, we denote the average as %
\begin{align}
\average[\domain]{\phi}\defeq\frac{1}{\abs{\domain}}\int_{\domain}{\phi\fof{\vx}\dvol}.\label{eqn:average_definition_continuous}%
\end{align}
\paragraph{Topological properties of geometric entities}
The assignment of an intersection status to a face $\vfface_f$ requires a hierarchically consistent evaluation of the topological properties of its vertices $\set{\vx^\vfface_{f,m}}$ and edges $\set{\vfedge_{f,m}}$. With respect to the orientable hypersurface $\plicplane$, described by the level-set function given in \refeqn{plane_parametrization}, the logical status $\vfstatus{}$ of a vertex is either interior ($\vfstatus{}=-1$), intersected ($\vfstatus{}=0$) or exterior ($\vfstatus{}=1$). %
\begin{remark}[Robustness]%
For the purpose of numerical robustness, the status assignment employs a tubular neighborhood of thickness $2\zerotol$ around $\plicplane$, corresponding to the interval $(-\zerotol,\zerotol)$, i.e.%
\begin{align}%
\vfstatus{\vx}\defeq%
\begin{cases}%
\phantom{-}0&\text{if}\quad\abs{\pliclvlset{}\fof{\vx}}<\zerotolerance,\\%
\sign{\pliclvlset{}\fof{\vx}}&\text{if}\quad\abs{\pliclvlset{}\fof{\vx}}\geq\zerotolerance.%
\end{cases}%
\end{align}
In other words, any point whose absolute distance to $\plicplane$ falls below $\zerotol$ is considered to be on $\plicplane$. The choice of an appropriate tolerance strongly depends on the absolute value of the characteristic length scale $h$ of the mesh $\mathcal{M}_\domain$, especially for $h\ll1$. Throughout this work, we have $h\approx1$ and let $\zerotolerance\defeq\num{e-14}$. In fact, all zero-comparisons are implemented in this way. %
\end{remark}%
The hierarchically superior entity is an edge $\vfedge$, whose status is a function of the status of its associated vertices. %
\begin{table}[htbp]
\centering%
\caption{Edge status as function of the status of the associated vertices.}%
\label{tab:intersection_status_edge}%
\renewcommand{\arraystretch}{1.5}%
\begin{tabular}{c||ccc|cc|cc|c}
&\multicolumn{3}{c|}{}&\multicolumn{5}{c}{\textbf{degenerate}}\\%
&\rotatebox{90}{exterior}&\rotatebox{90}{interior}&\rotatebox{90}{intersected}&\multicolumn{2}{c|}{\rotatebox{90}{exterior}}&\multicolumn{2}{c|}{\rotatebox{90}{interior}}&\rotatebox{90}{intersected}\\%
&\includegraphics[page=4]{edgestatus_illustration}&\includegraphics[page=2]{edgestatus_illustration}&\includegraphics[page=1]{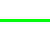}&\multicolumn{2}{c|}{\includegraphics[page=5]{edgestatus_illustration}}&\multicolumn{2}{c|}{\includegraphics[page=3]{edgestatus_illustration}}&\includegraphics[page=6]{edgestatus_illustration}\\%
\hline%
$\vfstatus{\vx^\vfface_{k,m}}$&$1$&$-1$&$\pm1$&$1$&$0$&$-1$&$\phantom{-}0$&$0$\\%
$\vfstatus{\vx^\vfface_{k,m+1}}$&$1$&$-1$&$\mp1$&$0$&$1$&$\phantom{-}0$&$-1$&$0$\\%
\hline%
$\vfstatus{\vfedge_{k,m}}$&$1$&$-1$&$0$&\multicolumn{2}{c|}{$2$}&\multicolumn{2}{c|}{$-2$}&$3$%
\end{tabular}%
\begin{tabular}{c}
\includegraphics[page=7]{edgestatus_illustration}%
\end{tabular}
\end{table}%
For a neighborhood, the intersection status gathers the status of all cells, cf.~\reffig{intersection_status}. %
\begin{figure}[htbp]%
\null\hfill%
\includegraphics[page=1]{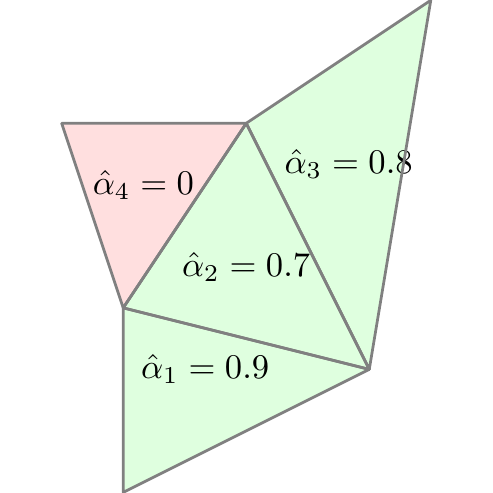}%
\hfill%
\includegraphics[page=2]{intersection_status_illustration}%
\hfill%
\includegraphics[page=3]{intersection_status_illustration}%
\hfill\null%
\caption{Intersection status of a neighborhood induced by data ($\polyvofdata_k$, left) and by an intersecting plane $\plicplane$ ($\polyvof_k$) with status match (center) and status mismatch (right), respectively. The meaning of the colors is analogous to \reftab{intersection_status_edge}.}%
\label{fig:intersection_status}%
\end{figure}
%
%
\subsection{Problem formulation, strategy \& objectives}\label{subsec:strategy}%
\begin{note}[Mathematical problem formulation]%
Let $\neighborhood{}=\set{\vfcell_k}$ be a neighborhood of $N$ cells with associated volume fractions $\set{\polyvofdata_k}$. We want to find a plane $\plicplane^\ast$, defined by the normal $\plicnormalref$ (corresponding to $\varphi^\ast$ and $\theta^\ast$) and an associated signed distance $\signdistref$, that minimizes the quadratic deviation of the induced $\set{\polyvof_k}$ and prescribed $\set{\polyvofdata_k}$ volume fractions, while the deviation in the center cell (label $k=0$) is zero. %
Employing the parametrization of \refeqn{plane_parametrization}, one obtains a minimization problem in $\fof{\signdist,\varphi,\theta}\in\setR^3$ with the error functional %
\begin{align}
\error\fof{\signdist,\varphi,\theta}\defeq%
\frac{1}{2}\sum\limits_{k=1}^{N}{\mu_k\brackets*{\polyvof\fof{\vfcell_k;\signdist,\plicnormal\fof{\varphi,\theta}}-\polyvofdata_k}^2}%
\quad\text{with weights}\quad%
\mu_k\in\setR^+,%
\label{eqn:error_functional}%
\end{align}
which is subject to the scalar constraint for the volume conservation %
\begin{align}
\centerconstraint{\signdistref;\plicnormalref}\defeq\polyvof\fof{\vfcell_0;\signdistref,\plicnormalref}-\polyvofdata_0\stackrel{!}{=}0.%
\tag{C}%
\label{eqn:center_constraint}%
\end{align}
In what follows, we omit the asterisks for ease of notation. %
\end{note}
\paragraph{Strategy}%
By exploiting the properties of the constraint in \refeqn{center_constraint}, we can transform the constrained optimization problem into an unconstrained one with intrinsic volume conservation. %
The following considerations convey the basic idea: %
for a given normal $\plicnormal$, the volume fraction $\polyvof_k$ associated to a cell $\vfcell_k$ is a piecewise cubic polynomial in the signed distance $\signdist$, whose coefficients are piecewise smooth functions of $\plicnormal$, i.e.~of $\varphi$ and $\theta$. %
Under mild restrictions, the signed distance $\signdistref$ that induces the prescribed volume $\polyvofdata_0$ in the center cell $\vfcell_0$ (corresponding to the \textit{unique} solution of \refeqn{center_constraint}), is itself a continuously differentiable function of $\varphi$ and $\theta$. %
Thus, replacing the variable $\signdist$ (signed distance) in the parametrization of the plane $\plicplane$ by $\signdistref\fof{\varphi,\theta}$ reduces the problem dimension by one, while simultaneously enforcing volume conservation in the center cell. %
Due to the periodicity of the parametrization in spherical angles, we obtain a piecewise smooth unconstrained minimization problem for $\brackets{\varphi,\theta}\in\setR^2$. %
As we will argue in subsection~\ref{subsec:regularity_considerations} below, the resulting error functional is continuously differentiable under mild restrictions on the mesh geometry. Furthermore, it is twice continuously differentiable in all but finitely many points (modulo $2\pi$ and $\pi$, respectively, for $\varphi$ and $\theta$). For our purposes, the \textsc{Gauss-Newton} method hence constitutes an appropriate minimization scheme; cf.~\citet[Theorem~10.1]{numopt_2006}. %
The required derivatives of $\signdistref$ with respect to $\varphi$ and $\theta$ can be easily obtained from differentiating the cubic polynomial governing the volume conservation in \refeqn{center_constraint} and rearranging. %

Thus, building on the work of \citet{JCP_2021_fbip}, we develop an iterative algorithm to minimize \refeqn{error_functional} with an implicit treatment of the volume conservation constraint. %
The two conceptually disjoint tasks, namely %
\begin{enumerate}
\item the efficient computation of the volume fractions $\polyvof_k$ and their derivatives for arbitrary polyhedral cells $\vfcell_k$; %
\item the design of an appropriate \textsc{Newton}-type minimization scheme which implicitly accounts for the volume conservation constraint in \refeqn{center_constraint}%
\end{enumerate}
will be the subject of section~\ref{sec:mathematical_foundation}. %
Subsequently, in section~\ref{sec:numerical_results}, we %
\begin{enumerate}
\item assess the algorithm by conducting a series of numerical experiments for a variety of hypersurfaces on hexahedral/tetrahedral meshes and %
\item provide local in-depth information about the minimization process. %
\end{enumerate}%
\paragraph{Novelties}%
The present work introduces the first algorithm that is capable of reconstructing the interface from volume fractions at \textbf{second-order accuracy} on \textbf{unstructured meshes} with \textbf{arbitrary polyhedral cells} in three spatial dimensions with \textbf{intrinsic volume conservation}. %

%% file: 02_literature_review.tex
\section{Literature review}\label{sec:literature_review}%
There are several methods for the reconstruction of the approximate interface, resorting to different approaches. E.g., on structured \textsc{Cartesian} meshes, \citet{JCP_2009_aaas} represents the interface as (the graph of a) \textit{height function} and approximates the discrete heights as a weighted sum of the volume fractions along the coordinate axes. From the literature, it is well known that the accuracy of the height function degreades if the components of the normal admit approximately equal magnitude in two or three of the coordinate directions. 
%
%
 
A less complex alternative are the variants of
the algorithm proposed by \citet{MOC_1989_tsos}, which iteratively increases accuracy of the interface normal based on the interface centroids in the surrounding cells.
\citet{JCP_2014_ahpc} and \citet{JCP_2015_aina} propose a generalized height function method for normal and curvature reconstruction on unstructured meshes in two and three spatial dimensions, respectively, by either exploiting the reduced complexity of the problem in 2D or embedding a supplementary structured cartesian mesh in the unstructured mesh by an interpolation procedure regarding the volume fractions to then make use of the traditional height function method in 3D. %
%
%
In his seminal paper, \citet{Youngs1982} approximates the normal $\plicnormal$ from the discrete gradient of the volume fraction $\polyvof$, i.e.%
\begin{align}
\plicnormal\approx-\frac{\grad{\polyvof}}{\norm{\grad{\polyvof}}}.\label{eqn:youngs_normal_approximation}%
\end{align}
On a structured \textsc{Cartesian} mesh, the evaluation of $\grad{\polyvof}$ resorts to a finite difference scheme employing the up to 27 cells in a $3\times3\times3$ neighborhood. After computing the normal from \refeqn{youngs_normal_approximation}, the plane $\plicplane$ is shifted in normal direction such that its negative halfspace truncates from the center cell a sub-polyhedron whose volume matches a prescribed value. %
As has been shown by \citet{JCP_2004_soav}, this method is first-order accurate. %
%
%
\citet{ISCFD_1991_avof} introduces a \textbf{L}east Squares \textbf{V}OF \textbf{I}nterface \textbf{R}econstruction \textbf{A}lgorithm (LVIRA) to solve the two-dimensional version of the problem stated in \refnote{conceptual_problem_formulation}. \citet{ISCFD_1991_avof} parametrizes the line representing the reconstructed interface in terms of its slope, employs a central difference algorithm to obtain an initial guess and then uses \textsc{Brent}’s algorithm to minimize the error functional. The claim that "\textit{[...] this method reconstructs all linear interfaces exactly}" is somewhat misleading, since the orientation remains ambiguous; cf.~\reffig{3x3_stencil_error_illustration}. %
\citet{JCP_2004_soav} propose a modification with improved efficiency (\textit{efficient} LVIRA or ELVIRA for short) by approximating the continuous minimization problem by a reduction to a finite number of normal candiates, from which the algorithm selects the one with the smallest induced error. %
\citet{JCP_2019_aaes} extend the isoAdvector-library of \citet{RSOS_2016_acmf} by introducing a reconstruction scheme based on \textit{reconstructed distance functions}. For spheres on \textsc{Cartesian} and polyhedral meshes in three spatial dimensions, they report second-order convergence of both the interface normal in an $L_\infty$-type norm. %
\citet{JCP_2008_romm} introduce a Moment-of-Fluid (MoF) method, which takes into account the centroids of the phases in each cell. While methods of this type show, ceteris paribus, lower erros, the additional information required may induce considerable increase of computational costs. %
%
%
\begin{figure}[htbp]%
\null\hfill%
\includegraphics[page=1]{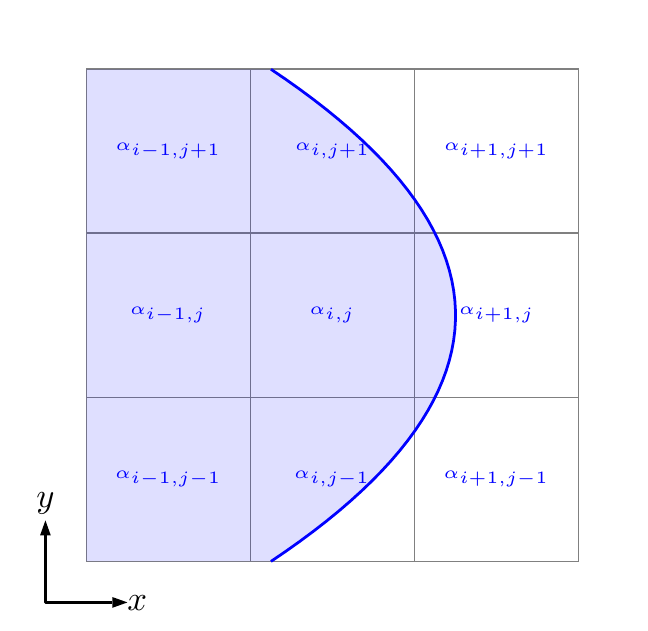}%
\hfill%
\includegraphics[page=3]{literature_review_campbell}%
\hfill\null%
\\%
\null\hfill%
\includegraphics[page=2]{literature_review_campbell}%
\hfill%
\includegraphics[page=4]{literature_review_campbell}%
\hfill\null%
\caption{Two-dimensional sketch of the spline interpolation scheme based on column-cumulated volume integrals in coordinate-aligned columns (here: $H_j=\sum_{i}{\polyvof_{i,j}\Delta x}$) proposed by \citet{JCP_2021_aaho}. The interface $\iface$ (blue) is required to be the graph over one of the coordinate planes, where the extension of the graph base, i.e.~the spatial resolution (number of nodes \textcolor{red}{$\bullet$}), determines the interpolation order.}%
\label{fig:campbell_literature}%
\end{figure}
%
%
For a comprehensive review of unstructured geometrical VOF-methods, the reader is referred to the recent review of \citet{JCP_2020_uusg}. %
\citet{JCP_2007_irwl} present two interface reconstruction procedures within a 3D standalone VOF scheme on cartesian meshes. One employs purely geometrical criterias in a static manner while the other reconstruction procedure resolves around the minimization of  a distance functional whose solution is obtained via inversion of a linear equation. Both methods heavily rely on the block structure of $n\times n\times n$ \textsc{Cartesian} stencils. %
\citet{JCP_2002_actd} aim to reduce computational effort by considering minimization of a constrained least squares problem over a discrete set of plane parameters that is obtained by applying ELVIRA rather than the continuous parameter superset. %
%
%
\citet{IJNMF_2008_anvo} propose a piecewise-planar interface reconstruction with cubic-\textsc{B\'ezier} fit. %
\citet{TVCG_2011_acog} give a comparison of gradient estimation methods on unstructured meshes. %
\citet{CF_2006_a3du} propose the CVTNA reconstruction scheme on orthogonal meshes that reportedly achieved high-order accuracy by considering a filtering procedure for outliers. %
%
%
\citet{JCP_2021_anie} introduce an iso-surface reconstruction algorithm on unstructured meshes that relies on owner information on the level of vertices. %
The spline-based reconstruction of \citet{JCP_2021_aaho} is described in \reffig{campbell_literature}. %

%% file: 03_mathematical_details.tex
%
%
\section{Mathematical foundations of the method}\label{sec:mathematical_foundation}%
As stated in subsection~\ref{subsec:strategy}, the purpose of the proposed algorithm consists in finding a normal $\plicnormal$ and signed distance $\signdist$ (corresponding to the base point $\xref=\xbase+\signdist\plicnormal$), such that the volume fractions $\set{\polyvof_k}$ induced by the plane $\plicplane$ (cf.~\refeqn{volume_fraction_halfspace}) admit a minimal quadratic deviation from the prescribed ones $\set{\polyvofdata_k}$ in a given neighborhood $\neighborhood{}$, while the deviation is zero in the center cell $\vfcell_0$. %
\begin{note}[Stencil data]%
In a numerical flow simulation scheme, one applies the reconstruction of the approximative interface $\plicplane$ only to intersected cells. Throughout this manuscript, the associated volume fractions induce the aforementioned status, say $\voftolerance<\polyvof<1-\voftolerance$ with some positive $\voftolerance\ll1$. For the present section the numerical value of $\voftolerance$ is of secondary relevance. However, let $\voftolerance\defeq\num{e-9}$. %
In what follows, we assume that the stencil cell $\vfcell_0$ is intersected (i.e.~$\polyvofdata_0\in[\voftolerance,1-\voftolerance]$), where it is crucial to note that this implies $\plicplane\cap\vfcell_0\neq\emptyset$, i.e.~the center cell contains a non-empty patch of the plane $\plicplane$. %
Beyond that, no assumptions concerning $\set{\polyvofdata_k}$ are made, which especially means that the neighborhood does not necessarily contain bulk cells ($\polyvofdata_k\in[0,\voftolerance]\cup[1-\voftolerance,1]$). %
\end{note}%
Therefore, this section will %
provide an efficient method for the computation of volume fractions and gradients with respect to the spherical angles $\varphi$ and $\theta$ (subsection~\ref{subsec:computation_volume_fractions}), %
discuss the regularity of the functions under consideration (subsection~\ref{subsec:regularity_considerations}), %
derive the implicit treatment of the volume constraint (subsection~\ref{subsec:intrinsic_volume_conservation}) and finally %
formulate the minimization problem along with a suitable solution strategy (subsection~\ref{subsec:minimization_strategy}). %
\input{03_01_computation_volume_fractions}%
%
%
\input{03_02_regularity}%
%
%
\input{03_03_intrinsic_volume_conservation}%
\input{03_04_minimization_strategy}%

%% file: 03_01_computation_volume_fractions.tex
\subsection{Efficient computation of volume fractions}\label{subsec:computation_volume_fractions}%
Within a PLIC interface positioning scheme, \citet{JCP_2021_fbip} introduce an efficient algorithm for the computation of the volume of a truncated arbitrary polyhedron, which we employ with a slight adaption that facilitates the computation of the derivatives. %
By combining a dynamic choice of the origin with the application of the \textsc{Gaussian} divergence theorem in appropriate form, the volume of a truncated polyhedron $\vfcell$ can be cast as %
\begin{align}
\polyvof\fof{\vfcell;\signdist,\plicnormal}=%
\frac{\abs{\vfcell\cap\neghalfspace{\signdist,\plicnormal}}}{\abs{\vfcell}}=%
\frac{1}{3\abs{\vfcell}}\sum\limits_{f}{\brackets*{\volcoeffconst{f}+\signdist\volcoefflin{f}\fof{\plicnormal}}\immersedarea{f}\fof{\signdist;\plicnormal}}\label{eqn:volume_computation_summation}%
\end{align}
with the coefficients %
\begin{align}
\volcoeffconst{f}=\iprod{\vx^\vfface_{f,1}-\xbase}{\vn^\vfface_f}%
\quad\text{and}\quad%
\volcoefflin{f}=-\iprod{\vn^\vfface_f}{\plicnormal},\label{eqn:volume_coefficients}%
\end{align}%
where $\iprod{\vx}{\vy}=\vx\transpose\vy$ denotes the inner product of two real vectors $\vx,\vy\in\setR^N$. %
\refEqn{volume_computation_summation} implies that the volume fraction can be represented as a sum of binary products of signed distances to an arbitrary but spatially fixed reference $\xbase$. 
%
%
For faces $\vfface_f$ that are not parallel to the plane $\plicplane$, i.e.~for $\iprod{\vn^\vfface_f}{\plicnormal}\not\in\set{-1,1}$, in an analogous manner one obtains the immersed area %
\begin{align}
\immersedarea{f}=\abs{\vfface_f\cap\neghalfspace{}}=\frac{1}{2}\sum\limits_{m}{\brackets*{\areacoeffconst{f}{m}\fof{\plicnormal}+\signdist\areacoefflin{f}{m}\fof{\plicnormal}}\immersedlength{f}{m}\fof{\signdist;\plicnormal}},%
\end{align}
where $\immersedlength{f}{m}=\abs{\vfedge_{f,m}\cap\neghalfspace{}}$ denotes the immersed length associated to edge $\vfedge_{f,m}$. The coefficients are %
\begin{align}
\begin{split}%
\areacoeffconst{f}{m}&=\iprod{\vx^\vfface_{f,m}-\xrefface{f}}{\vN_{f,m}}=\iprod{\vx^\vfface_{f,m}}{\vN_{f,m}}+%
\frac{\iprod*{\iprod{\vx^\vfface_{f,1}}{\vn^\vfface_f}\vn^\vfface_f-\xbase}{\plicnormal}\iprod{\vN_{f,m}}{\plicnormal}}{1-\iprod{\vn^\vfface_f}{\plicnormal}^2},\\%
\areacoefflin{f}{m}&=-\frac{\iprod{\vN_{f,m}}{\plicnormal}}{1-\iprod{\vn^\vfface_f}{\plicnormal}^2},%
\end{split}\label{eqn:face_coefficients}%
\end{align}%
where by assumption ($\iprod{\vn^\vfface_f}{\plicnormal}\not\in\set{-1,1}$) the intersection of the face $\vfface_f$ and the plane $\plicplane$ is non-empty and contains %
\begin{align}
\xrefface{f}\fof{\signdist,\plicnormal}=%
\frac{\iprod{\xbase}{\plicnormal}+\signdist-\iprod{\vx^\vfface_{f,1}}{\vn^\vfface_f}\iprod{\vn^\vfface_f}{\plicnormal}}{1-\iprod{\vn^\vfface_f}{\plicnormal}^2}\plicnormal+%
\frac{\iprod{\vx^\vfface_{f,1}}{\vn^\vfface_f}-%
\brackets{\iprod{\xbase}{\plicnormal}+\signdist}\iprod{\vn^\vfface_f}{\plicnormal}}{1-\iprod{\vn^\vfface_f}{\plicnormal}^2}\vn^\vfface_f.\label{eqn:face_origin_explicit}%
\end{align}
\begin{remark}[Efficiency]%
Note that the computation of the area of a planar polygon in $\setR^3$ actually poses a two-dimensional problem. If only the polygon area is required, a projection onto one of the coordiante planes can be exploited to reduce the computational effort \cite{JCP_2021_toai,JCP_2021_fbip}. However, computing the derivatives of the area with respect to $\varphi$ and $\theta$ from such a projection yields cumbersome expressions degrading the prior gain. %
\end{remark}%
Abbreviating $\frac{\partial}{\partial \signdist}$ by $\partial_\signdist$ for ease of notation ($\varphi$, $\theta$ analogous), the derivatives of the volume fraction $\polyvof$ read %
\begin{align}%
\partial_\signdist\polyvof&=%
\frac{1}{3\abs{\vfcell}}\sum\limits_{f}{%
\volcoefflin{f}\immersedarea{f}%
+%
\brackets*{\volcoeffconst{f}+\signdist\volcoefflin{f}}\partial_\signdist\immersedarea{f}%
},\label{eqn:dalphads}\\%
\partial_\varphi\polyvof&=%
\frac{1}{3\abs{\vfcell}}\sum\limits_{f}{%
\signdist\partial_\varphi\volcoefflin{f}\immersedarea{f}%
+%
\brackets*{\volcoeffconst{f}+\signdist\volcoefflin{f}}\partial_\varphi\immersedarea{f}%
}\label{eqn:dalphadphi},%
\end{align}
with the derivatives of the immersed area %
\begin{align}%
\partial_\signdist\immersedarea{f}&=%
\frac{1}{2}\sum\limits_{m}{%
\areacoefflin{f}{m}\immersedlength{f}{m}%
+%
\brackets*{\areacoeffconst{f}{m}+\signdist\areacoefflin{f}{m}}\partial_\signdist\immersedlength{f}{m}%
}%
,\\%
\partial_\varphi\immersedarea{f}&=%
\frac{1}{2}\sum\limits_{m}{%
\brackets*{\partial_\varphi\areacoeffconst{f}{m}+\signdist\partial_\varphi\areacoefflin{f}{m}}\immersedlength{f}{m}%
+%
\brackets*{\areacoeffconst{f}{m}+\signdist\areacoefflin{f}{m}}\partial_\varphi\immersedlength{f}{m}%
}\label{eqn:dAdphi}%
.%
\end{align}
Note that \refeqs{dalphadphi} and \refeqno{dAdphi} hold analogously for $\theta$ and \refapp{derivatives_vol_coeff} contains the derivatives of the coefficients. %
\begin{note}[Higher-order derivatives]\label{note:higher_order_derivatives}%
While the \textsc{Hessian} of the volume fraction could be easily derived, we omit doing so in anticipation of the strategy to be outlined in subsection~\ref{subsec:minimization_strategy}: we seek to employ a \textsc{Gauss-Newton}-type approach for the minimization of an error functional of least-squares type. For this class, the \textsc{Hessian} of the functional can be approximated by the exterior product of the gradient, implying that no second derivatives must be numerically cmputed. %
\end{note}%

%% file: 03_02_regularity.tex
%
%
\subsection{Regularity considerations}\label{subsec:regularity_considerations}%
From an implementation point of view, the strategy outlined in subsection~\ref{subsec:computation_volume_fractions} proves advantageous for the computation of the volume fractions and their derivatives. However, alternatively one may employ the \textsc{Reynolds} transport theorem to obtain %
\begin{align}
\partial_i\polyvof=\frac{1}{3\abs{\vfcell}}\int\limits_{\plicplane\cap\vfcell}{\iprod{\partial_i\vf_\plicplane}{\plicnormal}\darea}\label{eqn:gradient_vof_reynolds}%
\end{align}
with the function %
\begin{align}
\vf_\plicplane=\xbase+\signdist\plicnormal+\vtau_\plicplane\vt,%
\quad\text{where}\quad%
\vtau_\plicplane=\brackets*[s]{\begin{matrix}-\sin\varphi&\cos\varphi&0\\\cos\varphi\cos\theta&\sin\varphi\cos\theta&-\sin\theta&\end{matrix}}\transpose.%
\end{align}
%
Besides the concise formulation, \refeqn{gradient_vof_reynolds} offers a geometric interpretation of the derivatives of the volume fractions. E.g., it is well known in literature that $\partial_\signdist\polyvof_k=\abs{\plicplane\cap\vfcell_k}\abs{\vfcell_k}^{-1}$, i.e.~the derivative with respect to the signed distance corresponds to the area of the enclosed plane segment. %
In general, \refeqn{gradient_vof_reynolds} states that all derivatives result from integrating a continuous integrand over a continuous domain. Hence, we may deduce the following regularity properties: %
As long as $\plicplane$ does not contain a face of $\vfcell_k$, the continuity of $\grad{\polyvof}$ is guaranteed. This can be easily deduced after formulating the surface integral as an integral in 2D \textsc{Euclidean} space and using standard knowledge about parameter integrals. %
Analogously, one can show that $\hessian{\polyvof}$ is discontinuous for points at which $\plicplane$ contains an edge of $\vfcell_k$. %

%% file: 03_03_intrinsic_volume_conservation.tex
%
%
\subsection{Intrinsic volume conservation}\label{subsec:intrinsic_volume_conservation}%
\begin{figure}[htbp]
\null\hfill%
\includegraphics[page=1]{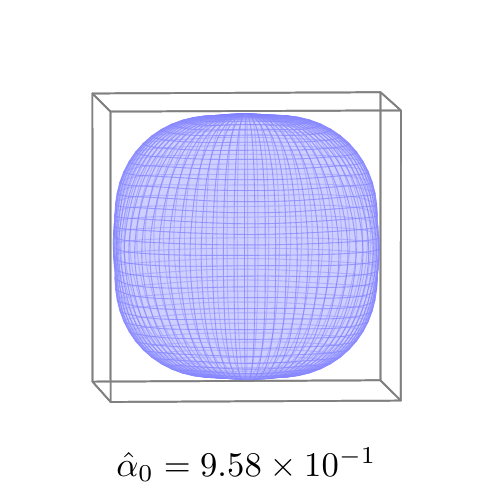}%
\hfill%
\includegraphics[page=1]{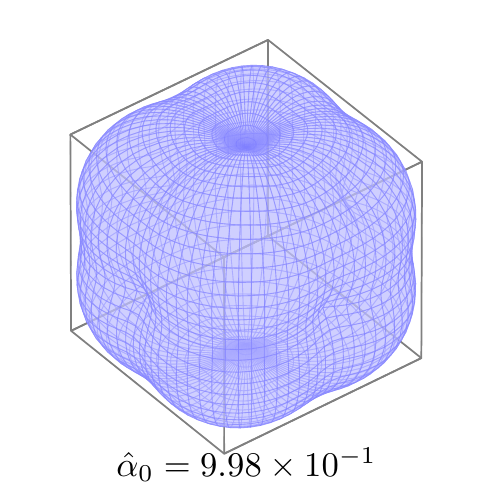}%
\hfill%
\includegraphics[page=1]{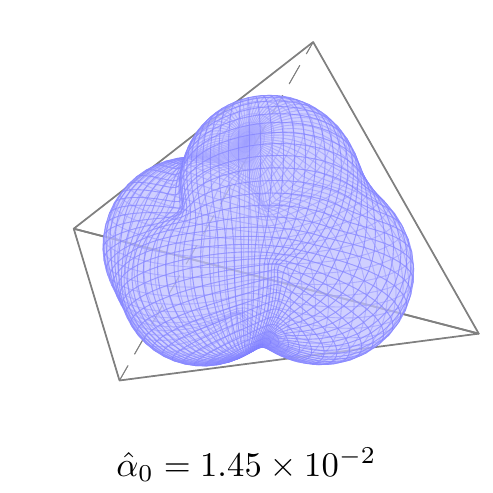}%
\hfill\null%
\caption{Iso-surface $\mathcal{S}_{\mathrm{iso}}\protect\fof{\polyvofdata_0}$ of plane base points from~\protect\refeqn{isoposition_surface} visualized for different cell types.}%
\label{fig:isoposition_surface}%
\end{figure}

\begin{figure}[htbp]
\null\hfill%
\includegraphics[page=3]{isoposition_square}%
\hfill%
\includegraphics[page=2]{isoposition_square}%
\hfill%
\includegraphics[page=1]{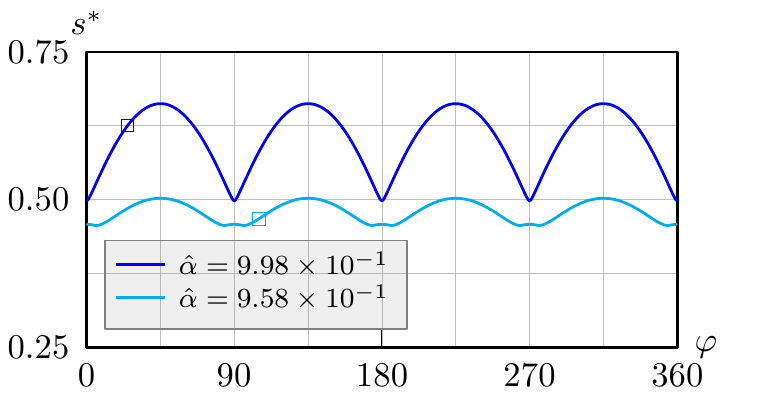}%
\hfill\null%
\caption{Two-dimensional iso-surface of plane base points from~\protect\refeqn{isoposition_surface_2d} for the unit square, where the left (\textcolor{cyan}{$\square$} in the right panel) and center (\textcolor{blue}{$\square$}) panel correspond to the section of $\mathcal{S}_{\mathrm{iso}}\protect\fof{\polyvofdata_0}$ in \protect\reffig{isoposition_surface} with the $xy$-plane. Note from the left panel that $\plicnormal$ is, in general, \textit{not} a normal of the iso-surface $\mathcal{S}_{\mathrm{iso}}$, since the "radius" $\signdistref$ of the spherical parametrization varies with $\varphi$ and $\theta$.%
}%
\label{fig:isoposition_surface_2d}%
\end{figure}

Using the face-based formulation of the volume fraction given in \refeqn{volume_computation_summation}, \refeqn{center_constraint} can be cast as %
\begin{align}
\centerconstraint{\signdist;\plicnormal}=%
\frac{1}{3}\sum\limits_{f}{\brackets{\volcoeffconst{f}+\signdist\volcoefflin{f}}\immersedarea{f}}-\polyvofdata_0.%
\label{eqn:center_volume_constraint_cubic}%
\end{align}
For some given normal $\plicnormal$, the volume fraction $\polyvof$ is a monotonous bounded continuous function of the signed distance $\signdist$, implying that the constraint in \refeqn{center_volume_constraint_cubic} admits a unique root, denoted $\signdistref\fof{\varphi,\theta;\polyvofdata_0}$. %
The set of basepoints %
\begin{align}
\mathcal{S}_{\mathrm{iso}}\fof{\polyvofdata_0}\defeq\set{\xbase+\signdistref\fof{\varphi,\theta;\polyvofdata_0}\plicnormal:\brackets{\varphi,\theta}\in\unitsphere}\label{eqn:isoposition_surface}%
\end{align}
forms a continuously differentiable hypersurface; cf.~\reffig{isoposition_surface} for an illustration. %
It is worth noting that away from points where $\plicplane$ intersects a vertex, $\mathcal{S}_{\mathrm{iso}}$ is arbitrarily smooth. In particular, $\mathcal{S}_{\mathrm{iso}}$ is globally continuously differentiable. %
Furthermore, the iso-surface is not necessarily contained in the parenting cell $\vfcell_k$, i.e.\ there are base points $\xref\in\mathcal{S}_{\mathrm{iso}}$ located outside of $\vfcell_k$, as can be seen from the center panel in \reffig{isoposition_surface_2d}. Configurations of this type are encountered for volume fractions $\polyvofdata_0$ close to zero or one. %

Evaluating \refeqn{center_volume_constraint_cubic} at $\signdist=\signdistref$, differentiating with respect to $\varphi$ and rearranging terms yields %
\begin{align}
\partial_\varphi\signdistref=-\evaluate{%
\brackets*{\sum_{f}{\signdist\partial_\varphi\volcoefflin{f}\immersedarea{f}+\brackets{\volcoeffconst{f}+\signdist\volcoefflin{f}}\partial_\varphi\immersedarea{f}}}%
\brackets*{\sum_{f}{\volcoefflin{f}\immersedarea{f}+\brackets{\volcoeffconst{f}+\signdist\volcoefflin{f}}\partial_\signdist\immersedarea{f}}}^{-1}}{\signdist=\signdistref}%
\label{eqn:plic_position_gradient}%
\end{align}
with an analogous expression for $\theta$. %
One could apply this strategy recursively to obtain the \textsc{Hessian} of $\signdistref$, which we omit for the reasons given in \refnote{higher_order_derivatives}. %
\paragraph{An example in two spatial dimensions}%
For the unit square $[0,1]^2$ with $\xbase=\frac{[1,1]\transpose}{2}$, $\plicnormal=[\cos\varphi,\sin\varphi]\transpose$ and $\refvof>\nicefrac{1}{2}$, the iso-surface from \refeqn{isoposition_surface} can be expressed analytically as %
\begin{align}
\signdistref\fof{\varphi}=%
\begin{cases}
\brackets{\refvof-\frac{1}{2}}\cos\varphi&\varphi\in[0,\varphi^\prime],\\%
\frac{\sin\varphi+\cos\varphi}{2}-\sqrt{2\brackets{1-\refvof}\cos\varphi\sin\varphi}&\varphi\in(\varphi^\prime,\frac{\pi}{4}]%
\end{cases}%
\quad\text{with}\quad%
\tan\varphi^\prime=2(1-\refvof),%
\label{eqn:isoposition_surface_2d}%
\end{align}
where the expressions for $\varphi\in(\frac{\pi}{4},2\pi]$ and $\refvof\leq\nicefrac{1}{2}$ are readily obtained from symmetry considerations. %
From \refeqn{isoposition_surface_2d}, it is easy to show that $\signdistref\fof{\cdot}$ is continuously differentiable at $\varphi^\prime$; cf.~\reffig{isoposition_surface_2d} for an illustration. %

%% file: 03_04_minimization_strategy.tex
\subsection{An iterative minimization strategy}\label{subsec:minimization_strategy}%
With the implicit treatment of \refeqn{center_constraint}, the vector of variables for the minimization reduces to $\vp\defeq\brackets[s]{\varphi,\theta}\transpose$, corresponding to the representation of the normal $\plicnormal$ in spherical coordinates. One obtains the unconstrained error function %
\begin{align}
\error\fof{\vp}\defeq%
\frac{1}{2}\sum\limits_{k=1}^{N}{\mu_k\brackets*{\polyvof_k-\polyvofdata_k}^2}%
\quad\text{with}\quad%
\polyvof_k\defeq\polyvof\fof{\vfcell_k;\signdistref\fof{\vp},\plicnormal\fof{\vp}},%
\label{eqn:error_functional_unconstrained}%
\end{align}
whose gradient and \textsc{Hessian}, respectively, with respect to $\vp$ (henceforth denoted as $\grad{}$ and $\hessian{}$) read %
\begin{align}
\grad{\error}=\sum\limits_{k=1}^{N}{\mu_k\brackets*{\polyvof_k-\polyvofdata_k}\grad{\polyvof_k}}%
\quad\text{and}\quad%
\hessian{\error}=\sum\limits_{k=1}^{N}{\mu_k\brackets*{\brackets*{\polyvof_k-\polyvofdata_k}\hessian{\polyvof_k}+\dyad{\grad{\polyvof_k}}{\grad{\polyvof_k}}}},%
\label{eqn:error_functional_unconstrained_gradient_hessian}%
\end{align}
where we defer to choice of the weights $\set{\mu_k}$ to subsection~\ref{subsubsec:choice_weights}. %
Evaluating the gradient of the volume fraction $\polyvof_k$ requires to account for $\signdistref$ (i.e.~the volume constraint) via the chainrule. With the partial derivatives from \refeqs{dalphads}, \refeqno{dalphadphi} and \refeqno{plic_position_gradient} one obtains %
\begin{align}
\grad{\polyvof_k}&=\evaluate{%
\brackets*[s]{%
\partial_\varphi\polyvof_k+\partial_\signdist\polyvof_k\partial_\varphi\signdistref,
\partial_\theta\polyvof_k+\partial_\signdist\polyvof_k\partial_\theta\signdistref%
}\transpose}{\signdist=\signdistref}.%
\end{align}
In the vicinity of the minimum, it is sensible to assume that the absolute deviation of the volume fractions is small, i.e.~$\max_k\abs{\polyvof_k-\polyvofdata_k}\ll1$, which allows to approximate the \textsc{Hessian} by the exterior product of the gradient, i.e. %
\begin{align}
\hessian{\error}\approx\sum\limits_{k=1}^{N}{\mu_k\dyad{\grad{\polyvof_k}}{\grad{\polyvof_k}}}\defeq*\tilde{\nabla}^2\error.%
\end{align}
Finally, we may formulate the minimization strategy: %
let $\vp^n$ and $\signdist^n\defeq\signdistref\fof{\vp^n}$ denote the $n$-th iteration of the normal and the associated constraint position (such that $\polyvof_0\fof{\vp^n;\signdist^n}=\polyvofdata_0$), respectively. Starting from $\vp^0$, whose computation we defer to subsection~\ref{subsubsec:initial_iteration}, we approximate a local minimum of \refeqn{error_functional_unconstrained} by a standard \textsc{Gauss-Newton}-method and iteratively update %
\begin{align}
\vp^{n+1}=\vp^n+\Delta\vp^n%
\quad\text{until}\quad%
\norm{\grad{\error}\fof{\vp^n}}_2<\epsilon,%
\label{eqn:newton_iteration}%
\end{align}
with the relative tolerance $\epsilon\defeq\num{e-4}$. %
In order to produce iterations that strictly reduce the error\footnote{If it is not possible to reduce the error, it holds that $\grad{\error}=0$ and the iteration terminates since a local minimum was found.}, the computation of the step $\Delta\vp$ employs a line-search approach along the direction %
\begin{align}
\Delta\bar{\vp}^n=-\evaluate{\brackets*{\brackets*[s]{\tilde{\nabla}^2\error}\inverse\grad{\error}}}{\vp=\vp^n,\signdistref=\signdist^n}.%
\label{eqn:newton_iteration_step}%
\end{align}
\refFig{flowchart_step_computation} contains a flowchart of the procedure, whose description requires two additional comments: %
\begin{enumerate}
\item The length of the original step $\Delta\bar{\vp}^n$ depends, among others, on the choice of the weights $\set{\mu_k}$, i.e.~$$\hot{\norm{\Delta\bar{\vp}^n}_2}=\max_k\set{\mu_k}.$$ Let $\Delta\bar{\vp}^n=a\ve_{\vp}$ with $a>0$ and $\norm{\ve_\vp}_2=1$. Due to the periodicity of the error functional with respect to the spherical angles $\varphi$ and $\theta$ only $a\operatorname{mod}2\pi$ will be decisive for the updated position, i.e.~$\error\fof{\vp^n+a\ve_{\vp}}\equiv\error\fof{\vp^n+\brackets{a\operatorname{mod}2\pi}\ve_{\vp}}$. %
Moreover, assuming that the initial iteration $\vp^0$ is located sufficiently close to the sought minimum suggests to further reduce the maximum step length. %
Intuition suggests to simply restrict the length to, say, $\frac{\pi}{4}$. However, considering the effect of the parametrization in spherical coordinates implies a loss of symmetry in the vicinity of the poles and, hence, renders this approach invalid. In order to properly account for the effect of the parametrization, we confine the set of reachable positions by a box around $\vp^n$, whose edges admit lengths $\Delta\hat{\theta}\defeq\frac{\pi}{4}$ and %
\begin{align}
\Delta\hat{\varphi}\fof{\theta}=\min\brackets*{2\pi,\Delta\hat{\theta}+\brackets{\pi-\Delta\hat{\theta}}\brackets*{\frac{2\theta-\pi}{\pi-\Delta\hat{\theta}}}^{N}}%
,\label{eqn:stepsize_box_clipping}%
\end{align}
where, henceforth let $N\defeq12$. As can be seen from the leftmost panel in \reffig{step_limitation_illustration}, choosing a large value for $N$ ensures that $\Delta\hat{\varphi}$ remains almost constant away of the equator but increases steeply in the vicinity of the poles. %
\begin{figure}[htbp]%
\null\hfill%
\includegraphics[page=4]{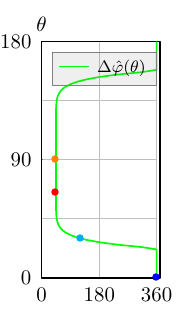}%
\includegraphics[page=3]{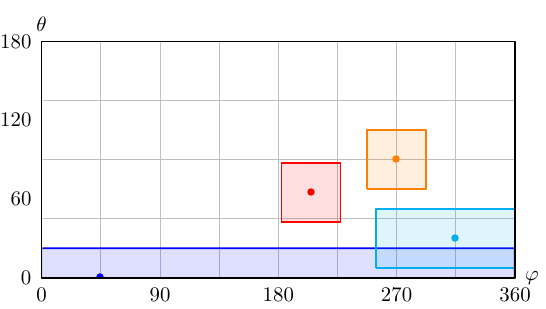}%
\includegraphics[page=5]{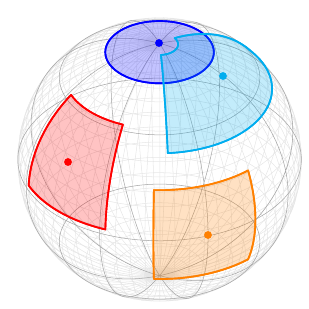}%
\hfill\null%
\caption{Step size limitiation by clipping $\Delta\hat{\vp}^n$ with a box of variable edge lengths (cf.~\protect\refeqn{stepsize_box_clipping}) centered at $\vp^n$ ($\bullet$), illustrated for different orientations on $\unitsphere$ with corresponding normals in $\setR^3$. The shaded regions show the positions that can be reached from the current orientation ($\bullet$) within a single limited step.}%
\label{fig:step_limitation_illustration}%
\end{figure}

\item There are configurations, mainly determined by bad initial conditions, for which the step size reduction cannot produce an iteration that reduces the error as desired. In these cases the step computation falls back to steepest descent, i.e.~$\Delta\bar{\vp}^n\defeq-\grad{\error}\fof{\vp^n}$, for which the line search is applied as well; cf.~\reffig{flowchart_step_computation}. %
\end{enumerate}

\begin{figure}[htbp]%
\null\hfill%
\includegraphics{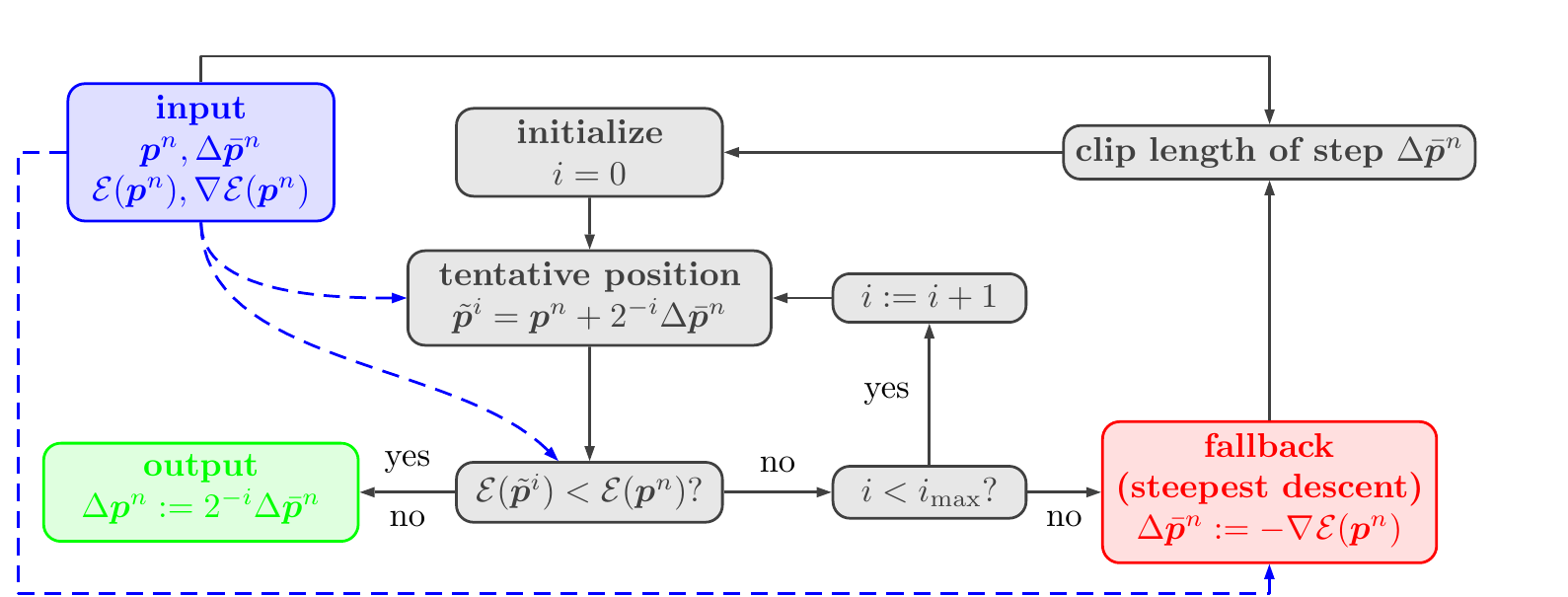}
\hfill\null%
\caption{Flowchart of the step computation; cf.~\protect\reffig{step_limitation_illustration} for details of the step computation. Throughout this work we employ $i_\mathrm{max}=6$.}%
\label{fig:flowchart_step_computation}%
\end{figure}

\begin{figure}[htbp]%
\null\hfill%
\includegraphics{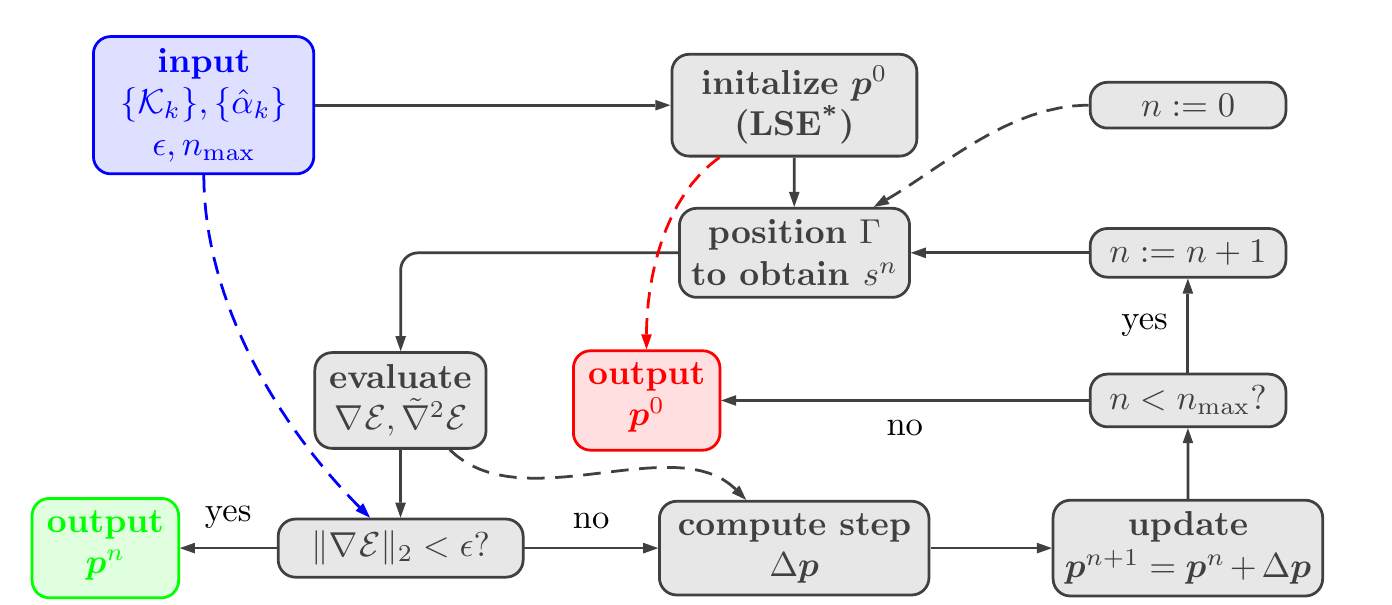}%
\hfill\null%
\caption{Flowchart of the main algorithm with step computation from \protect\reffig{flowchart_step_computation}.}%
\label{fig:flowchart}%
\end{figure}

After each iteration, the unique root of \refeqn{center_volume_constraint_cubic} yields the updated position $\signdist^{n+1}=\signdistref\fof{\vp^{n+1}}$, which we obtain using the highly efficient positioning algorithm of \citet{JCP_2021_fbip}; cf.~\refapp{implicit_bracketing} for a brief summary. %
\refFig{flowchart} contains the flowchart of the procedure described above. %

In terms of the numerical implementation, two comments are at order:%
\begin{enumerate}
\item Despite the fact that the $2\times2$ system of equations in \refeqn{newton_iteration} can be solved directly, we prefer to employ the LAPACK \cite{lapack99} routine \href{https://www.netlib.org/lapack/explore-html/d6/d0e/group__double_s_ysolve_ga9995c47692c9885ed5d6a6b431686f41.html}{\texttt{dsysv}} for reasons of numerical error control. %
\item In theory, all $\vp^n\in\setR^2$ are admissible. In practice, however, allowing for large numerical values of $\vp^n$ potentially degrades the accuracy of the arithmetic operations. Hence, after each iteration we assign $\vp^n\defeq\vp^n\operatorname{mod}2\pi$. %
\end{enumerate}

%
%
\subsubsection{On the choice of weights $\mu_k$}\label{subsubsec:choice_weights}%
The choice of the weights should ensure the desired intersection status (cf.~\reffig{intersection_status}) in the neighborhood\footnote{Recall that the term \textit{stencil} will be used synonymously.}, implying that the reconstructed plane $\plicplane$ must not intersect \textit{data-wise}\footnote{Within this work, we distinguish volume fractions that (i) are induced by intersecting a neighborhood with a plane $\plicplane$ of variable orientation (denoted $\polyvof_k$) and that (ii) are given (denoted $\polyvofdata_k$), where the term \textit{data-wise} refers to the latter.} bulk cells ($\polyvofdata_k\in[0,\voftolerance]\cup[1-\voftolerance,1]$). %
\begin{remark}[Geometrical interpretation]%
From a geometrical point of view, data-wise bulk cells carry crucial information in a stencil in that they allow to identify the a-priori geometrically admissible plane parameters. Thus, the choice of weights should ensure that the local reconstruction does not intersect immediate bulk neighbors of the center cell, corresponding to a penalty formulation of the constrained optimization problem. %
\end{remark}
While neighborhoods with sufficient spatial extension (type \textit{vertex}) typically contain a sufficient number of bulk cells, neighborhoods with small spatial extension (type \textit{face}) often contain -- if at all -- only a single bulk cell. %
On the other hand, if a given stencil does not admit any bulk cells, it must be regarded as a resolution problem since, geometrically speaking, any local linear reconstruction will inevitably convey a bias by the inability to intersect all cells of the stencil at once. %
A series of preliminary numerical experiments on tetrahedral meshes has shown that in those cases, the contribution to the gradient and the \textsc{Gauss-Newton} step is insufficent for qualitatively reasonable reconstruction. Hence, for \textit{face}-neighborhoods, bulk cells will be assigned a relatively large weight, say $\mu_k=\num{e9}$. %
\begin{remark}
It should be noted that, in general, the quadratic penalty method is not \textit{exact}, i.e.~the penalty parameters usually have to be chosen adaptively. However, our numerical experiments did not show any benefits of such an adaptive choice. On the contrary, for weights smaller than \num{e6}, a lot of unnecessary minimization steps with convergence to a non-compliant orientation were observed, suggesting indeed the existence of a treshold penalty parameter. 
\end{remark}

However, as stated above, the stencil does not necessarly contain a bulk cell. As we shall see below, the reconstruction is prone to produce outliers for those configurations. To avoid this, we dynamically extend the stencil, based on the volume fraction data $\set{\polyvofdata_k}$, i.e.~we include the first (by label) bulk cell in the iterated neighborhood %
\begin{align}
\neighborhood{}^\mathrm{ext}=\set{\neighborhood{u}:\vfcell_u\in\neighborhood{}\wedge\polyvofdata_u\in[0,\voftolerance]\cup[1-\voftolerance,1]}.\label{eqn:local_dynamic_stencil_extension}%
\end{align}
%
%
\subsubsection{On the initial value $\vp^0$}\label{subsubsec:initial_iteration}%
\begin{figure}[htbp]%
\null\hfill%
\includegraphics[page=1,height=4cm]{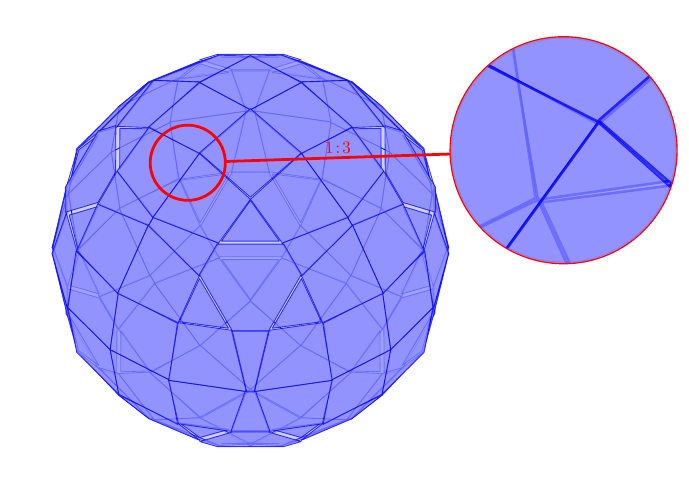}%
\includegraphics[page=2]{lse_weight_illustration}%
\includegraphics[page=2,height=4cm]{lse_cube_example}%
\includegraphics[page=1]{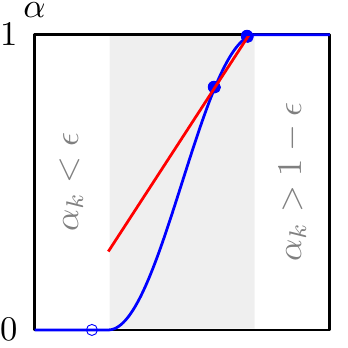}%
\hfill\null%
\\%
\null\hfill%
\includegraphics[page=3,height=4cm]{lse_cube_example}%
\includegraphics[page=4]{lse_weight_illustration}%
\includegraphics[page=4,height=4cm]{lse_cube_example}%
\includegraphics[page=3]{lse_weight_illustration}%
\hfill\null%
\caption{PLIC reconstruction with LSE normals with different weights (left colum: $\omega_k\equiv1$, right column: $\omega_k$ as in \protect\refeqn{lse_weights}) for sphere ($R_0=0.5445$, centered at $\vec{0}$) at equidistant spatial resolution (top row: $\Delta x=\nicefrac{1}{10}$, bottom row: $\Delta x=\nicefrac{1}{40}$). The panels qualitatively illustrate the effect of considering (\textcolor{blue}{$\bullet$}) or omitting (\textcolor{blue}{$\circ$}) bulk cells for the gradient approximation (red line) in one spatial dimension.}%
\label{fig:lse_reconstruction_weight_comparison}%
\end{figure}

The intial orientation $\vp^0$ is obtained by exploiting the fact that $\plicnormal\approx-\frac{\grad{\polyvof}}{\norm{\grad{\polyvof}}}$, where $\grad{}$ here denotes the gradient with respect to the spatial variable $\vx$. %
Once the inital value for the normal $\plicnormal$ is computed in the neighborhood center cell $\vfcell_0$, we solve the interface positioning problem \cite{JCP_2021_fbip} to obtain the signed distance $\signdist^0=\signdistref\fof{\vp^0}$. The estimation of $\grad{\polyvof}$ resorts to one of the following schemes: %
%
%
\paragraph{Least Square Error (LSE)}%
From the volume fraction data $\set{\polyvofdata_k}$ and the centroids $\set{\cellcentroid{k}}$, we assemble a linear least-squares problem for the average gradient of the volume fraction in the center cell $\average[\vfcell_0]{\grad{\polyvof}}$. Introducing $\Delta\polyvofdata_k\defeq\polyvofdata_k-\polyvofdata_0$ and $\Delta\cellcentroid{k}\defeq\cellcentroid{k}-\cellcentroid{0}$, one obtains %
\begin{align}
\average[\vfcell_0]{\grad{\polyvof}}&=%
\arg\min_{\vy}\frac{1}{2}\sum\limits_{k=1}^{N}{\omega_k\brackets*{\iprod*{\vy}{\Delta\cellcentroid{k}}-\Delta\polyvofdata_k}^2}\nonumber\\%
\label{eqn:lse_gradient}%
&=\tA\inverse\vb%
\quad\quad\quad\quad\quad\quad\text{with}\quad
\tA\defeq\sum\limits_{k=1}^{N}{\psi_k\dyad{\Delta\cellcentroid{k}}{\Delta\cellcentroid{k}}},%
\quad\text{and}\quad%
\vb\defeq\sum\limits_{k=1}^{N}{\psi_k\Delta\polyvofdata_k\Delta\cellcentroid{k}},%
\end{align}
where the weights %
\begin{align}%
\psi_k=%
\begin{cases}
0&\text{for}\quad\polyvofdata_k\in[0,\voftolerance)\cup(1-\voftolerance,1],\\%
1&\text{for}\quad\polyvofdata_k\in[\voftolerance,1-\voftolerance],%
\end{cases}\label{eqn:lse_weights}%
\end{align}
remove the influence of bulk cells. %
In the numerical implementation, we employ the LAPACK routine \href{http://www.netlib.org/lapack/explore-html/d6/d0e/group__double_s_ysolve_ga9995c47692c9885ed5d6a6b431686f41.html}{\texttt{dsysv}} to solve \refeqn{lse_gradient}. %
\begin{note}[Convention]\label{note:lse_name_convention}%
In section~\ref{sec:numerical_results}, the algorithm mentioned above will be referred to as (i) \textbf{LSE\textsuperscript{*}} for the choice of weights given in \refeqn{lse_weights} and as (ii) \textbf{LSE} for $\psi_k\equiv1$, i.e.~if bulk cells are considered. %
\end{note}
At this point it is crucial to note that the choice of weights may effect the reconstruction quality differently, depending on the spatial resolution. %
While the reconstruction quality profits from considering bulk cells in the LSE problem at lower resolutions, increasing the spatial resolution reduces and eventually spoils the gain. %
The set of graphs in \reffig{lse_reconstruction_weight_comparison} provides an illustration of the reason behind this phenomenon: %
apparently, when estimating the gradient from a small number of volume fractions which are unevenly distributed across $[0,1]$, taking into account non-intersected cells improves the result of the approximation (top row in \reffig{lse_reconstruction_weight_comparison}). On the other hand, if the volume fractions are distributed sufficiently densly in $[0,1]$, non-intersected cells adversely affect the accuracy of the estimation (bottom row). %

\begin{figure}[htbp]
\null\hfill%
\includegraphics{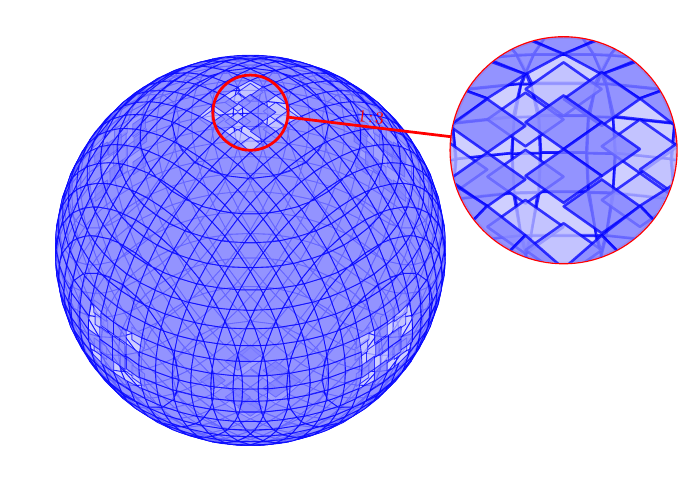}%
\hfill%
\includegraphics{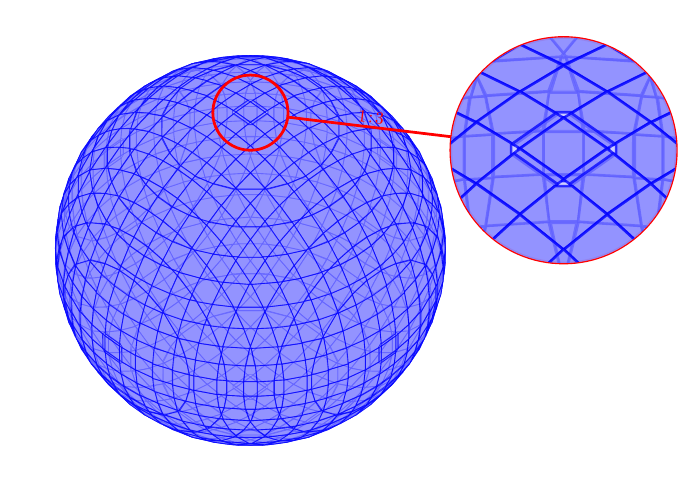}%
\hfill\null%
\caption{Seemingly random erroneous (left) and corrected (right) pattern of normal orientations produced by LSE gradient estimation.}%
\label{fig:lse_degenerate_coplanar}%
\end{figure}

For \textsc{Cartesian} structured meshes, there are configurations for which the weigths of \refeqn{lse_weights} induce a coplanar set of center coordinates of intersected stencil cells: %
e.g., if the interface intersects a $3\times3\times3$ stencil in a way that produces one interior, one intersected and one exterior $3\times3$ layer of cells, resembling a common configuration at the poles of a sphere. %
In those cases the matrix $\tA$ in \refeqn{lse_gradient} becomes singular, which, in an arithmetic context, yields $\abs{\det\tA}<\num{e-16}$. %
Carefully designed linear solvers, as those in LAPACK, still compute a meaningful normal in the sense that, locally, the associated plane is an approximation of the tangent space of the interface. %
However, the orientation, i.e.\ the sign of the normal, corresponds to $\sign{\det\tA}$ and, hence, cannot be determined reliably, resulting in a seemingly random pattern of inverted normals; cf.~the left panel in \reffig{lse_degenerate_coplanar}. %
The ambiguity can be removed by considering the inner product of the estimated normal with the vector connecting the centers of the cell with the smallest and largest volume fraction $\polyvofdata_k$, respectively: for negative values, the estimated normals must be inverted. In a geometrical sense, this ensures that the gradient points "outward". %

%
%
\paragraph{Gauss-Green (GG)}%
Employing the \textsc{Gaussian} divergence theorem, the cell-average of the gradient of the volume fraction $\polyvof$ can be expressed as a weighted sum of the associated face-averages, i.e.%
\begin{align}
\average[\vfcell_0]{\grad{\polyvof}}%
=\frac{1}{\abs{\vfcell_0}}\int\limits_{\vfcell_0}{\grad{\polyvof}\dvol}%
=\frac{1}{\abs{\vfcell_0}}\int\limits_{\partial\vfcell_0}{\polyvof\vn\darea}%
=\frac{1}{\abs{\vfcell_0}}\sum\limits_{f}{\vn_f\average[\vfface_f]{\polyvof}}.\label{eqn:gaussgreen_gradient}%
\end{align}
For the remainder of this work, the numerical approximation of the face-average $\average[\vfface_f]{\polyvof}$ resorts to a node-based scheme: %
for each node $\vfvert^\vfface_m$ of the face $\vfface_f$, all cells in the stencil containing that node contribute equally to the node average. The face average in turn emerges from the equally weighted average of the node averages. 

\begin{remark}[Continuing the literature survey]%
For a thorough overview on the topic of gradient estimation on unstructured meshes, the reader is referred to, e.g., \citet{TVCG_2011_acog}. %
\end{remark}

%% file: 04_design_experiments.tex
\section{Design of numerical experiments}\label{sec:design_experiments}%
\subsection{Meshes}\label{subsec:meshes}%
In what follows, we consider the domain $\domain=[-1,1]^3$, which is decomposed into cubes of equal size or tetrahedra, respectively. The latter are generated using the library \texttt{gmsh}, introduced in the seminal paper of \citet{IJNME_2009_gmsh}. For the purpose of the present paper, however, we only resort to some of the basic features of \texttt{gmsh}; see~\refapp{gmsh} and \reftab{mesh_characteristics} for further details. %

\begin{table}[htbp]%
\caption{Mesh characteristics.}%
\label{tab:mesh_characteristics}%
\null\hfill%
\subtable[Tetrahedron meshes; cf.~\refapp{gmsh}.]{%
\label{tab:tetrahedron_mesh_setup}%
\begin{tabular}{c|c|c|c}%
\textbf{resolution}&\textbf{char. length}&\# of cells&\# of faces\\%
$N$&$h=\frac{1}{N}$&&\\%
\hline%
\num{10}&$\num{1.00e-1}$&\num{4764}&\num{10261}\\%
\num{15}&$\num{6.66e-2}$&\num{15266}&\num{32160}\\%
\num{20}&$\num{5.00e-2}$&\num{33744}&\num{70322}\\%
\num{25}&$\num{4.00e-2}$&\num{64165}&\num{132750}\\%
\num{30}&$\num{3.33e-2}$&\num{108582}&\num{223550}\\%
\num{35}&$\num{2.85e-2}$&\num{171228}&\num{351057}%
\end{tabular}%
}%
\hfill%
\subtable[Equidistant cube meshes.]{%
\label{tab:cube_mesh_setup}%
\begin{tabular}{c|c|c|c}%
\textbf{resolution}&\textbf{char. length}&\# of cells&\# of faces\\%
$N$&$h=\frac{1}{N}$&$N^3$&$3N^2(N+1)$\\%
\hline%
\num{15}&\num{6.66e-02}&\num{3375}&\num{10800}\\%
\num{20}&\num{5.00e-02}&\num{8000}&\num{25200}\\%
\num{25}&\num{4.00e-02}&\num{15625}&\num{48750}\\%
$\vdots$&$\vdots$&$\vdots$&$\vdots$\\%
\num{70}&\num{1.42e-02}&\num{343000}&\num{1043700}\\%
\end{tabular}%
}%
\hfill\null%
\end{table}%
%
%
\subsection{Hypersurfaces}\label{subsec:experiments_hypersurfaces}%
In this work, we consider %
a sphere of radius $R_0=\nicefrac{4}{5}$, %
an oblate and prolate ellipsoid with semiaxes $\brackets{\nicefrac{4}{5},\nicefrac{4}{5},\nicefrac{2}{5}}$ and $\brackets{\nicefrac{1}{4},\nicefrac{1}{2},\nicefrac{3}{4}}$, respectively, as well as %
perturbed spheres. %
%
%
The latter can be parametrized in spherical coordinates as %
\begin{align}
\iface=\set{\vx_0+R\fof{\varphi,\theta}\ve_r:\fof{\varphi,\theta}\in\unitsphere}%
\quad\text{with}\quad%
R\fof{\varphi,\theta;\vciface}=\brackets*{\sum\limits_{l=0}^{\Liface}{\sum\limits_{m=-l}^{l}{\ciface{lm}\sphericalharm{}}}}^{\!\frac{1}{3}},%
\end{align}
where the description of the radius $R$ employs tesseral spherical harmonics $\sphericalharm{noarg}$ up to and including order $\Liface\in\setN$. %
The reason for expanding the third power of the radius instead of the radius itself is that the computation of the enclosed volume is considerably simplified, because $V_\iface=\ciface{00}\nicefrac{\sqrt{4\pi}}{3}$ then. %
The $(\Liface+1)^2$ coefficients $\ciface{lm}\sim\mathcal{N}(0,\sigma_0)$ are computed by the method of \citet{AMS_1958_anot}, i.e.
\begin{align}
\ciface{lm}=%
\begin{cases}
\sqrt{4\pi}R_0^3&l=0,\\%
\sqrt{\sigma_0}\sqrt{-2\log\gamma_1}\cos(2\pi\gamma_2)&l>0,%
\end{cases}\qquad\text{with}\qquad\gamma_{1,2}\sim\mathcal{U}(0,1),%
\end{align}
where the uniformly distributed random numbers $\gamma_{1,2}$ are generated by the \texttt{fortran} subroutine \texttt{random\_number()}. %
In this work, we consider perturbed spheres with base radius $R_0=\nicefrac{4}{5}$, modes $\Liface\in\set{3,6,9}$ and variance $\sigma_0=\num{5e-4}$; cf.~\reffig{hypersurface_illustration} for an illustration. %
\begin{figure}[htbp]%
\null\hfill%
\includegraphics[bb=0 0 1190 1100,height=4cm]{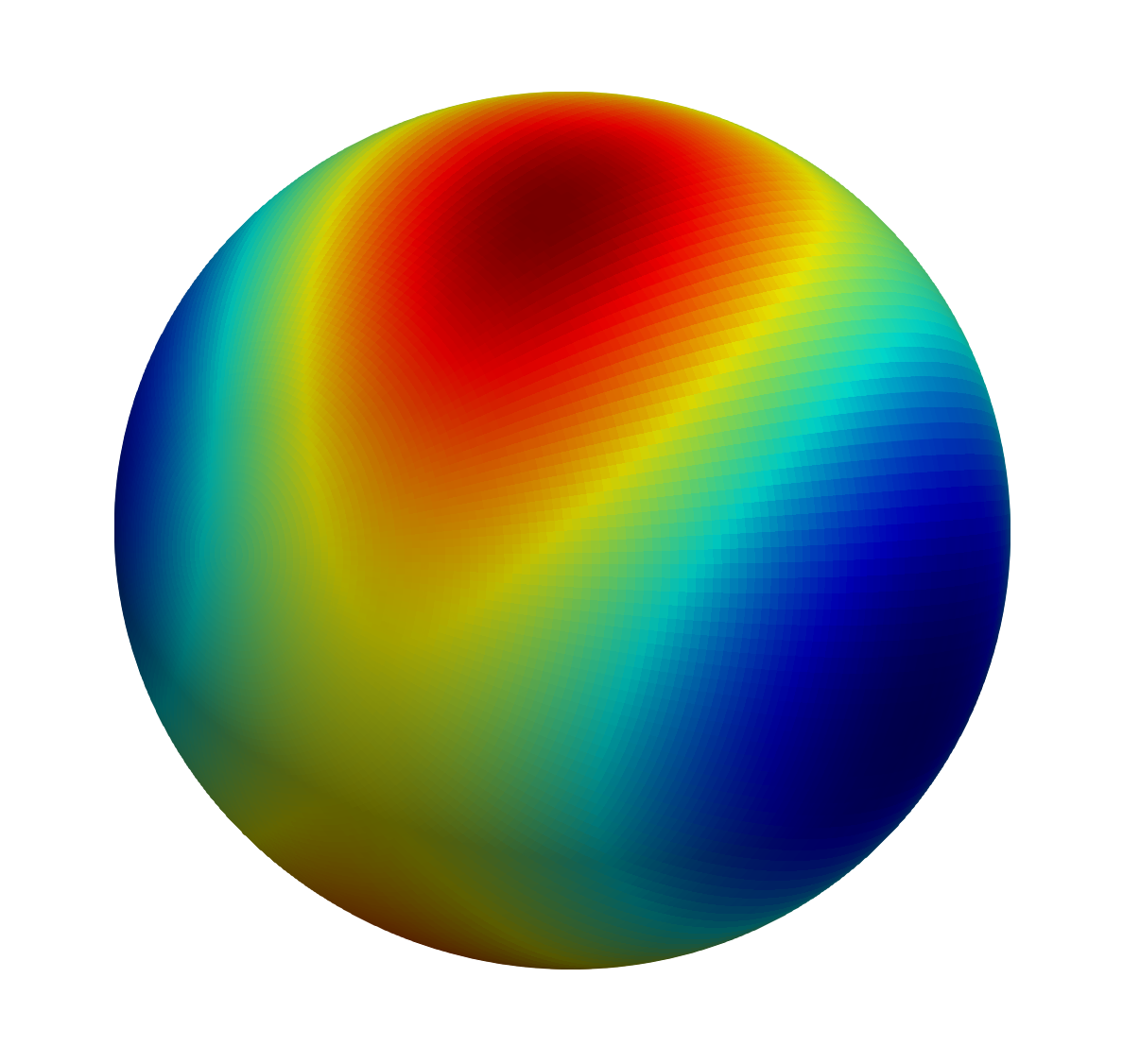}%
\includegraphics[page=1,height=4cm]{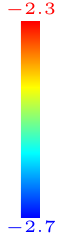}%
\hfill%
\includegraphics[bb=0 0 1190 1100,height=4cm]{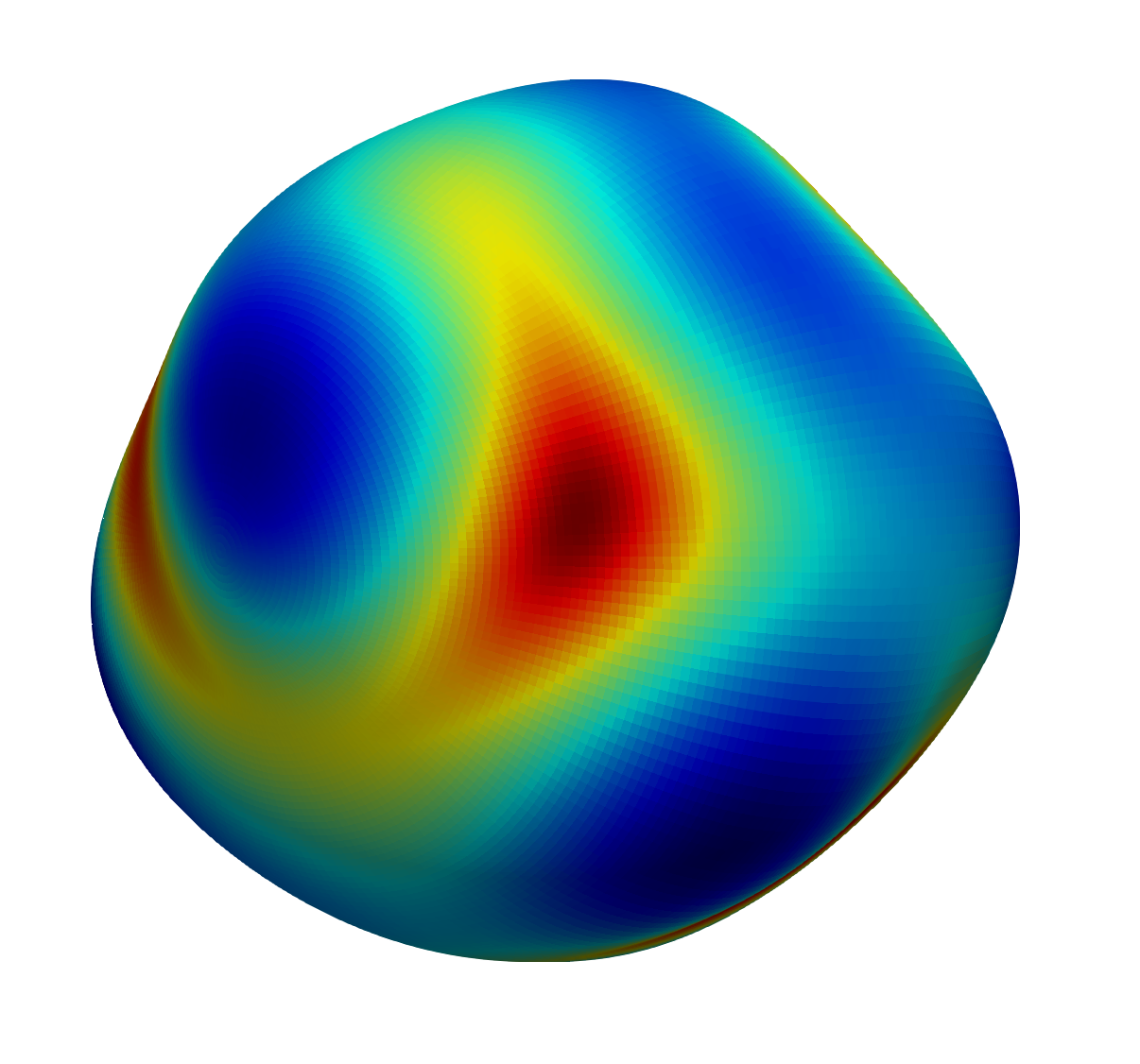}%
\includegraphics[page=2,height=4cm]{colorcurvature}%
\hfill\null%
\includegraphics[bb=0 0 1190 1100,height=4cm]{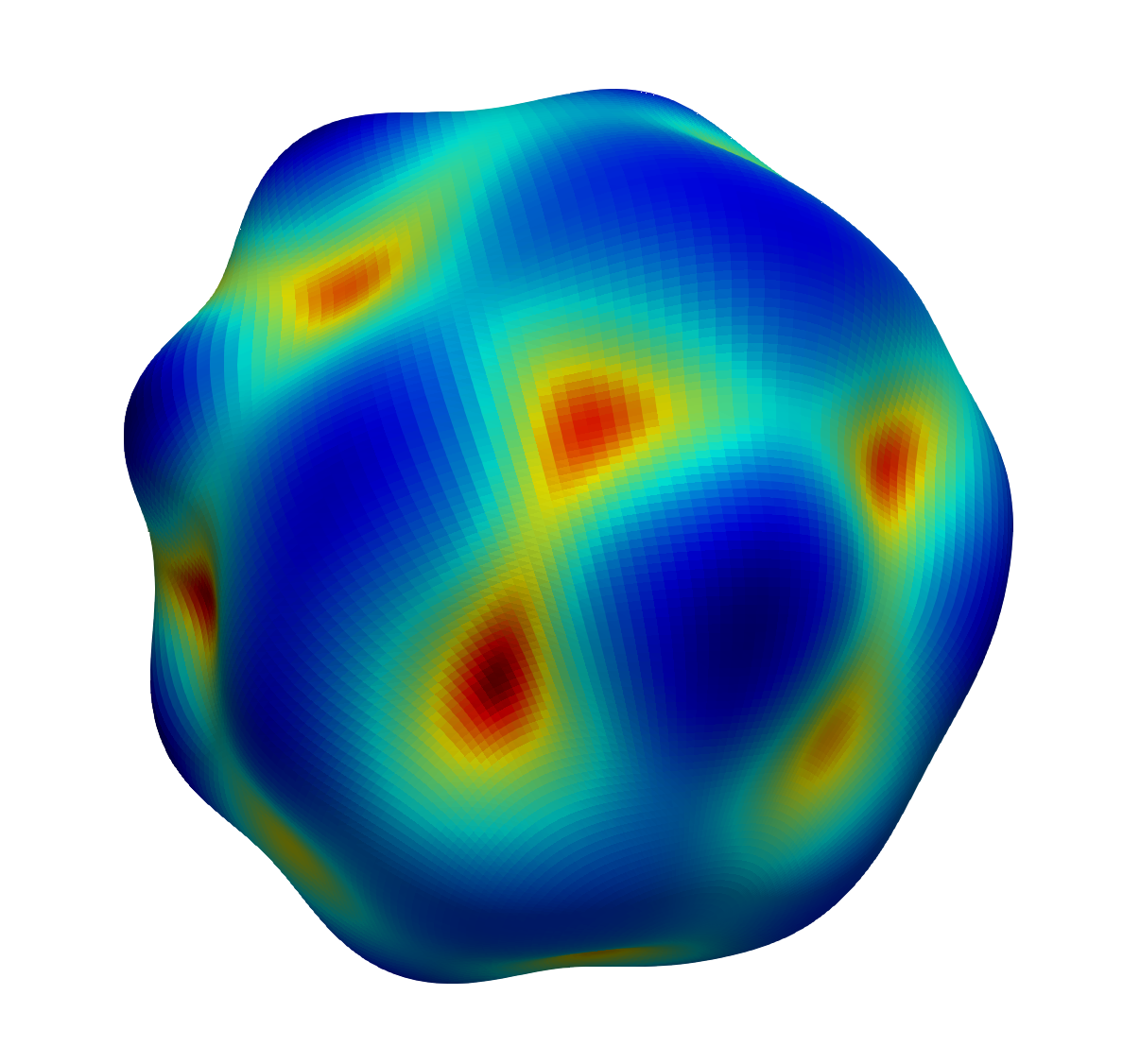}%
\includegraphics[page=5,height=4cm]{colorcurvature}%
\hfill\null%
\caption{Illustration of perturbed spheres ($R_0=\nicefrac{4}{5}$, $\Liface=\protect\set{3,6,9}$ and $\sigma_0=\num{5e-4}$), where the color indicates twice the mean curvature $2\curv_\iface=\curv_1+\curv_2$.}%
\label{fig:hypersurface_illustration}%
\end{figure}

All hypersurfaces are centered at $\vx_0=\vec{0}$ to ensure sufficient distance to the domain boundary. 
%
%
\subsection{Volume fractions}\label{subsec:initial_volume_fractions}%
For the initialization of the volume fraction, we use the method of \citet{JCP_2021_toai}, which can be summarized as follows: %
locally, i.e.~in each mesh cell, a principal coordinate system is computed in which the hypersurface can be approximated as the graph of an osculating paraboloid. The volume integrals are then transformed to curve integrals by recursive application of the \textsc{Gaussian} divergence theorem and can be easily approximated numerically by means of standard \textsc{Gauss-Legendre} quadrature. This face-based formulation enables the applicability to unstructured meshes. %
The initialization algorithm assigns a logical status (interior/intersected/exterior) to each computational cell of the mesh $\mathcal{M}$, which could be utilized to identify the interface cells within which the interface is reconstructed. In practice, however, those cells are identified based on their associated volume fraction $\polyvofdata_k$. Let %
\begin{align}
\mathcal{I}^\iface_\mathcal{M}\defeq\set{k:\vfcell_k\in\mathcal{M}\text{ such that }%
\voftolerance\leq\polyvofdata_{k}\leq1-\voftolerance%
}\label{eqn:mesh_interface_labels_logical}%
\end{align}
denote the list of labels of intersected cells with size $\Niface=\abs{\mathcal{I}^\iface_\mathcal{M}}$. %
\begin{remark}[Choice of the reconstruction tolerance]%
The choice of the tolerance $\voftolerance$ may have severe influence on the performance, especially if very small volume fraction values $\polyvofdata_k$ occur. %
\end{remark}
\begin{figure}[htbp]
\null\hfill%
\includegraphics[page=1]{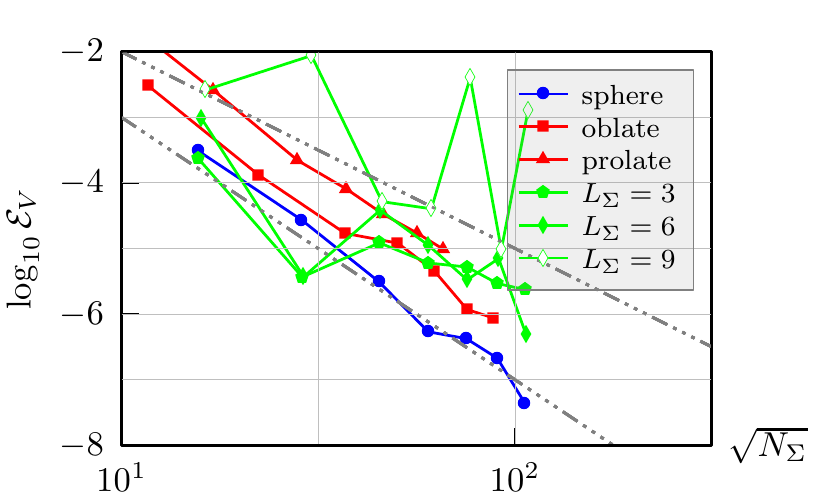}%
\hfill%
\includegraphics[page=2]{volerror_initial}%
\hfill\null%
\caption{Global relative volume error as a function of the number of intersected cells $\Niface$ for cube (left) and tetrahedron mesh (right); cf.~\reftab{mesh_characteristics} for detailed mesh information and \citet[Figure~10]{JCP_2021_toai}. The number of dots corresponds to the order of convergence; cf.~remark~\ref{remark:outlier}.}%
\label{fig:volume_fraction_initial_error}%
\end{figure}%
As can be seen from \reffig{volume_fraction_initial_error}, the algorithm of \citet{JCP_2021_toai} provides third- to fourth-order accurate initial volume fractions. Furthermore, the algorithm provides a locally linear or quadratic approximation of the hypersurface, implying that different qualities of initial volume fractions can be considered. %
\begin{remark}[Outliers in \reffig{volume_fraction_initial_error}]\label{remark:outlier}%
For non-spherical hypersurfaces, the quality of the local parabolic approximation underlying the volume computation varies with the hypersurface parameters and resolution: a combination of strong variation of the principal curvatures with large absolute values corresponds to a locally underresolved configuration. This is the case for the perturbed sphere with $\Liface=9$ (green $\diamondsuit$): at the largest resolution, a very small number of cells admits a large local volume error with coinciding sign, which accumulates to a comparatively large mesh-global volume error. %
\end{remark}


%% file: 05_numerical_results.tex
%
%
\section{Numerical results}\label{sec:numerical_results}%
The present section provides a thorough assessment of the proposed method resorting to the meshes and hypersurfaces described in section~\ref{sec:design_experiments}. In order to assess the accuracy of the reconstruction, subsection~\ref{subsec:convergence_measures} introduces some appropriate measures. %
Before subsection~\ref{subsec:mesh_convergence} investigates the curved hypersurfaces described in subsection~\ref{subsec:experiments_hypersurfaces}, we commence by applying the proposed algorithm to a halfspace in subsection~\ref{subsec:halfspace_reconstruction}. Since its planar boundary can be reconstructed exactly in theory, this setup provides an intrinsic benchmark. %

\begin{note}[Local error maps]\label{note:local_error_map}%
In what follows, we will frequently show local error maps associated to some neighborhood $\neighborhood{k}$, whose type (face/edge/vertex) and label $k$ will be given where necessary. %
Unless explicitly stated otherwise, they share the following properties: %
\begin{itemize}
\item The logarithm of the error from \refeqn{error_functional_unconstrained} is shown as a function of $\vp=\brackets[s]{\varphi,\theta}\transpose$ on an equidistant $2M\times M$ grid covering $\unitsphere$ with $M=60$, i.e.~ $\vp_{ij}=\brackets[s]{\pi\frac{2i-1}{M},\pi\frac{2j-1}{2M}}\transpose$ for $1\leq i\leq2M$ and $1\leq j\leq M$. The color legend at the top contains the bounds of the underlying discrete evaluation, i.e.~$\min_{i,j}{\log_{10}\error\fof{\vp_{ij}}}$ and $\max_{i,j}\log_{10}\error\fof{\vp_{ij}}$. %
\item The arrows at $\vp_{ij}$ indicate the direction of the \textsc{Gauss-Newton} step $\Delta\bar{\vp}^n$ from \refeqn{newton_iteration_step}. For the sake of clarity their length is fixed, i.e.~it does \textbf{not} correspond to $\norm{\Delta\bar{\vp}^n\fof{\vp_{ij}}}_2$. %
\item Starting from LSE initial conditions, the sequence of steps $(\vp^n)$ of the minimization procedure are marked as $\times$, where a dashed arrow connects the initial (\textbullet) and final iteration ($\oplus$). %
\item If shown, the line to the right of the error plot contains the distribution of volume fraction data $\polyvofdata_k$ in the stencil ($\times$), where $\otimes$ marks the center volume fraction $\polyvofdata_0$. %
\item \begin{minipage}[t]{.55\textwidth}%
In order to restrict the visualization of the iterations to $\unitsphere$, let %
\begin{align}
\fof{\varphi^n,\theta^n}\defeq%
\begin{cases}%
(\varphi^n+\pi,-\theta^n)&\theta^n<0,\\%
(\varphi^n+\pi,2\pi-\theta^n)&\theta^n>\pi,\\%
(\varphi^n,\theta^n)&\text{otherwise}.%
\end{cases}
\end{align}
As can be seen from the figure, the corresponding sequence of normals $(\plicnormal^n)$ is invariant under this transformation, while the corresponding sequence of spherical angles $(\vp^n)$ may exhibit "jumps". %
\end{minipage}
\hspace{.5cm}%
\begin{minipage}[t]{.3\textwidth}%
\null\hfill%
\includegraphics{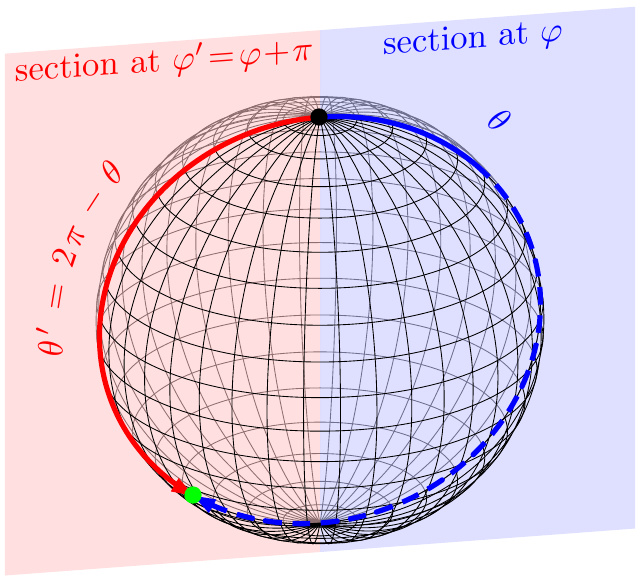}%
\hfill\null%
\end{minipage}
\end{itemize}
\end{note}
%
%
\input{05_01_convergence_measures}%
%
%
\input{05_02_halfspace_reconstruction}%
%
%
\input{05_03_convergence_curved}%
\clearpage%
%
%
\input{05_04_local_analysis}%
%
%
\input{05_05_influence_weights}%

%% file: 05_01_convergence_measures.tex
\subsection{On appropriate convergence measures}\label{subsec:convergence_measures}%
Since we wish to approximately reconstruct a hypersurface, it is crucial to define appropriate quality measures. After the normal reconstruction, the associated plane is positioned to match the prescribed volume fraction in each cell $\vfcell_k$. %
%
%
%
\paragraph{Normal alignment}%
As stated in subsection~\ref{subsec:initial_volume_fractions}, the volume fractions are initialized using the algorithm of \citet{JCP_2021_toai}, which for each cell $\vfcell_k$ computes a local quadratic approximation of the underlying hypersurface $\iface$, to a tangential principal coordinate system. Thus, as a by-product of the initialization, we obtain a discrete normal field $\ininormal{k}$ to which we compare the reconstructed normal $\recnormal{k}$ via %
\begin{align}
\Delta\vn_k\defeq\abs{1-\iprod{\recnormal{k}}{\ininormal{k}}}.%
\end{align}
Recall that the PLIC reconstruction approximates the continuous hypersurface $\iface$ by a finite set of planes $\set{\plicplane_k}$ associated to the mesh cells $\vfcell_k$, i.e.~$\iface\approx\bigcup_k{\plicplane_k}$. This suggests to compute a weighted average %
\begin{align}
\average[\plicplane]{\Delta\vn}\defeq\frac{\sum_{k=1}^{\Niface}{\abs{\plicplane_k}\Delta\vn_k}}{\sum_{k=1}^{\Niface}{\abs{\plicplane_k}}}%
\quad\text{with}\quad%
\abs{\plicplane_k}=\evaluate{\partial_\signdist\polyvof_k}{\signdist=\signdistref}\abs{\vfcell_k}%
.%
\label{eqn:normal_alignment}%
\end{align}
%
%
\paragraph{Symmetric volume difference}%
The symmetric volume difference consitutes a common measure in the literature \cite{JCP_2007_mmir,JCP_2009_asoa}. For two given planes, say $\primaryplicplane$ and $\secondaryplicplane$, it corresponds to the volume of the part of $\vfcell_k$ for which the halfspace assignment mismatches, i.e. %
\begin{align}
\Delta\polyvol\fof{\vn_1,\vn_2}\defeq%
\abs{\vfcell_k\cap\poshalfspace[1]{}\cap\neghalfspace[2]{}}+\abs{\vfcell_k\cap\poshalfspace[2]{}\cap\neghalfspace[1]{}}%
\quad\text{with}\quad%
\Delta\polyvol\fof{\vn_1,\vn_2}=\Delta\polyvol\fof{\vn_2,\vn_1}\geq0,%
\end{align}
where the numerical computation is carried out using the algorithm work of \citet{JCP_2021_etmp}. %
\refFig{illustration_symmetric_volume_difference} provides an illustration in two spatial dimensions. %
\begin{figure}[htbp]
\null\hfill%
\includegraphics[page=1]{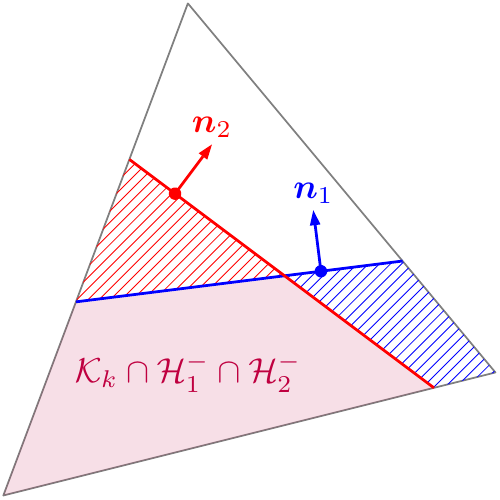}%
\hfill%
\includegraphics[page=2]{illustration_halfspace}%
\hfill%
\includegraphics[page=3]{illustration_halfspace}%
\hfill\null%
\caption{Illustration of volume overlap (purple shaded in the left panel) induced by two planes parametrized by their normals $\vn_{i}$ with associated positive ($\poshalfspace[i]{}$) and negative ($\neghalfspace[i]{}$) halfspaces. The symmetric difference (hatched) corresponds to the union of sequential truncations with order $\protect\brackets{\textcolor{blue}{1},\textcolor{red}{2}}$ (center) and $\protect\brackets{\textcolor{red}{2},\textcolor{blue}{1}}$ (right).}%
\label{fig:illustration_symmetric_volume_difference}%
\end{figure}
In consistency with the normal alignment $\Delta\vn_k$ introduced in \refeqn{normal_alignment}, the remainder of this work resorts to the symmetric volume fraction difference with respect to the initialization normal $\ininormal{k}$. In order to obtain a non-dimensional measure we normalize by the volume enclosed by the hypersurface, i.e. 
\begin{align}
\average[\plicplane]{\Delta\polyvol}\defeq\frac{\sum_{k=1}^{\Niface}\Delta\polyvol_k}{V_\iface}%
\quad\text{with}\quad%
V_\iface=\abs{\set{\vx\in\domain:\lvlset_\iface\fofx\leq0}}.%
\label{eqn:symmetric_volume_difference}%
\end{align}
It is crucial to note at this point that the evaluation of \refeqs{normal_alignment} and \refeqno{symmetric_volume_difference} for the normal alignment and symmetric volume difference, respectively, (i) is carried out only for \textit{data-wise} intersected cells and (ii) does not depend on the neighborhood type. %

%% file: 05_02_halfspace_reconstruction.tex
\subsection{Halfspace reconstruction}\label{subsec:halfspace_reconstruction}%
The present subsection investigates the reconstruction of the boundary of a halfspace parametrized by a base point $\xref^\mathrm{ref}$ and normal $\plicnormal^\mathrm{ref}$ (corresponding to $\varphi^\mathrm{ref}$ and $\theta^\mathrm{ref}$ with $\vp^\mathrm{ref}=\brackets[s]{\varphi^\mathrm{ref},\theta^\mathrm{ref}}\transpose$); cf.~\refeqn{plane_parametrization}. %
As stated in subsection~\ref{subsec:initial_volume_fractions}, the initial volume fractions are locally third-order accurate for polyhedral meshes, implying that halfspaces can be represented exactly up to machine precision. %
As a consequence, the corresponding error in \refeqn{error_functional_unconstrained} exhibits a \textit{global} minimum at $\vp^\mathrm{ref}$ with $\error\fof{\vp^\mathrm{ref}}=0$ up to machine precision, irrespective of the choice of the weights $\set{\mu_k}$. %
In order to assess the proposed algorithm, we consider the plane parametrized by the base point $\xref^\mathrm{ref}=\brackets[s]{0.4534,0.5442,0.4330}\transpose$ and normal $\plicnormal^\mathrm{ref}=\frac{\brackets[s]{1,-3,6}\transpose}{\sqrt{46}}$. We exemplarily show the detailled results for a cuboid and tetrahedron mesh at the lowest resolution $N=5$, while the findings hold in a qualitative sense for all meshes and resolutions given in \reftab{mesh_characteristics}. %

\begin{figure}[htbp]%
\centering%
\subfigure[cuboid mesh ($N=5$) with vertex neighborhood]{%
\label{fig:halfspace_vert_cuboid}%
\null\hfill%
\includegraphics[page=2,trim=0cm 1cm 0cm .25cm,height=7cm]{halfspace_vert_cuboid_illustration}%
\hfill%
\includegraphics[page=1,trim=0cm 1cm 0cm .25cm,height=7cm]{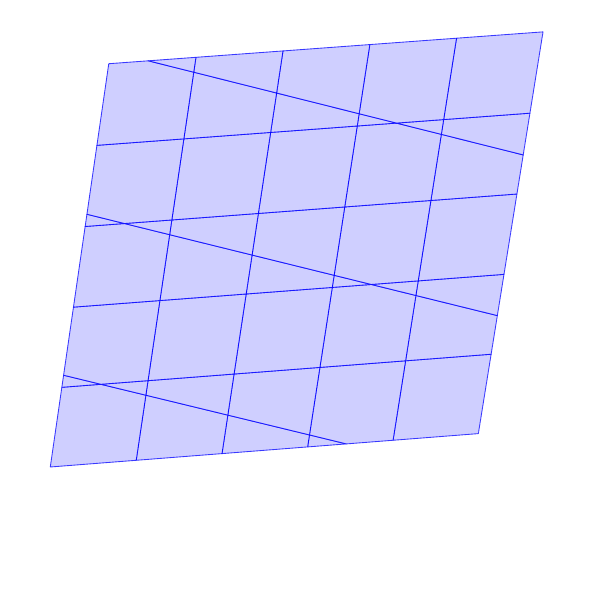}%
\hfill\null%
}%
\\%
\subfigure[tetrahedral mesh ($N=5$) with edge neighborhood]{%
\label{fig:halfspace_edge_tet}%
\null\hfill%
\includegraphics[page=2,trim=0cm 1cm 0cm .25cm,height=7cm]{halfspace_edge_tet_illustration}%
\hfill%
\includegraphics[page=1,trim=0cm 1cm 0cm .25cm,height=7cm]{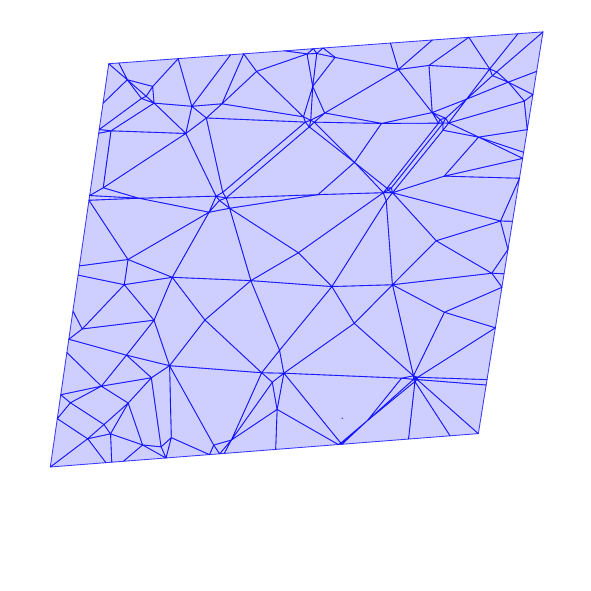}%
\hfill\null%
}%
\\%
\subfigure[tetrahedral mesh ($N=5$) with face neighborhood; cf.~\protect\reffig{halfspace_face_tet_closeup} for details]{%
\label{fig:halfspace_face_tet}%
\null\hfill%
\includegraphics[page=2,trim=0cm 1cm 0cm .25cm,height=7cm]{halfspace_face_tet_illustration}%
\hfill%
\includegraphics[page=1,trim=0cm 1cm 0cm .25cm,height=7cm]{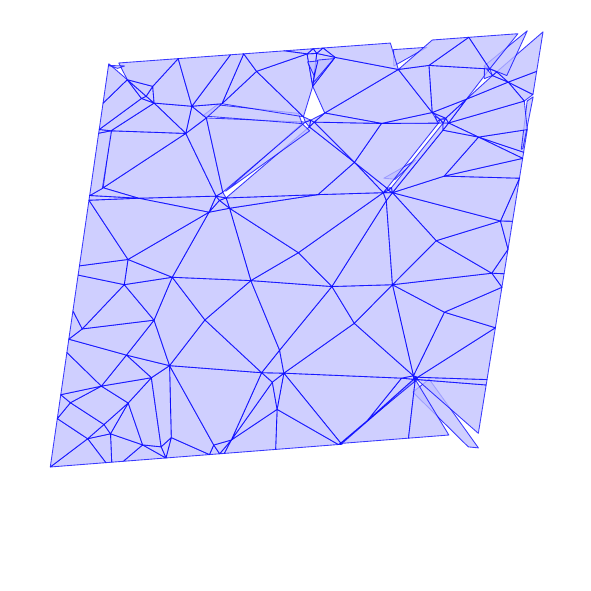}%
\hfill\null%
}%
\caption{Reconstruction of halfspace boundary with LSE (left) and present approach (right).}%
\label{fig:halfspace_reconstruction_illustration_tet_cube}%
\end{figure}

\begin{figure}[htbp]
\null\hfill%
\includegraphics[page=3]{halfspace_face_tet_contour_magnified}%
\hfill%
\includegraphics[page=3,trim=0cm 1cm 0cm 1cm,height=5cm]{halfspace_face_tet_outlier}%
\hfill\null%
\\%
\null\hfill%
\includegraphics[page=1]{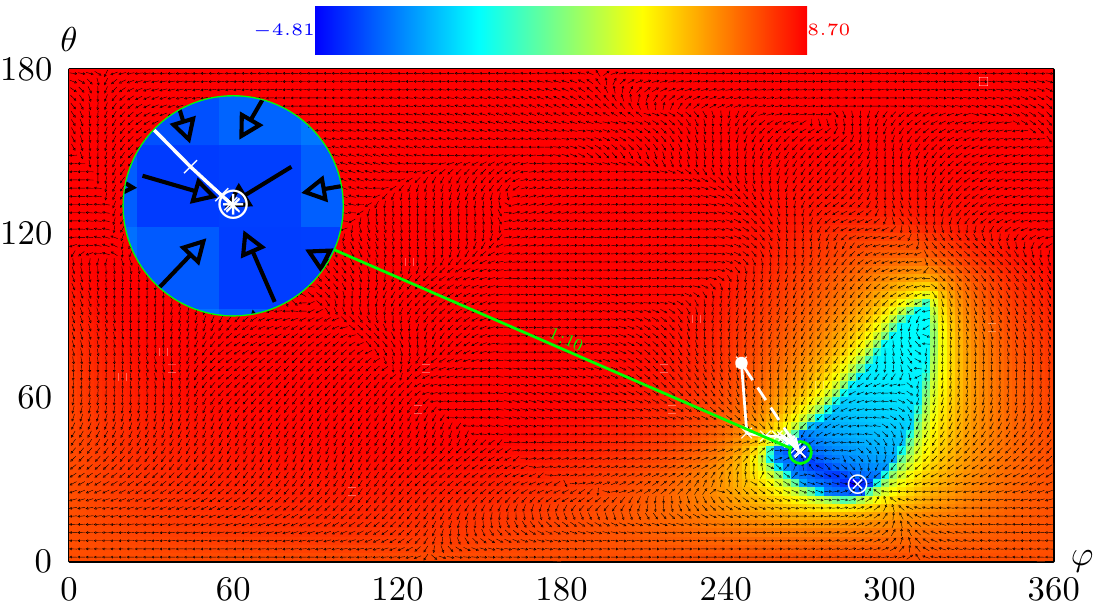}%
\hfill%
\includegraphics[page=1,trim=0cm 1cm 0cm 1cm,height=5cm]{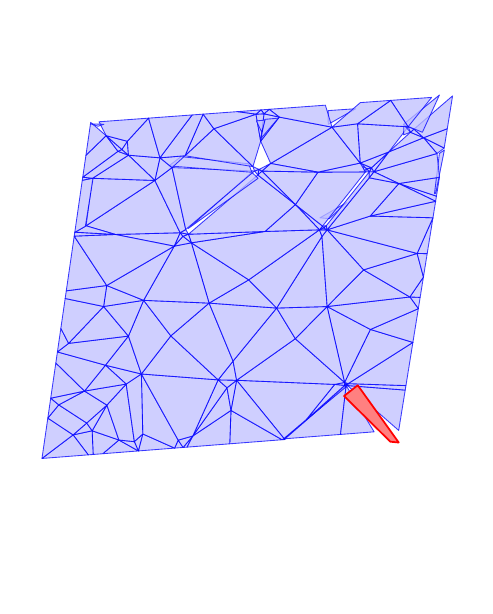}%
\hfill\null%
\\%
\null\hfill%
\includegraphics[page=2]{halfspace_face_tet_contour_magnified}%
\hfill%
\includegraphics[page=2,trim=0cm 1cm 0cm 1cm,height=5cm]{halfspace_face_tet_outlier}%
\hfill\null%
\caption{Local error maps (cf.~\refnote{local_error_map}) of halfspace boundary reconstruction on tetrahedral mesh ($N=5$, face neighborhood, cf.~\protect\reffig{halfspace_face_tet}) with reference orientation $\vp^\mathrm{ref}\approx\protect\brackets[s]{288^{\circ},28^{\circ}}\transpose$ ($\otimes$) and associated PLIC patches. While the patch in the top row exhibits the \textit{compliant} (green) orientation, the \textit{outliers} (red) in the center and bottom row are caused by convergence to a non-compliant local minimum or exceeding the maximum number of iterations without convergence, respectively.}%
\label{fig:halfspace_face_tet_closeup}%
\end{figure}

\refFig{halfspace_reconstruction_illustration_tet_cube} illustrates and compares the PLIC reconstruction obtained from LSE and the proposed approach. %
The main findings can be summarized as follows: %
\begin{itemize}
\item For the cuboid mesh with vertex neighborhood, see~\reffig{halfspace_vert_cuboid}, there are no outliers and the errors are in the order of $\num{e-16}$. This can be expected from the large variations of the error in the vicinity of the minimum, forming a comparatively strong attractor (left panel in \reffig{halfspace_minimum_comparison}). %
\item The tetrahedron mesh with edge neighborhood shown in \reffig{halfspace_edge_tet} admits no outliers as well. However, the errors exhibit values of order \num{e-6}, which can be attributed to a less pronounced (compared to the cuboid case) minimum of the error functional. %
\begin{figure}[htbp]%
\null\hfill%
\includegraphics[page=2]{halfspace_minimum_comparison}%
\hfill%
\includegraphics[page=1]{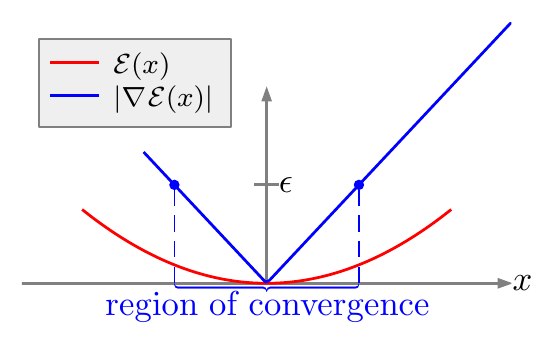}%
\hfill\null%
\caption{Extension of region of convergence for functions $\error\fof{x}$ with different variation in the vicinity of the minimum.}%
\label{fig:halfspace_minimum_comparison}%
\end{figure}%

As can be seen from the schematic illustration in \reffig{halfspace_minimum_comparison}, this induces a larger region around the minimum within which an iteration will be considered converged. From a geometric point of view, the smaller spatial extension of the neighborhood $\neighborhood{}$ (\textit{vertex} for the cuboid vs.~\textit{edge} for the tetrahedron) makes the minimization problem less stiff (right panel in \reffig{halfspace_minimum_comparison}). %
\item The reconstruction on the tetrahedron mesh with face neighborhood in \reffig{halfspace_face_tet} shows several outliers, which are caused by two mechanisms: (i) convergence to a non-compliant local minimum due to qualitatively bad initial conditions and (ii) non-convergence after the maximum number of iterations. \refFig{halfspace_face_tet_closeup} illustrates the aforementioned mechanisms. As can be seen, non-convergent cases feature iterations $\vp^n$ that follow valleys of the error function that exhibit very small gradients and, thus, small \textsc{Gauss-Newton} steps. %
\end{itemize}

%% file: 05_03_convergence_curved.tex
\subsection{Mesh convergence for curved hypersurfaces}\label{subsec:mesh_convergence}%
\begin{figure}[htbp]
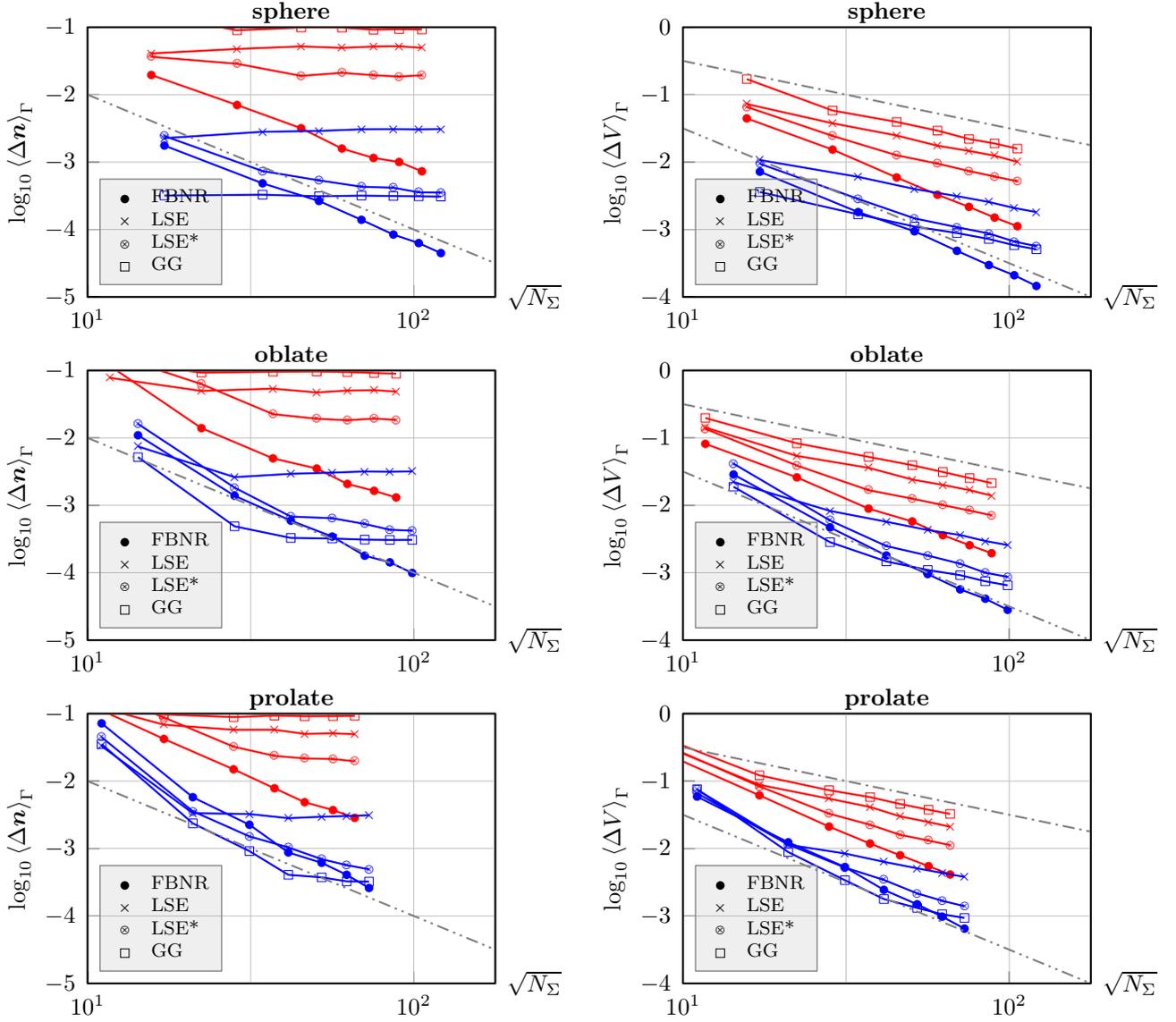

\multido{\i=1+1}{3}{%
\null\hfill%
\includegraphics[page=\i]{normaldeviation}%
\hfill%
\includegraphics[page=\i]{symvoldiff}%
\hfill\null%
\ifthenelse{\equal{\i}{3}}{}{\\}%
}%
\caption{Normal alignment (left; cf.~\protect\refeqn{normal_alignment}) and symmetric volume difference (right; cf.~\protect\refeqn{symmetric_volume_difference}) for different reconstruction algorithms as a function of spatial resolution (cf.~subsection~\ref{subsec:meshes} for details) for sphere, oblate and prolate ellipsoid on tetrahedron (red) and cuboid (blue) mesh.}%
\label{fig:mesh_convergence_ellipsoids}%
\end{figure}

\begin{figure}[htbp]
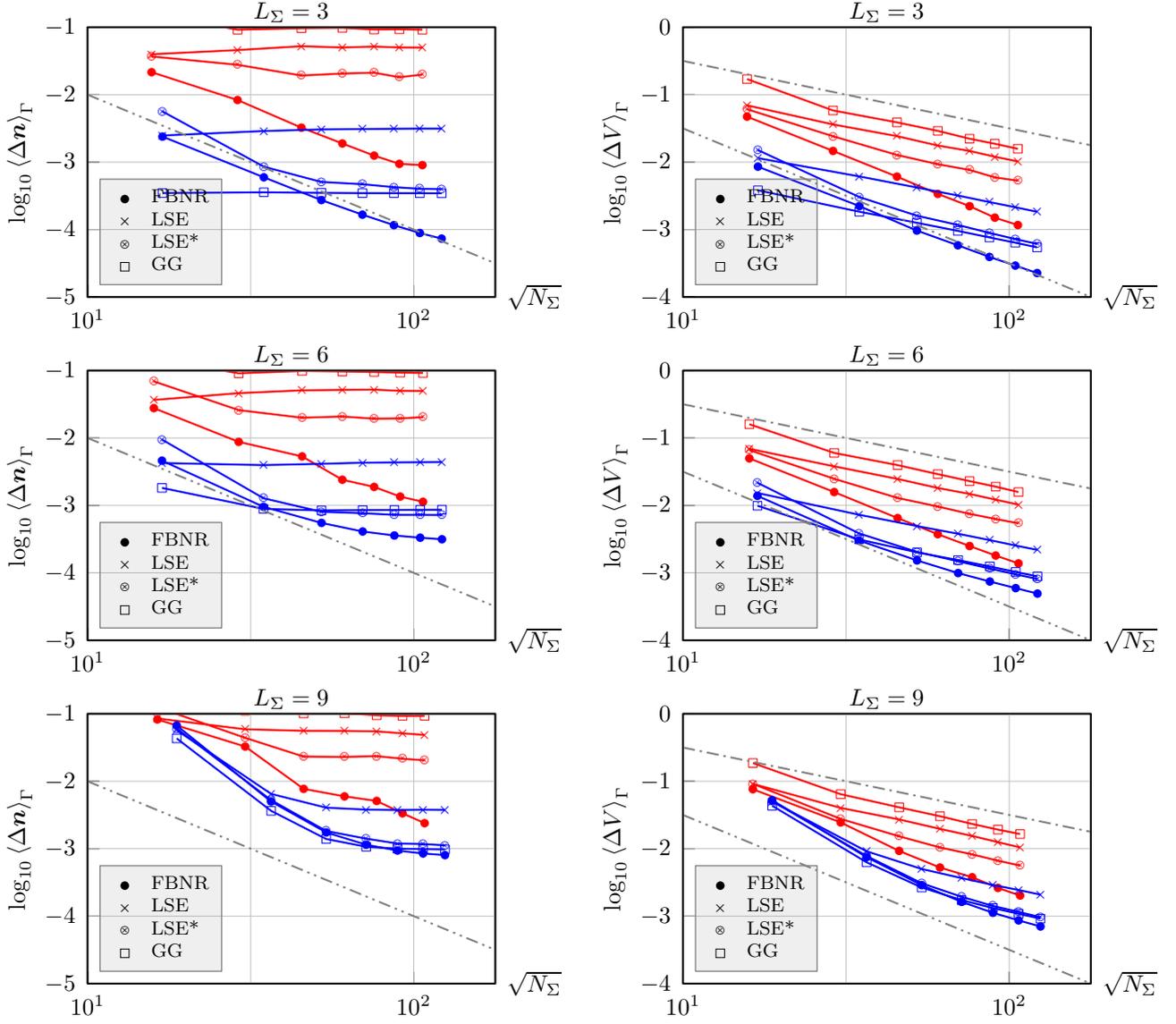

\multido{\i=4+1}{3}{%
\null\hfill%
\includegraphics[page=\i]{normaldeviation}%
\hfill%
\includegraphics[page=\i]{symvoldiff}%
\hfill\null%
\ifthenelse{\equal{\i}{6}}{}{\\}%
}%
\caption{Normal alignment (left; cf.~\protect\refeqn{normal_alignment}) and symmetric volume difference (right; cf.~\protect\refeqn{symmetric_volume_difference}) for different reconstruction algorithms as a function of spatial resolution (cf.~subsection~\ref{subsec:meshes} for details) for perturbed sphere ($R_0=\nicefrac{4}{5}$, $\vx_0=\vec{0}$, $\sigma_0=\num{5e-4}$) on tetrahedron (red) and cuboid (blue) mesh.}%
\label{fig:mesh_convergence_perturbed_sphere}%
\end{figure}

\begin{figure}[htbp]
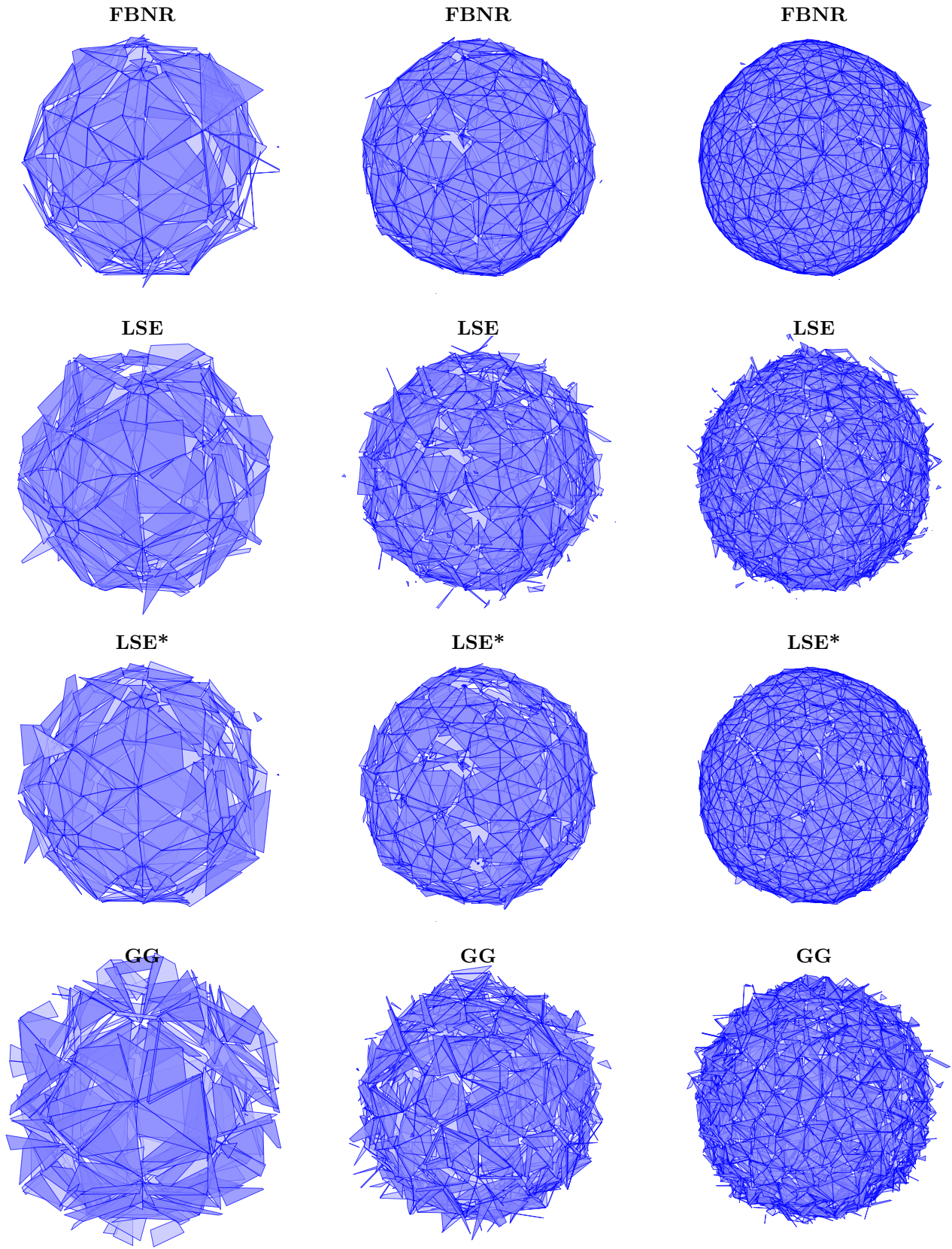

\multido{\i=1+1}{4}{
\null\hfill%
\includegraphics[page=\i]{illustration_glsasurf3D_R0_0-dot-800_L06_var_5-dot-0E-04_face_1-dot-0E+09_tetmesh_0005}%
\hfill%
\includegraphics[page=\i]{illustration_glsasurf3D_R0_0-dot-800_L06_var_5-dot-0E-04_face_1-dot-0E+09_tetmesh_0010}%
\hfill%
\includegraphics[page=\i]{illustration_glsasurf3D_R0_0-dot-800_L06_var_5-dot-0E-04_face_1-dot-0E+09_tetmesh_0015}%
\hfill\null%
\ifthenelse{\equal{\i}{4}}{}{\\}%
}%
\caption{Comparison of algorithms of perturbed sphere ($\Liface=6$) on tetrahedron mesh with $N\in\protect\set{5,10,15}$.}
\label{fig:illustration_tet_algo_L06}%
\end{figure}

\begin{figure}[htbp]
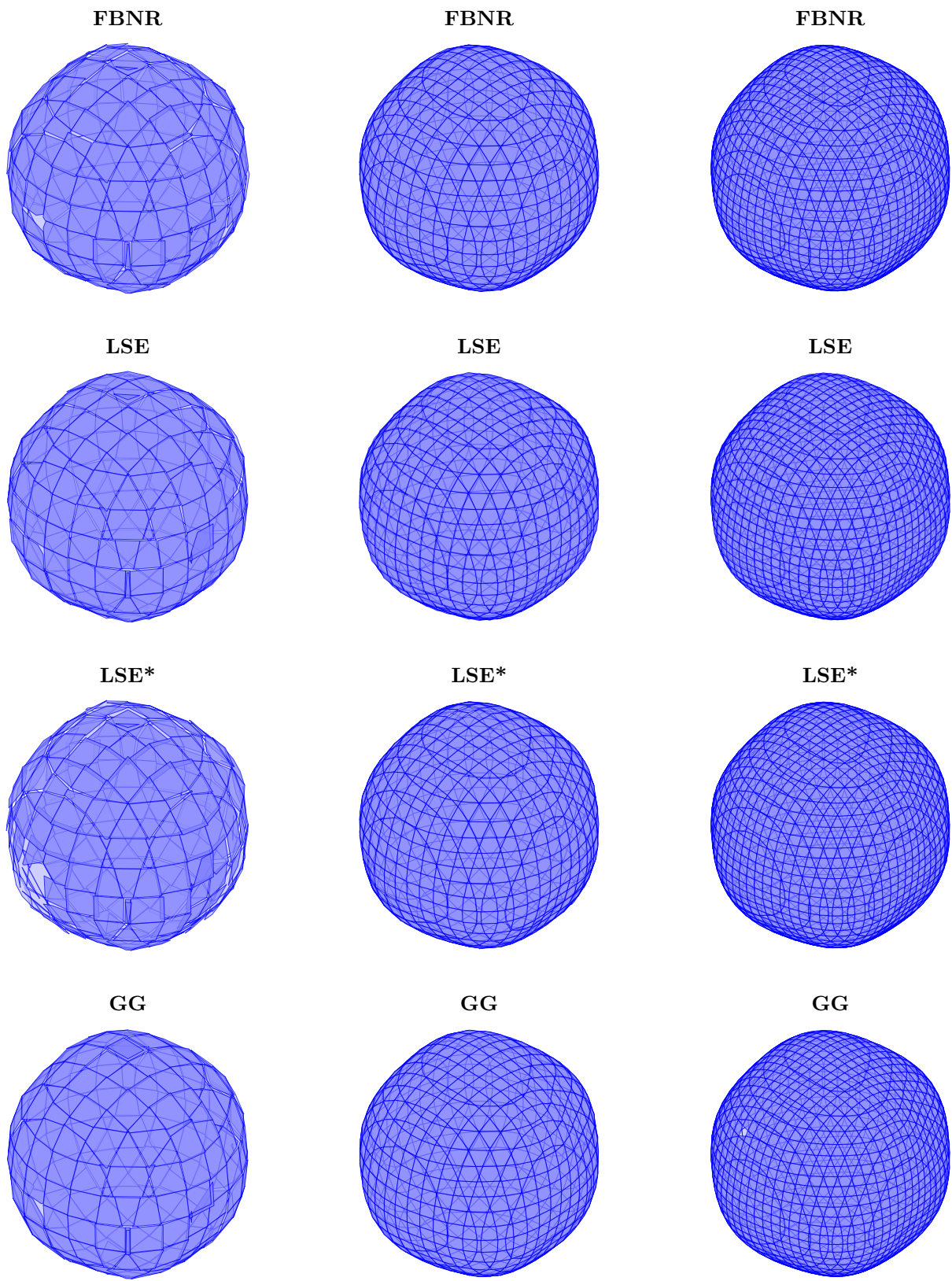

\multido{\i=1+1}{4}{
\null\hfill%
\includegraphics[page=\i]{illustration_glsasurf3D_R0_0-dot-800_L06_var_5-dot-0E-04_vert_1-dot-0E+00_cuboidmesh_0010}%
\hfill%
\includegraphics[page=\i]{illustration_glsasurf3D_R0_0-dot-800_L06_var_5-dot-0E-04_vert_1-dot-0E+00_cuboidmesh_0020}%
\hfill%
\includegraphics[page=\i]{illustration_glsasurf3D_R0_0-dot-800_L06_var_5-dot-0E-04_vert_1-dot-0E+00_cuboidmesh_0030}%
\hfill\null%
\ifthenelse{\equal{\i}{4}}{}{\\}%
}%
\caption{Comparison of algorithms of perturbed sphere ($\Liface=6$) on cuboid mesh with $N\in\protect\set{10,20,30}$.}
\label{fig:illustration_cube_algo_L06}%
\end{figure}

\begin{figure}[htbp]%
\null\hfill%
\includegraphics[page=1]{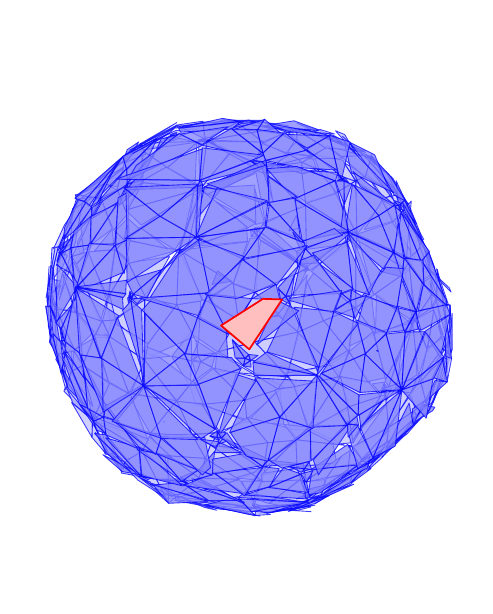}%
\hfill%
\includegraphics[page=2]{L06_lowres_outlier}%
\hfill%
\includegraphics[page=3]{L06_lowres_outlier}%
\hfill\null%
\\%
\null\hfill%
\includegraphics[page=1]{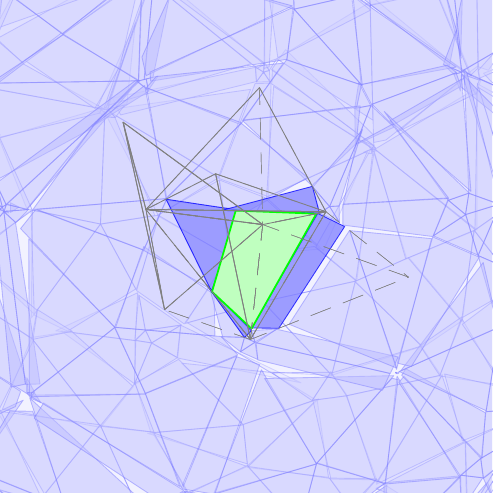}%
\hfill%
\includegraphics[page=2]{L06_lowres_stencil_ini}%
\hfill%
\includegraphics[page=3]{L06_lowres_stencil_ini}%
\hfill\null%
\\[12pt]%
\null\hfill%
\includegraphics[page=1]{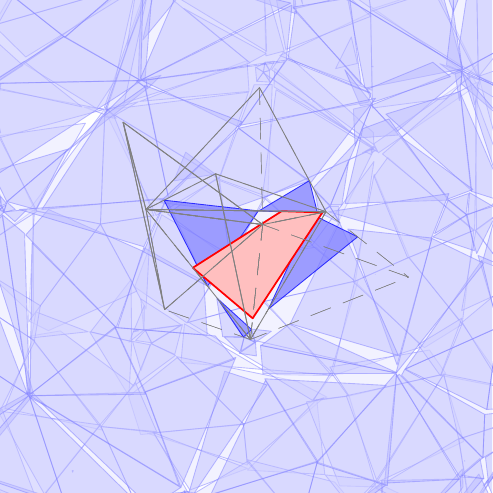}%
\hfill%
\includegraphics[page=2]{L06_lowres_stencil}%
\hfill%
\includegraphics[page=3]{L06_lowres_stencil}%
\hfill\null%
\caption{Local comparison of the reconstruction from the referential (green, center row) and orientation obtained from our approach (bottom row), respectively, for perturbed sphere ($\Liface=6$) for PLIC patches shown in top row.}%
\label{fig:l06_stencil_ini_comparison}%
\end{figure}

\begin{figure}[htbp]
\null\hfill%
\includegraphics[width=6cm]{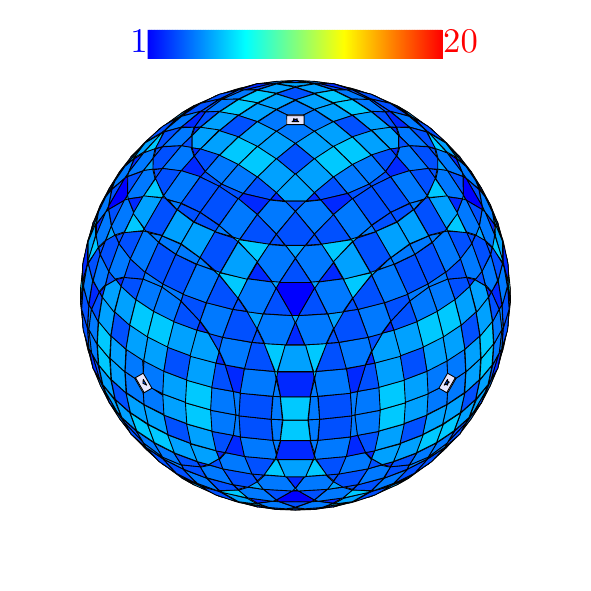}%
\hfill%
\includegraphics[width=6cm]{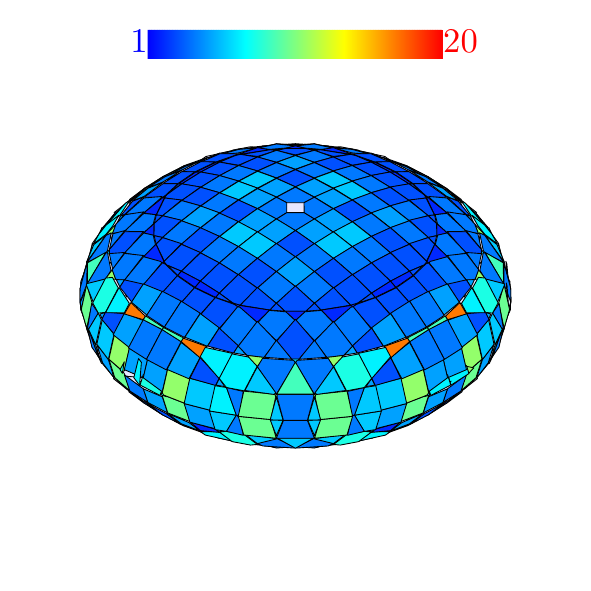}%
\hfill%
\includegraphics[width=6cm]{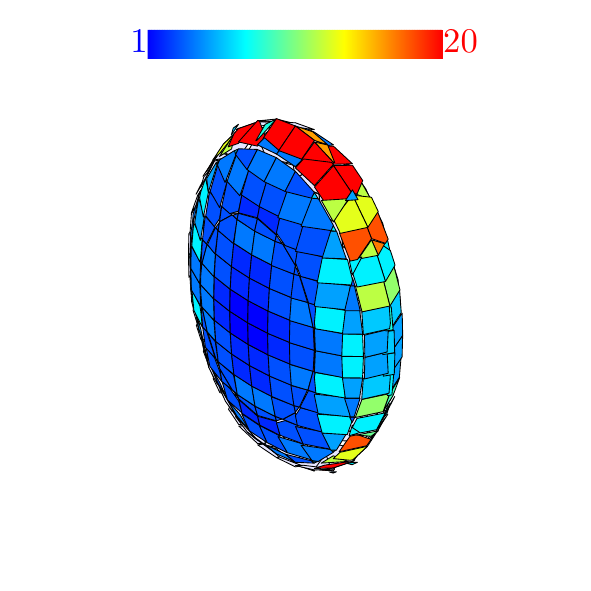}%
\hfill\null%
\caption{Comparison of number of iterations ($N_\mathrm{iter}^\mathrm{max}=20$) for sphere, oblate and prolate ellipsoid on cube mesh ($N=20$). As can be seen from the rightmost panel, regions with relatively large curvature require more iterations.}%
\label{fig:comparison_iterations_cube}%
\end{figure}

\begin{figure}[htbp]
\null\hfill%
\includegraphics[page=1]{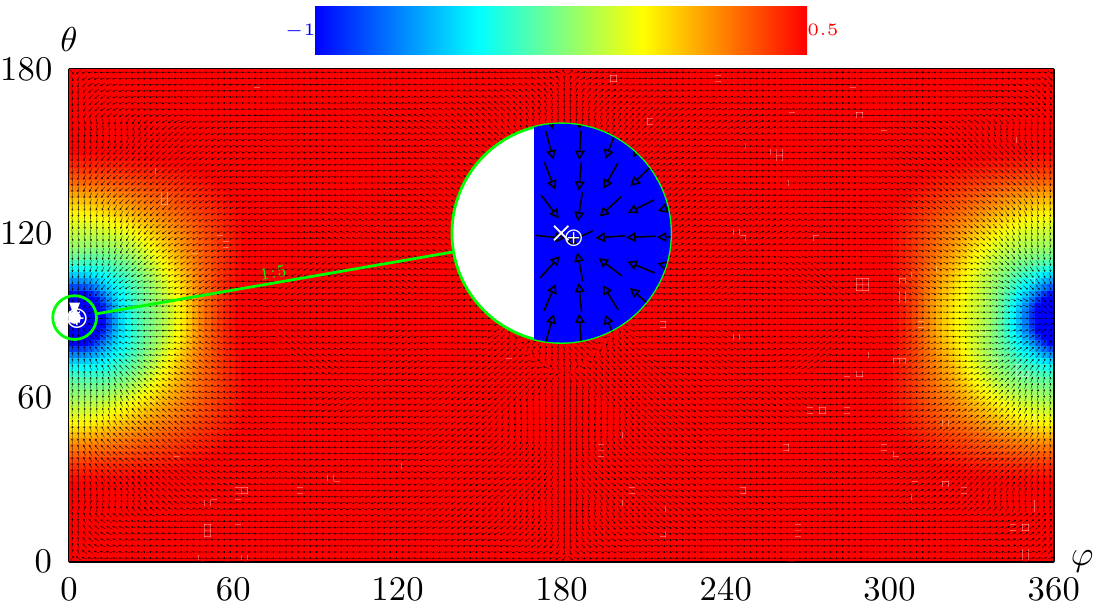}%
\includegraphics[page=1]{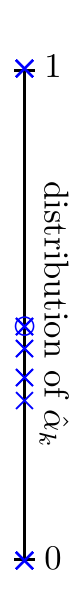}
\hfill%
\includegraphics[page=1]{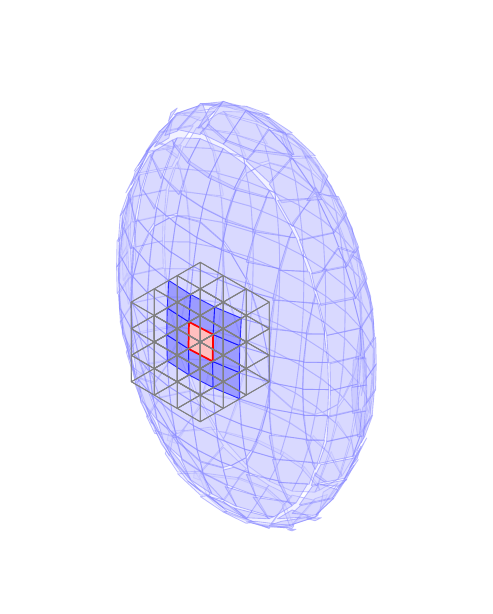}%
\hfill\null%
\\
\null\hfill%
\includegraphics[page=2]{prolate_cube_localanalysis_contour_magnified}%
\includegraphics[page=2]{prolate_cube_localanalysis_vof_distribution}
\hfill%
\includegraphics[page=2]{prolate_cube_localanalysis_illustration_local}%
\hfill\null%
\caption{Comparison of number of iterations for two stencils for a prolate ellipsoid on a cuboid mesh ($N=20$) with different local curvature resolution. The local error maps illustrate the effect of the insufficient resolution of the bottom configuration: the reduced variation of the error in the vicinity of the minimum slows down the minimization and inhibits convergence within the prescribed number of iterations ($N_\mathrm{iter}^\mathrm{max}=20$); cf.~\protect\reffig{comparison_iterations_cube}.}%
\label{fig:comparison_iterations_cube_closeup}%
\end{figure}

The present subsection gathers the error measures introduced in subsection~\ref{subsec:convergence_measures} for all combinations of meshes (tetrahedral with \textit{face} neighborhood / cuboidal with \textit{vertex} neighborhood) and hypersurfaces obtained with the proposed \textbf{F}ace \textbf{B}ased \textbf{N}ormal \textbf{R}econstruction (\textbf{FBNR} for short) algorithm. For comparison, we show the results obtained with existing algorithms: \textsc{Gauss-Green} (GG) as well as least squares estimation with (LSE) and without considering bulk cells (LSE\textsuperscript{*}), respectively. %
\begin{remark}[Comparison to literature]%
To the best of our knowledge, non-spherical hypersurfaces with non-trivial curvature have not been investigated in the literature to this date. %
For spheres, e.g.~\citet{JCP_2015_aina} and \citet{JCP_2019_aaes} provide analogous convergence studies of the normal alignment on different meshes. However, due to the conceptual difference in the employed error measures we can only compare to their work qualitatively. %
\end{remark}
Figures~\ref{fig:illustration_tet_algo_L06}~--~\reffigno{comparison_iterations_cube_closeup} provide further information. %
The main findings can be summarized as follows: %
\begin{enumerate}
%
%
\item In general, \reffigs{mesh_convergence_ellipsoids} and \reffigno{mesh_convergence_perturbed_sphere} show that the proposed algorithm exhibts second-order convergence with respect to both the symmetric volume difference $\average[\plicplane]{\Delta\polyvof}$ and normal alignment $\average[\plicplane]{\Delta\vn}$, for all configurations under consideration. This is in qualitative accordance with the findings of \citet[Fig.~7]{JCP_2019_aaes}, who report second-order convergence of an $L_\infty$-type error\footnote{While convergence in an $L_\infty$-type measure implies convergence for the $L_2$-pendant, we are convinced that the practical assessment of the problem at hand does not necessarily profit from such a strict measure.} of the normal alignment for spheres on unstructured meshes. %
For non-spherical hypersurfaces on cuboid meshes both error measures exhibit a decreasing slope for increasing spatial resolution. %
The rationale behind this phenomenon can be understood by recalling that the computation of the error measures $\average[\plicplane]{\Delta\vn}$ and $\average[\plicplane]{\Delta\polyvof}$ resorts to (i) the reconstruced normal $\recnormal{k}$, which, loosely speaking, approximates an average of the exact normal $\nS$ over the neighborhood $\neighborhood{k}$, as well as (ii) the normal $\ininormal{k}$ used for the initialization, which corresponds to a \textit{local} evaluation of the exact normal $\nS$ at, say, $\vx^\iface_{\mathrm{k}}$. %
For hypersurfaces of constant curvature, the averaging character of the reconstruction has virtually no effect, which directly corresponds to the uniform convergence observed for spheres in top row in \reffig{mesh_convergence_ellipsoids}. For non-spherical hypersurfaces, however, it implies that one compares a tangential approximation at a \textit{fixed point} $\vx^\iface_{\mathrm{k}}$ on the hypersurface to a planar approximation of $\iface$ based on a \textit{spatial average} over a neighborhood containing $\vx^\iface_{\mathrm{k}}$. If the neighborhood contains a non-symmetric segment of the hypersurface, which will almost always be the case for non-coinciding principal curvatures, one thus cannot expect the reconstructed normal to approximate the referential one, i.e.~$\ininormal{k}$. In other words, convergence can only be observed in a specific range of resolutions, depending on the mesh, the neighborhood type and the class and parameters of the hypersurface. %
Qualitatively, our investigation shows that the threshold solution, above which convergence degrades, strictly increases with increasing variation in the maximum difference of the principal curvature. %
While $\average[\plicplane]{\Delta\vn}$ and $\average[\plicplane]{\Delta\polyvof}$ constitute reliable error measures nonetheless, this raises the question whether the cell-based analysis of errors obtained from minimization over a neighborhood constitutes a suitable measure for convergence in general. However, adressing this issue in detail poses a highly non-trivial task in itself and goes beyond the scope of the present work.%
%
%
\item Figures~\ref{fig:illustration_tet_algo_L06} and \reffigno{illustration_cube_algo_L06} examplarily illustrate the PLIC reconstruction obtained by the proposed algorithm for a perturbed sphere ($\Liface=6$). %
\refFig{illustration_cube_algo_L06} visually supports the intuition of the convergence plots for cuboid meshes shown in \reffig{mesh_convergence_perturbed_sphere}: due to the well-behaved character of the minimization configuration (comparatively large spatial extension of the neighborhood as well as highly accurate volume fractions), the existing algorithms perform perfectly well. This can be expected since, in some sense, they are tailored to configurations of this type. On the other hand, increasing spatial resolution does neither imply qualitative nor quantitative improvement. %
As can be seen from \reffig{illustration_tet_algo_L06}, the contrary holds for tetrahedron meshes, where the gain of the proposed algorithm shows its full potential: while the convergence of the existing reconstruction algorithms reduced to first order, the proposed method pertains the second-order convergence. %
Furthermore, note that in \reffig{illustration_tet_algo_L06} the apparent qualitative similarity between the proposed algorithm (first row) and LSE\textsuperscript{*} (third row) degrades on closer inspection, showing noisy orientation for the LSE\textsuperscript{*} which are virtually independent of the resolution. This directly reflects to the respective convergence behavior in \reffig{mesh_convergence_perturbed_sphere}. %
%
%
\item \refFig{illustration_tet_algo_L06} shows that the quality of the existing reconstruction algorithms on tetrahedral meshes is insufficient, in the sense that no meaningful geometric transport of volume fraction can be performed based on their PLIC patches. Note that, even for the favorable case of highly accurate volume fraction data, increasing the spatial resolution does no alleviate the noisiness of the hypersurface approximation. %
%
%
\item Recall from subsection~\ref{subsubsec:initial_iteration} that the initial iterations for our approach correspond to LSE\textsuperscript{*}. From comparing the second and third row in \reffigs{illustration_tet_algo_L06} and \reffigno{illustration_tet_algo_L06} it is evident that the conjecture formulated for cuboidal meshes (i.e.~the least squares estimation profits from omitting the bulk cells with increasing resolution; cf.~\reffig{lse_reconstruction_weight_comparison}), analogously applies to tetrahedral meshes. %
%
%
\item For a prototypical stencil, \reffig{l06_stencil_ini_comparison} compares the PLIC patch reconstructed by the proposed algorithm to the referential orientation $\ininormal{k}$. While the algorithm converged in all cases shown, the volume-based minimization, which takes into account the neighborhood of the cell $\vfcell_k$, cannot reproduce the initial normal. %
%
%
\item For cuboid meshes (blue), the results shown in \reffig{mesh_convergence_ellipsoids} and \reffig{mesh_convergence_perturbed_sphere} exhibit a decrease in the rate of convergence with increasing maximum curvature and order $\Liface$, respectively, for all algorithms including the proposed one. This saturation-type behaviour can be explained by the fact that, locally within the neighborhood, the interface becomes more planar with increasing resolution, implying that the gain of the proposed algorithm diminishes. Note that increasing the maximum curvature (e.g.~by increasing $\Liface$) shifts the onset of the saturation to lower resolutions. E.g., for $\Liface=3$, GG shows virtually constant errors, implying that the spatial resolution is sufficient to reach the smallest possible error level of around $\average[\plicplane]{\Delta\vn}\approx\num{e-2}$. For $\Liface=9$, the error produced by GG decreases with increasing spatial resolution, where an onset of the aforementioned saturation occurs for $\sqrt{\Niface}\approx\num{e2}$. These findings apply analogously to LSE and LSE\textsuperscript{*}. %
%
%
\item As can be seen from \reffig{comparison_iterations_cube}, the number of iterations crucially depends on the volume fraction data in the stencil. \refFig{comparison_iterations_cube_closeup} illustrates the minimization process for the prolate ellipsoid on a cuboidal mesh ($N=20$) in the local error map, where the two selected cells can be considered prototypical: the first neighborhood contains an almost planar segment of the hypersurface $\iface$ (i.e.~with relatively small curvature), for which (i) the initial condition obtained form LSE\textsuperscript{*} already provides a very good approximation, (ii) the error steeply decreases in the vicinity of the minimum, which (iii) lies close to the referential orientation ($\oplus$). Contrary, the second stencil covers a segment of $\iface$ with comparatively large curvature. Also, the minimum exhibits an increased distance to the referential orientation ($\oplus$). The identical color legends highlight that the associated error exhibits far less variation in the second stencil, implying that the minimization produces a larger number of smaller steps. In the present case, convergence could not be obtained within the prescribed maximum number of iterations ($N_\mathrm{iter}^\mathrm{max}=20$). Nonetheless, the final iteration $\vp^{N_\mathrm{iter}^\mathrm{max}}$ is located in the vicinity of the sought minimum. %
\end{enumerate}

%
%

%
%


%% file: 05_04_local_analysis.tex
%
%
\subsection{Local attraction characteristics}%
\begin{figure}[htpb]%
\null\hfill%
\includegraphics[page=2]{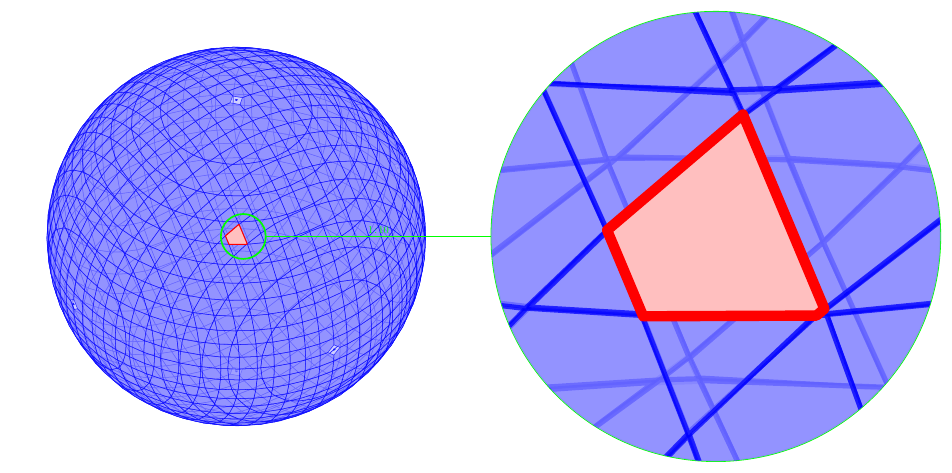}%
\hfill\null%
\\%
\null\hfill%
\includegraphics[page=2]{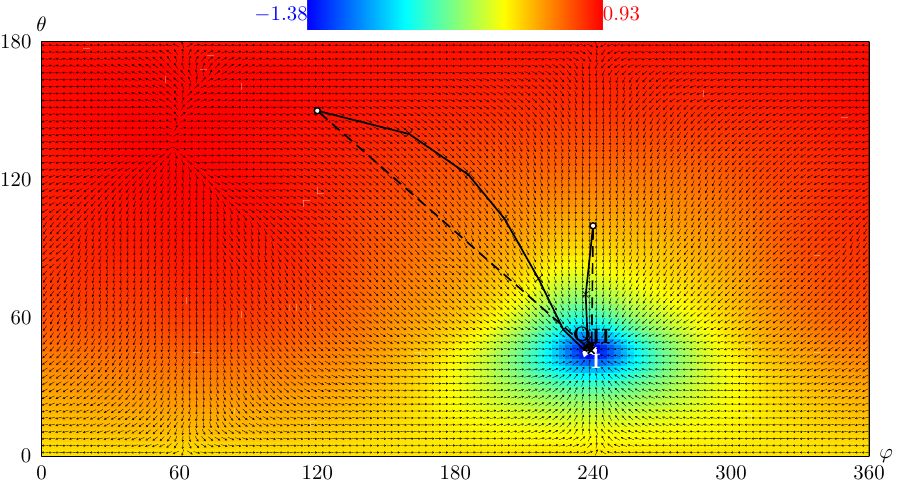}%
\includegraphics[page=2]{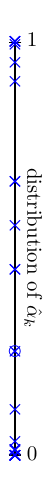}%
\hfill\null%
\\%
\null\hfill%
\includegraphics[page=8]{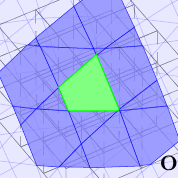}%
\hfill%
\includegraphics[page=5]{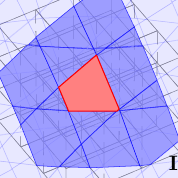}%
\hfill%
\includegraphics[page=6]{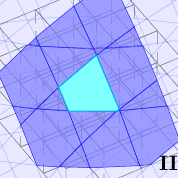}%
\hfill\null%
\caption{Reconstructed sphere on a cuboidal mesh ($N=30$; cf.~\reftab{cube_mesh_setup}) with exemplary \textit{compliant} patch, along with the local error map (cf.~\refnote{local_error_map}) exhibiting a \textbf{unique} minimum. The referential orientation (\textbf{O}) along with the iterations ($\times$) and final orientations obtained from the present approach are illustrated for different initial conditions (\textbf{I}: LSE\textsuperscript{*}, \textbf{II}: artificial). Obviously, due to the uniqueness of the attractor, the converged iterations coincide.}%
\label{fig:cuboidmesh_local_attractor_unique}%
\end{figure}

The numerical reason behind the existence of \textit{outliers}, as shown in \reffig{illustration_tet_algo_L06}, can be found by inspecting the characteristics of the error functional from \refeqn{error_functional_unconstrained_gradient_hessian} for the respective cells, which we shall discuss for two exemplary but conceptually prototypical configurations. %
%
%
\paragraph{Unique attractor}%
\refFig{cuboidmesh_local_attractor_unique} depicts a configuration for which the error functional in \refeqn{error_functional_unconstrained_gradient_hessian} exhibts a unique minimum. Hence, irrespective of the initial condition $\vp^0$, the minimization produces the desired (henceforth also referred to as \textit{compliant}) orientation; however, the number of iterations may differ. Apparently, configurations of this class emerge for (i) sufficiently resolved hypersurfaces on (ii) cuboid meshes with vertex neighborhood within which (iii) the volume fraction data $\polyvofdata_k$ admit a rather uniform distribution across $[0,1]$ as well as (iv) a value away from zero or one in the center cell. %

\begin{figure}[htpb]%
\null\hfill%
\includegraphics[page=1]{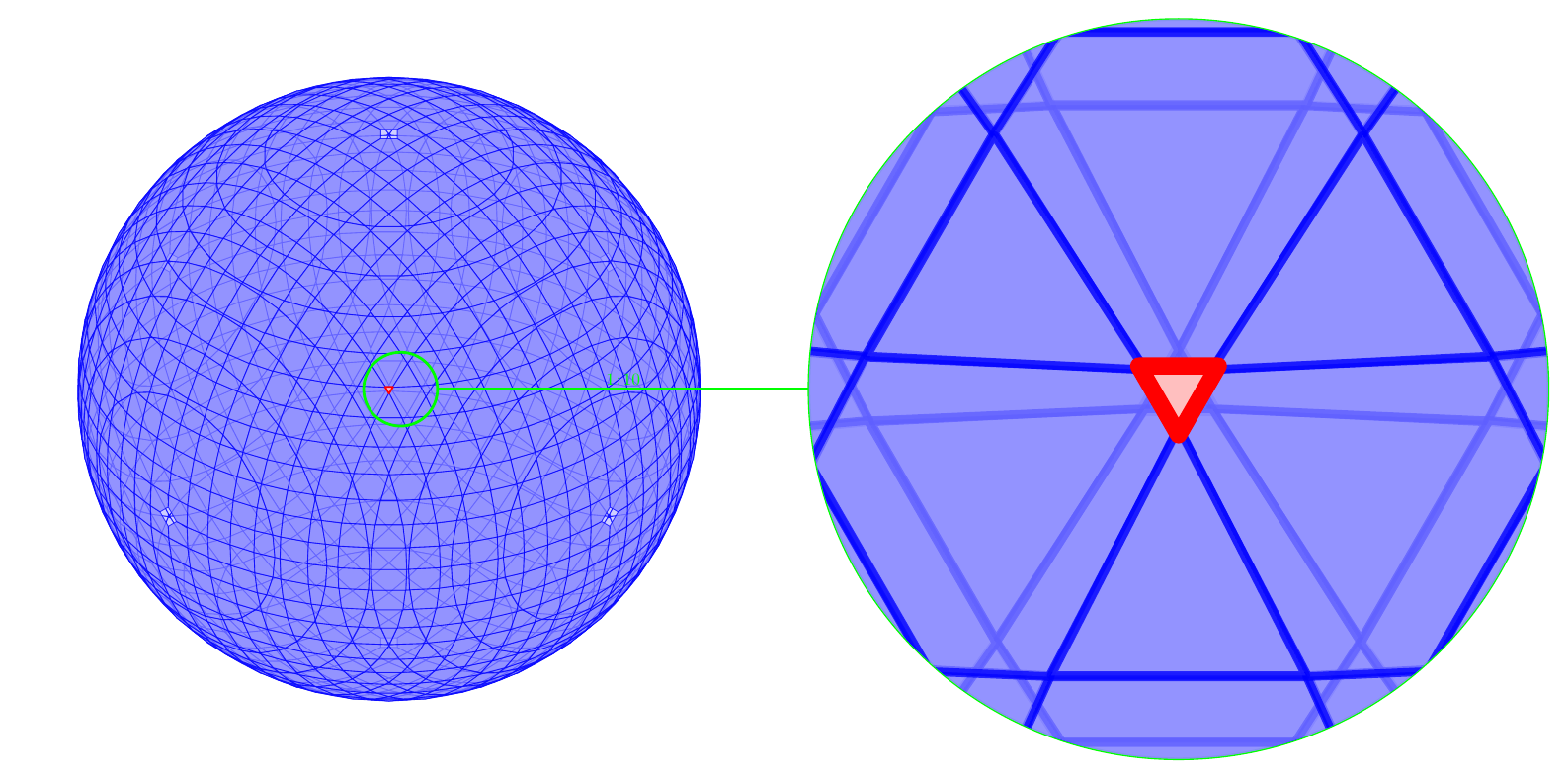}%
\hfill\null%
\\%
\null\hfill%
\includegraphics[page=1]{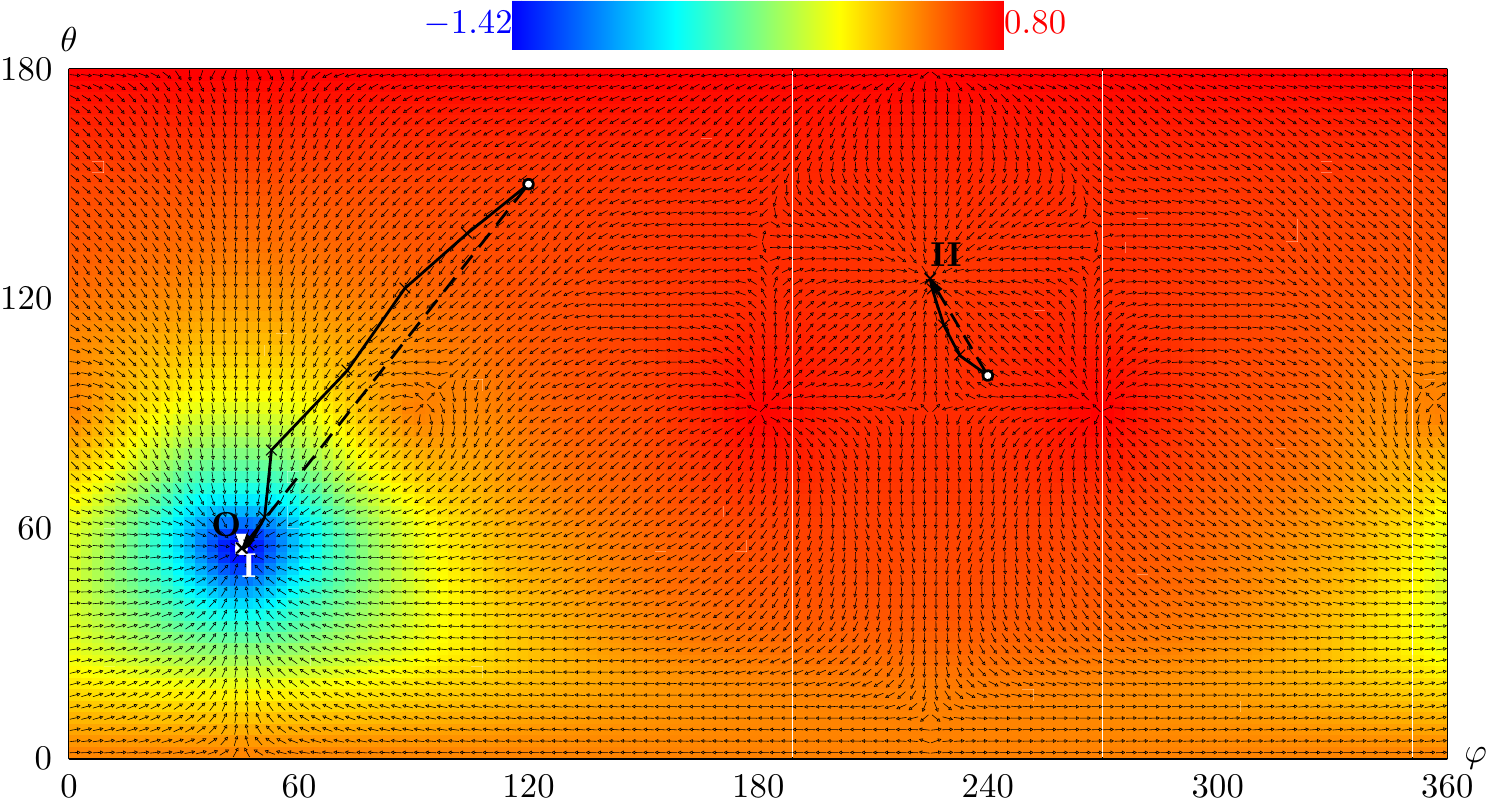}%
\includegraphics[page=1]{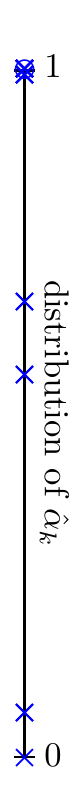}%
\hfill\null%
\\%
\null\hfill%
\includegraphics[page=4]{attractortype_cube_1-dot-0E+00_stencil}%
\hfill%
\includegraphics[page=1]{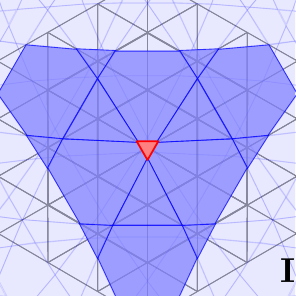}%
\hfill%
\includegraphics[page=2]{attractortype_cube_1-dot-0E+00_stencil}%
\hfill\null%
\caption{Reconstructed sphere on a cuboidal mesh ($N=30$; cf.~\reftab{cube_mesh_setup}) with exemplary \textit{compliant} patch, along with the local error map (cf.~\refnote{local_error_map}) exhibiting a pair of \textbf{diametrical} minima. The referential orientation (\textbf{O}) along with the iterations ($\times$) and final orientations obtained from the present approach are illustrated for different initial conditions (\textbf{I}: LSE\textsuperscript{*}, \textbf{II}: artificial).}%
\label{fig:cuboidmesh_local_attractor_opposing}%
\end{figure}

\begin{figure}[htbp]
\null\hfill%
\includegraphics[page=3]{error_square_illustration}%
\hfill%
\includegraphics[page=1]{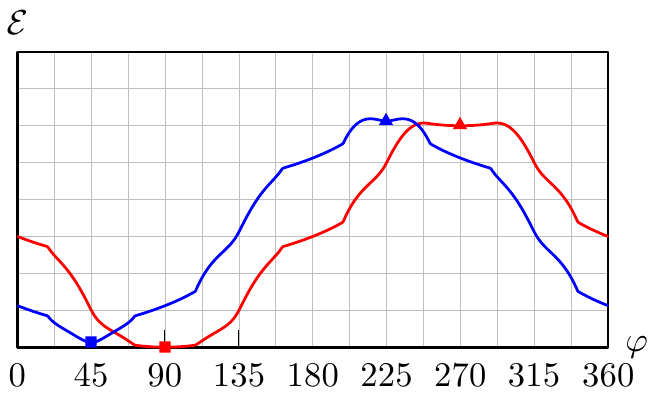}%
\hfill%
\includegraphics[page=2]{error_square_illustration}%
\hfill\null%
\caption{Illustration of the error functional from \protect\refeqn{error_functional_unconstrained_gradient_hessian} for $3\times3$ square stencil in two spatial dimensions with different sets of volume fraction data $\protect\set{\polyvofdata_k}$ with \textit{compliant} ($\blacksquare$) and diametrical attractor ($\blacktriangle$, shifted by $\pi$).}%
\label{fig:3x3_stencil_error_illustration}%
\end{figure}

\paragraph{Diametrical attractor}%
\refFig{cuboidmesh_local_attractor_opposing} illustrates a peculiarity of cuboidal meshes: if the volume fraction of the center cell $\polyvofdata_0$ (marked by $\otimes$ in the distribution plot) admits a value close to zero or one (ceteris paribus, cf.~(i)-(iii) above), the attractor associated to the \textit{compliant} orientation may become non-unique due to the emergence of a diametrical attractor. 
\begin{remark}[Quasi-diametrality]%
Technically, two points $\vx$ and $\vy$ on the boundary of the unit sphere are diametrical iff $\iprod{\vx}{\vy}=-1$. Here, however, we employ the term diametrical to qualitatively describe the orientations, i.e.~two orientations are considered diametrical if $\iprod{\vx}{\vy}\approx-1$. %
\end{remark}
In a geometrical sense, this phenomenon is closely related to the orientation ambiguity for the least-squares estimation outlined in subsection~\ref{subsubsec:initial_iteration}. %
The two-dimensional analogon ($3\times3$ neighborhood of squares) shown in \reffig{3x3_stencil_error_illustration} nicely conveys the rationale behind this phenomenen: the diametrical minima (marked by $\blacktriangle$) correspond to a plane with inverted orientation, which locally also minimizes the error, however inducing a larger value. In the center layer of the stencil in the rightmost panel, the volume fractions and their gradient with respect to $\varphi$ coincide for the two orientations (\textcolor{red}{$\blacktriangle$} and \textcolor{red}{$\blacksquare$}). Recalling that only intersected cells contribute to the gradient of the error functional in \refeqn{error_functional_unconstrained_gradient_hessian} implies the existence of a diametrical minimum. %
Moreover, note that in two spatial dimensions, the diametrical minimum exists also for center volume fractions away from 0 and 1. %

\begin{figure}[htpb]%
\null\hfill%
\includegraphics{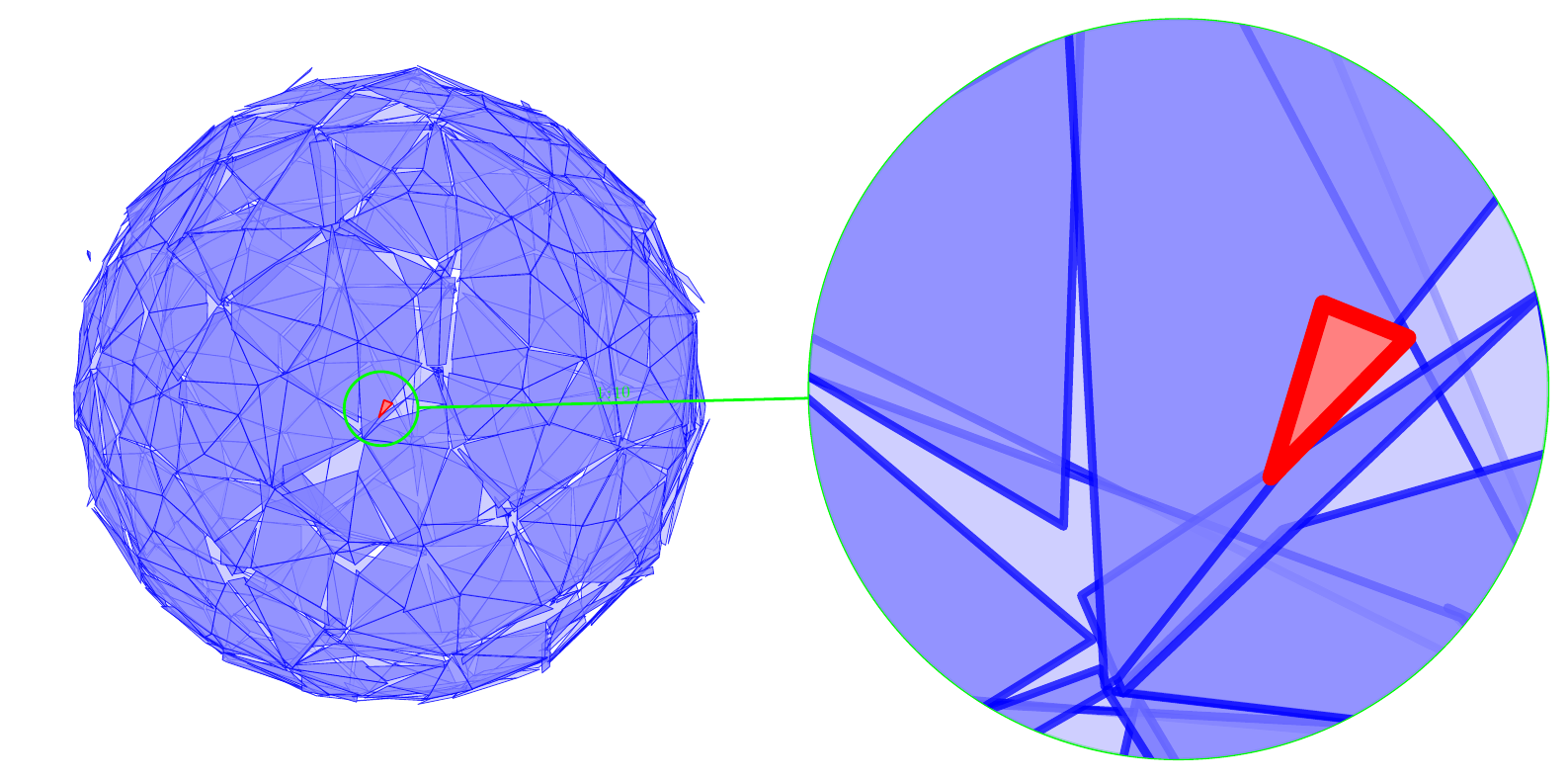}%
\hfill\null%
\\%
\null\hfill%
\includegraphics[page=1]{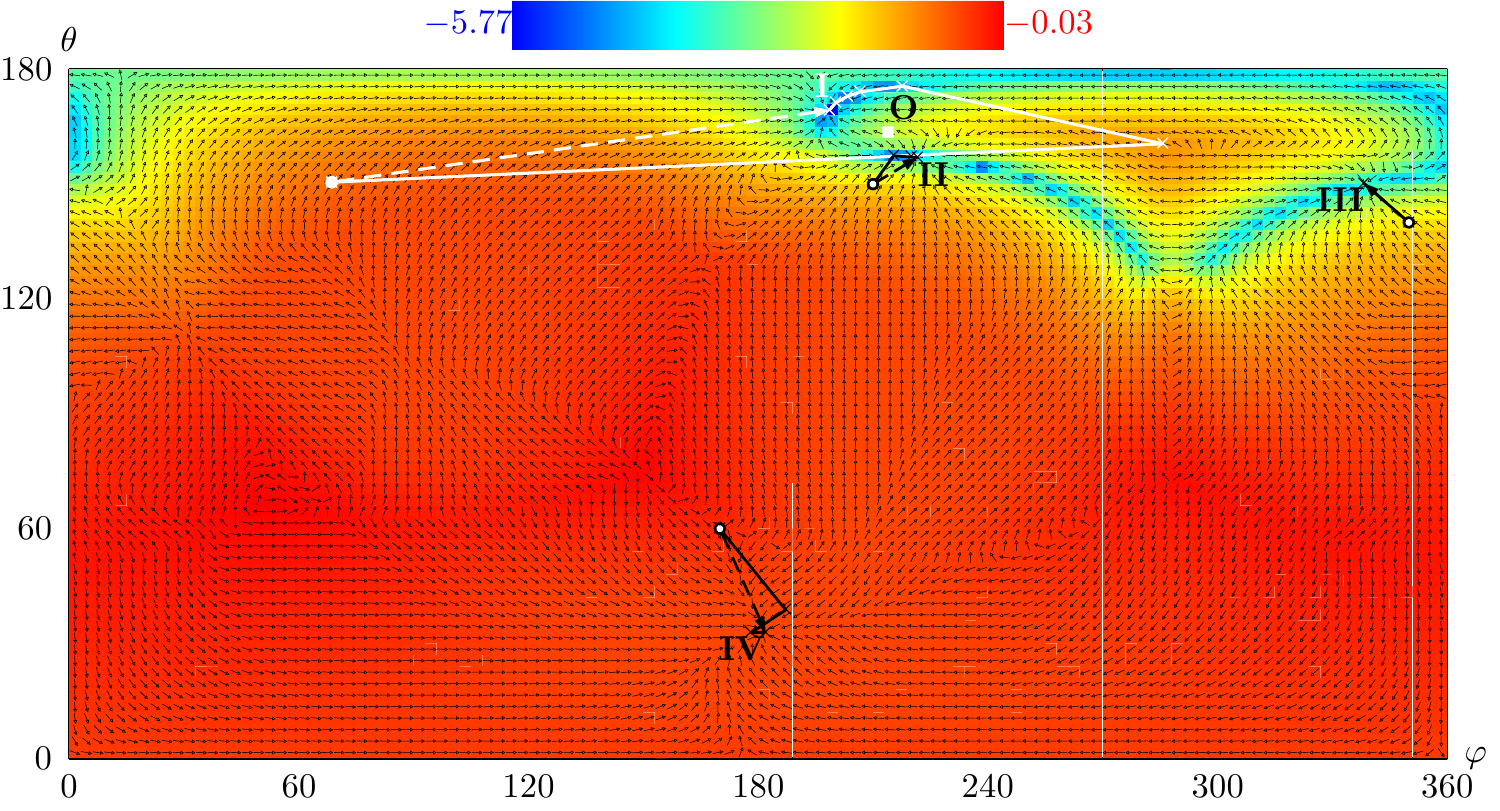}%
\hfill%
\includegraphics[page=1]{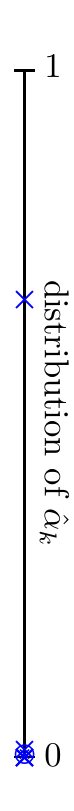}%
\hfill\null%
\\%
\null\hfill%
\includegraphics[page=5]{attractortype_tet_1-dot-0E+00_stencil}%
\hfill%
\includegraphics[page=1]{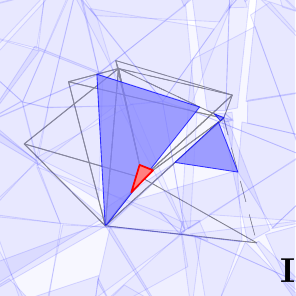}%
\hfill%
\includegraphics[page=2]{attractortype_tet_1-dot-0E+00_stencil}%
\hfill%
\includegraphics[page=3]{attractortype_tet_1-dot-0E+00_stencil}%
\hfill%
\includegraphics[page=4]{attractortype_tet_1-dot-0E+00_stencil}%
\hfill\null%
\caption{Reconstructed sphere on tetrahedral mesh ($N=10$; cf.~\reftab{tetrahedron_mesh_setup}) with exemplary \textit{outlier} patch, along with the local error map (cf.~\refnote{local_error_map}) exhibiting four \textbf{non-unique} minima. The referential orientation (\textbf{O}) along with the iterations ($\times$) and final orientations obtained from the present approach are illustrated for different initial conditions (\textbf{I}: LSE\textsuperscript{*}, \textbf{II}-\textbf{IV}: artificial).}%
\label{fig:tetmesh_local_attractor}%
\end{figure}

%
%
\paragraph{Non-unique attractors}%
Moving to tetrahedral meshes with face neighborhood entirely degrades the well-behaved characteristics obtained for cuboid meshes. As a result, the error functional and, hence, its minimization become far more demanding. In particular, general meshes will not admit unique attractors. %
\refFig{tetmesh_local_attractor} illustrates a prototypical configuration with non-unique attractors. %
For these general cases, we summarize our findings as follows: %
\begin{enumerate}
\item The error functional admits multiple local minima, whose location crucially depends, among others, on the volume fraction data, neighborhood type and spatial extension of the stencil. %
\item The local minima may be located in close vicinity to each other, highlighting the sensitivity of the minimization process with respect to the choice of the initial value $\vp^0$ for the iteration. As a consequence, comparison of cuboidal and tetrahedral meshes visually suggests a slightly increased noise of the overall reconstruction for the latter; cf.~\reffigs{illustration_tet_algo_L06} and \reffigno{illustration_cube_algo_L06}. %
\item The referential orientation (marked by \textbf{O} in \refFig{tetmesh_local_attractor}) does not necessarily coincide with one of the minima. In fact, as can be seen from the location of \textbf{I} and \textbf{II} in \reffig{tetmesh_local_attractor}, respectively, one cannot assign it to one of the attractors even qualitatively. This further highlights the challenging character of the minimization.  
\end{enumerate}
In fact, this motivates the choice of the bulk weights $\set{\mu_k}$ outlined in the next subsection~\ref{subsec:influence_weights}. %

%% file: 05_05_influence_weights.tex
\subsection{Influence of the weights $\mu_k$}\label{subsec:influence_weights}%
\begin{figure}[htpb]%
\null\hfill%
\includegraphics{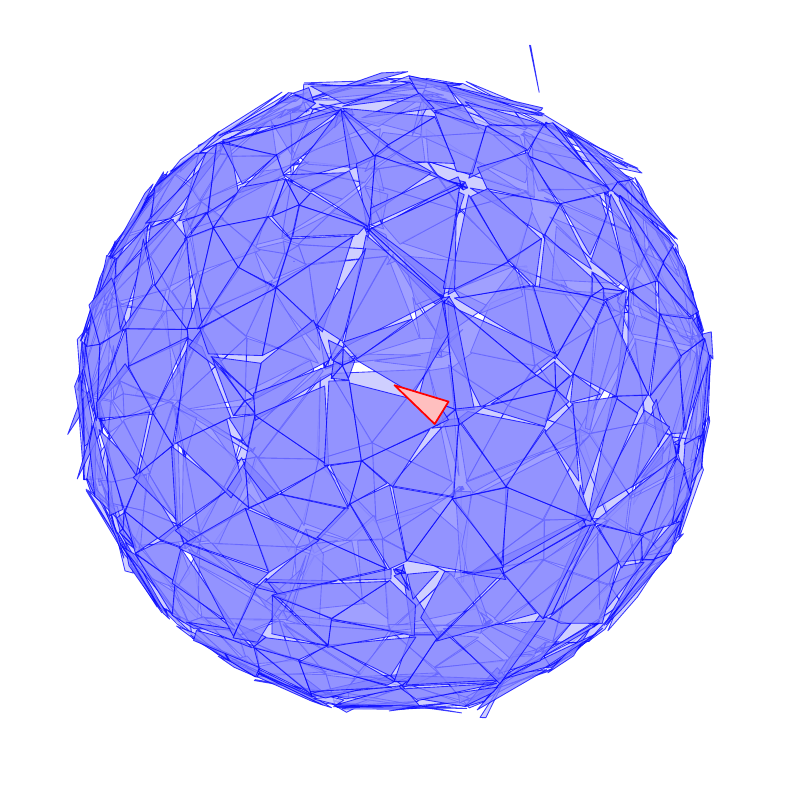}%
\hfill%
\includegraphics{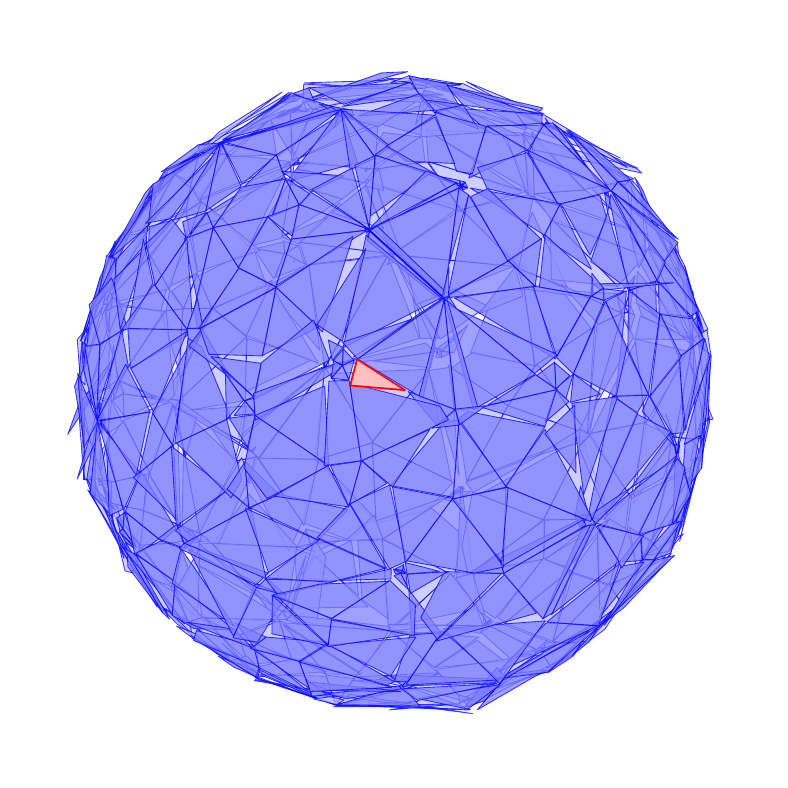}
\hfill\null%
\\%
\null\hfill%
\includegraphics{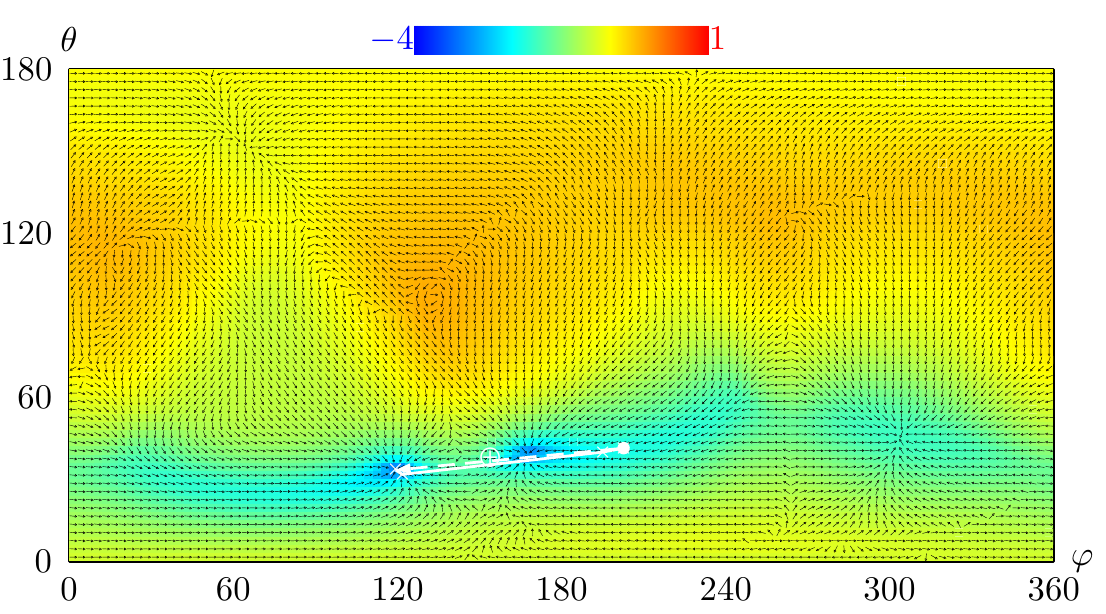}%
\includegraphics{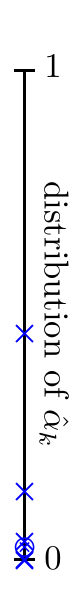}%
\hfill%
\includegraphics{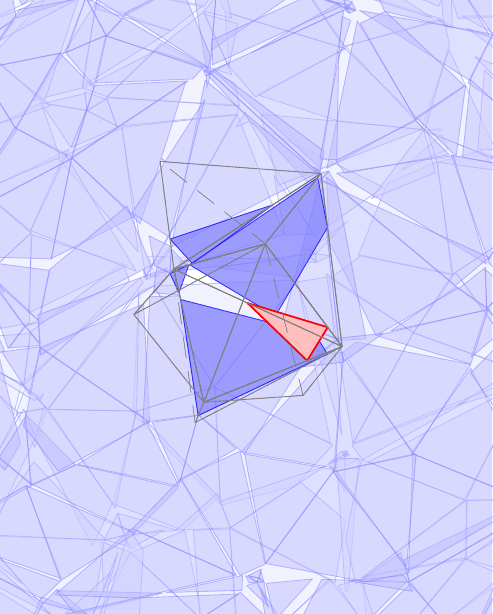}
\hfill\null%
\\%
\null\hfill%
\includegraphics{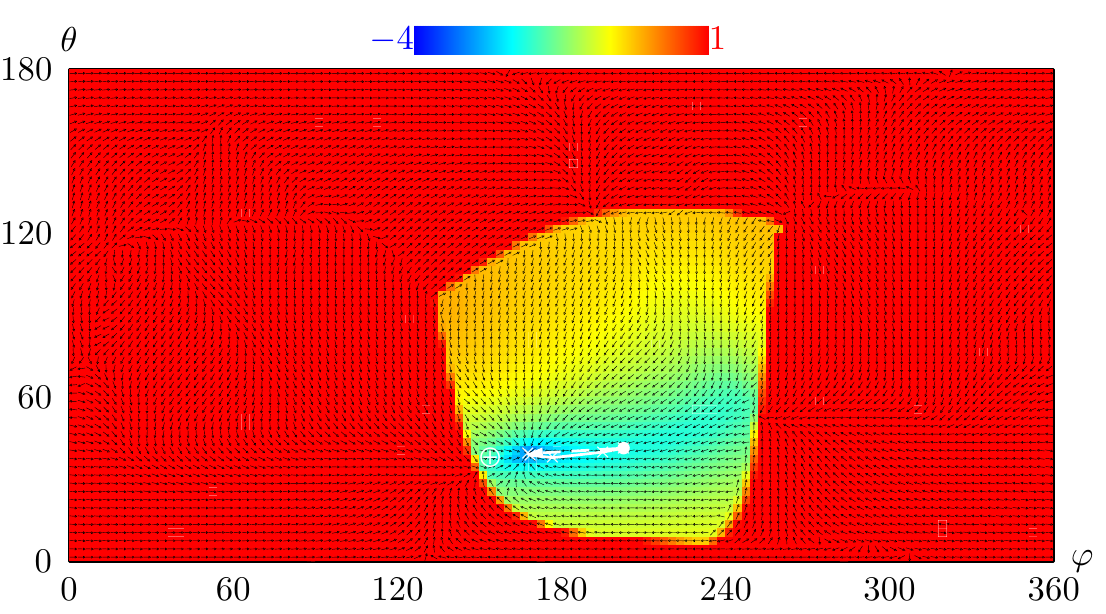}%
\includegraphics{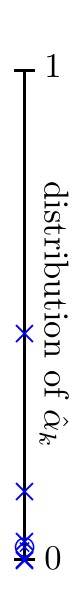}%
\hfill%
\includegraphics{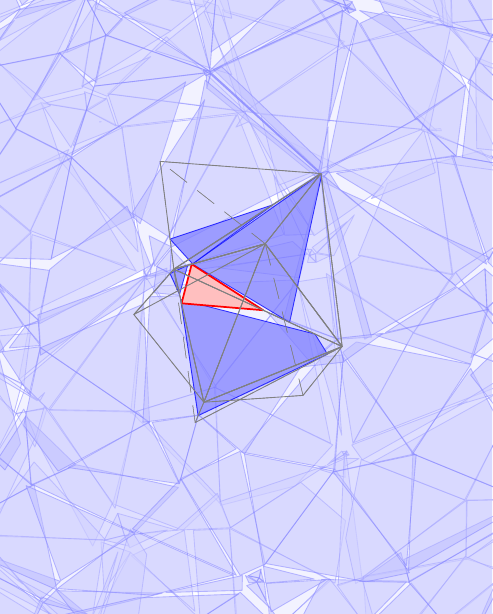}
\hfill\null%
\caption{Comparison of local error map and PLIC reconstruction for a sphere on tetrahedron mesh ($N=10$) for a stencil with bulk weights $\mu_k=1$ (center row/left panel in top row) and $\mu_k=\num{e9}$ (bottom row/right panel in top row). As intended, the large numerical values of $\mu_k$ for data-wise bulk cells remove the \textit{non-compliant} attractor at approx.~$(120^\circ,40^\circ)$, towards which the minimization with $\mu_k=1$ converges (center row).}%
\label{fig:comparion_influence_weight}%
\end{figure}

As outlined in subsection~\ref{subsubsec:choice_weights}, choosing large numerical values for the weights $\set{\mu_k}$ associated to the bulk cells in the error functional corresponds to imposing an additional constraint for the topological intersection status of the neighborhood; cf.~\reffig{intersection_status}. %
\refFig{comparion_influence_weight} compares the effect of the choice of weights for a prototypical configuration: %
for small neighborhoods, choosing large bulk weights $\set{\mu_k}$ removes some of the undesired local minima of the error, which drives the minimization procedure to converge to the \textit{compliant} orientation. On spatially extensive neighborhoods, such as the \textit{vertex} neighborhood for cuboid meshes, the opposite effect can be observed. In general, the bulk weights $\set{\mu_k}$ should be inversely correlated to the spatial extension of the neighborhood. %

%% file: 06_summary.tex
\section{Summary \& conclusion}\label{sec:summary}%
This work introduces an iterative minimization algorithm to approximate a cell-based discrete normal field from volume fractions in a given neighborhood on an unstructured mesh composed of arbitrary polyhedral cells. Beyond a comprehensive formulation of the numerical scheme, the regularity of the underlying error functional has been thoroughly discussed in detail. The proposed method has been assessed in an extensive set of numerical experiments for spherical and non-spherical hypersurfaces on cuboidal and tetrahedral meshes. Beyond the usual convergence studies, this investigation provides deep insights into the minimization process. %
Altogether, the following conclusions can be drawn: %
\begin{enumerate}%
\item For the combinations of hypersurfaces and resolutions considered, second-order convergence with spatial resolution for both, $\average[\plicplane]{\Delta\polyvof}$ and $\average[\plicplane]{\Delta\vn}$, was obtained. %
\item For cuboid meshes with vertex neighborhood, the proposed algorithm performs slighly better than existing standard approaches, such as least-square error and \textsc{Gauss-Green}. Moreover, no outlying PLIC patches were observed. %
\item For tetrahedral meshes, the performance gain of the proposed algorithm in comparison to existing methods was shown to be substantial: while standard approaches do not show mesh convergence, our approach produces second-order convergence. The small number of outliers, which occur for the challenging face-neighborhood configurations, can be removed by assigning weights to bulk cells whose numerical value should be chosen inversely to the spatial extension of the neighborhood. %
\item We obtain both quantitively (in terms of normal alignment and symmetric volume difference) and qualitatively (in terms of visual patch alignment) convincing discrete hypersurface reconstruction. %
\end{enumerate}%

%% file: 99_appendix.tex
\begin{appendix}
\section{Implicit bracketing}\label{app:implicit_bracketing}%
\begin{figure}[htbp]%
\null\hfill%
\includegraphics[page=1]{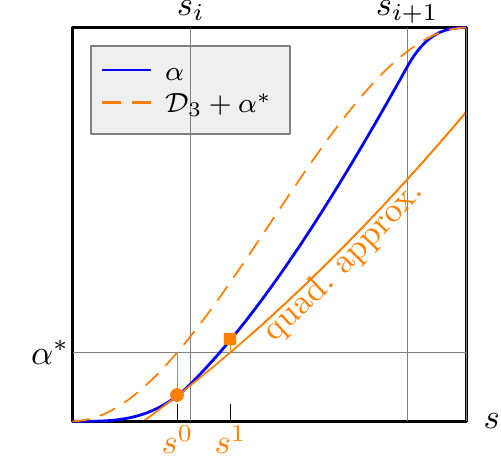}
\hfill%
\includegraphics[page=2]{implicit_bracketing}
\hfill%
\includegraphics[page=3]{implicit_bracketing}
\hfill\null%
\caption{Implicit bracketing and locally quadratic approximation of volume fraction $\polyvof\fof{\signdist}$.}%
\label{fig:implicit_bracketing_illustration}%
\end{figure}

\begin{center}
\textit{This section heavily draws from a previous work of the authors; \cite[sections~2.3 and 2.4]{JCP_2021_fbip}.}%
\end{center}
Within a bracket $\mathcal{B}_i\defeq[\signdist_i,\signdist_{i+1}]$, the volume fraction $\polyvof\fof{\signdist}$ is an increasing cubic polynomial, denoted by $\mathcal{S}_i\fof{\signdist}$. For any given $\signdist^n\in\mathcal{B}_i$, the polynomial reads %
\begin{align}
\mathcal{S}_i\fof{z;\signdist^n}=%
\frac{\polyvof^{\prime\prime\prime}\fof{\signdist^n}}{6}\brackets*{z-\signdist^n}^3+%
\frac{\polyvof^{\prime\prime}\fof{\signdist^n}}{2}\brackets*{z-\signdist^n}^2+%
\polyvof^{\prime}\fof{\signdist^n}\brackets*{z-\signdist^n}+%
\polyvof\fof{\signdist^n}.\label{eqn:third_order_approximation}%
\end{align}%
Hence, the truncation of the polyhedron $\polyhedron*$ at any $\signdist^n$ (implicitly) provides the full information of $\polyvof\fof{\signdist}$ within the containing bracket $\mathcal{B}_i$. Exploiting that $\polyvof_i=\mathcal{S}_i\fof{\signdist_i;\signdist^n}$ and $\polyvof_{i+1}=\mathcal{S}_i\fof{\signdist_{i+1};\signdist^n}$ suggests the following strategy: %
\begin{enumerate}%
\item if the current iteration $\signdist^n$ is \textit{not} located in the target bracket $\refbracket$ ($\refvof<\polyvof_i$ or $\polyvof_{i+1}<\refvof$), the next iteration is obtained from locally quadratic approximation (\reffig{implicit_bracketing_illustration}, left). %
\item if the current iteration $\signdist^n$ is located in the target bracket $\refbracket$ ($\polyvof_i\leq\refvof\leq\polyvof_{i+1}$), the sought $\signdistref$ corresponds to the root of $\mathcal{S}_i-\refvof$, requiring no further truncation (\reffig{implicit_bracketing_illustration}, center). %
\end{enumerate}%
The initial guess is obtained from a global cubic spline interpolation via %
\begin{align}
s^0\fof{\refvof}\defeq\projvertposmin+\brackets{\projvertposmax-\projvertposmin}\brackets*{\frac{1}{2}-\cos\brackets*{\frac{\arccos\fof{2\refvof-1}-2\pi}{3}}}%
\end{align}
with $\projvertposmin=\min\fof{\hat{\mathcal{S}}}$ and $\projvertposmax=\max\fof{\hat{\mathcal{S}}}$. %
\refFig{implicit_bracketing_illustration} illustrates the components of the strategy outlined above: with $\signdist^0\not\in\mathcal{B}_i$ (left), a quadratic approximation yields $\signdist^{1}$, which lies within the target bracket $\refbracket=\mathcal{B}_i$ (center). %
The rightmost configuration already starts with $\signdist^0\in\refbracket$, such that the spline interpolation directly yields the sought position $\signdistref$, implying that only a single truncation is required. %
\section{Derivatives of the volume and area coefficients}\label{app:derivatives_vol_coeff}%
With %
\def\subone{a}%
\def\subtwo{b}%
\def\subthree{c}%
\begin{align}
\subone\defeq\iprod{\vn^\vfface_f}{\plicnormal},\quad%
\subtwo\defeq-\iprod{\vN_{k,m}}{\plicnormal}\quad\text{and}\quad%
\subthree\defeq\iprod*{\iprod{\vx^\vfface_{k,1}}{\vn^\vfface_f}\vn^\vfface_f-\xbase}{\plicnormal},%
\end{align}
for notational convenience, the coefficients in \refeqn{face_coefficients} become %
\begin{align}
\volcoefflin{f}=-\subone,\quad%
\areacoeffconst{f}{m}=\iprod{\vx^\vfface_{k,m}}{\vN_{k,m}}-\subthree\areacoefflin{f}{m}\quad\text{and}\quad%
\areacoefflin{f}{m}=\frac{\subtwo}{1-\subone^2}.%
\end{align}
With $i,j\in\set{\varphi,\theta}$, $\partial_if$ and $f_i$ both denoting $\frac{\partial f}{\partial p_i}$, the derivatives of the coefficients in \refeqs{volume_coefficients} and \refeqno{face_coefficients} read %
\begin{align}%
\partial_{i}\volcoeffconst{f}=0,\quad%
\quad\partial_{i}\volcoefflin{f}=-\subone_i,\quad%
\partial_{i}\areacoeffconst{f}{m}=-\brackets{\subthree_i\areacoefflin{f}{m}+\subthree\partial_{i}\areacoefflin{f}{m}},\quad%
\partial_{i}\areacoefflin{f}{m}=\frac{\subtwo_i}{1-\subone^2}+2\frac{\subtwo\subone\subone_i}{\brackets{1-\subone^2}^2}.%
\end{align}%
\section{Mesh generation with \texttt{gmsh}}\label{app:gmsh}%
The tetrahedral meshes used in section~\ref{sec:numerical_results} were generated with \href{https://gmsh.info}{\texttt{gmsh} 4.7.1}. For $h=\frac{1}{N}$, cf.~\reftab{tetrahedron_mesh_setup}, the file \texttt{mesh.geo} gathers the relevant information. %
\lstset{language=C++,commentstyle=\color{blue!50},frame=single,caption={Example geometry file ($N=20$)}}%
\begin{lstlisting}
    // add points
    Point(1) = {0, 0, 0, h};Point(2) = {1, 0, 0, h};
    Point(3) = {1, 1, 0, h};Point(4) = {0, 1, 0, h};
    Point(5) = {0, 0, 1, h};Point(6) = {1, 0, 1, h};
    Point(7) = {1, 1, 1, h};Point(8) = {0, 1, 1, h};
    // add lines
    Line(1) = {1, 2};Line(2) = {2, 3};Line(3) = {3, 4};Line(4) = {4, 1};
    Line(5) = {1, 5};Line(6) = {2, 6};Line(7) = {3, 7};Line(8) = {4, 8};
    Line(9) = {5, 6};Line(10) = {6, 7};Line(11) = {7, 8};Line(12) = {8, 5};
    // loops and surfaces
    Line Loop(13) = {5,-12,-8,4};Line Loop(14) = {2,7,-10,-6};
    Line Loop(15) = {-4,-3,-2,-1};Line Loop(16) = {9,10,11,12};
    Line Loop(17) = {1,6,-9,-5};Line Loop(18) = {3,8,-11,-7};
    Plane Surface(1) = {13};Plane Surface(2) = {14};
    Plane Surface(3) = {15};Plane Surface(4) = {16};
    Plane Surface(5) = {17};Plane Surface(6) = {18};
    Surface Loop(1) = {6, 3, 1, 5, 2, 4};
    Volume(1) = {1};
\end{lstlisting}
With the above geometry file, the mesh is generated by invoking %
\begin{center}
\verb+gmsh -refine -smooth 100 -optimize_netgen -save -3 -format vtk -o mesh.vtk mesh.geo+
\end{center}

\end{appendix}